\documentclass[a4paper,11pt]{article}
\usepackage[a4paper,text={16cm,22.5cm},centering]{geometry}
\usepackage{amsmath,amssymb,mathtools}
\usepackage{empheq}
\usepackage{mathrsfs,bbm,yhmath}
\usepackage[utf8]{inputenc}
\usepackage[TS1,T1]{fontenc}
\usepackage{lmodern}
\usepackage{microtype}
\usepackage[english]{babel}
\usepackage{tikz}

\usetikzlibrary{arrows} 
\usetikzlibrary{patterns} 
\usetikzlibrary{decorations.markings} 
\usetikzlibrary{decorations.pathmorphing}

\definecolor{vert}{rgb}{0.1,0.7,0.15}

\tikzset{hatch distance/.store in=\hatchdistance, hatch distance=10pt, hatch thickness/.store in=\hatchthickness, hatch thickness=2pt}

\makeatletter 
\pgfdeclarepatternformonly[\hatchdistance,\hatchthickness]{flexible hatch} {\pgfqpoint{-1pt}{-1pt}} {\pgfqpoint{\hatchdistance}{\hatchdistance}} {\pgfpoint{\hatchdistance-1pt}{\hatchdistance-1pt}} {\pgfsetcolor{\tikz@pattern@color} 
\pgfsetlinewidth{\hatchthickness} \pgfpathmoveto{\pgfqpoint{-1pt}{-1pt}} \pgfpathlineto{\pgfqpoint{\hatchdistance}{\hatchdistance}} \pgfusepath{stroke}} 
\makeatother

\rmfamily
\DeclareFontShape{T1}{lmr}{b}{sc}{<->ssub*cmr/bx/sc}{}
\DeclareFontShape{T1}{lmr}{bx}{sc}{<->ssub*cmr/bx/sc}{}

\DeclareMathOperator*{\esup}{ess~sup}
\DeclareMathOperator*{\argmin}{arg~min}

\newtheorem{defi}{Definition}[section]
\newtheorem{prop}[defi]{Proposition} 
\newtheorem{lemm}[defi]{Lemma} 
\newtheorem{theo}[defi]{Theorem}

\numberwithin{equation}{section}
\newcommand{\qed}{\hglue 0pt\hfill$\square$\par}

\title{Transition from Continuous to Jumping Solutions in 2D Quasi-static Elastic Contact Problems with Coulomb Friction:\\ the Mathematics Underlying the Onset of Brake Squeal}
\author{Patrick Ballard,\\[1.5ex]
{\small Institut Jean Le Rond d'Alembert,}\\ 
{\small Sorbonne Université, \textsc{cnrs}, Université de Paris}\\
{\small 4, place Jussieu, 75252 Paris Cedex 05, France}\\
{\small \texttt{patrick.ballard@dalembert.upmc.fr}}\\[2.0ex]
Flaviana Iurlano,\\[1.5ex]
{\small Dipartimento di Matematica,}\\
{\small Università degli Studi di Genova,}\\
{\small Via Dodecaneso 35, 16146 Genova, Italia}\\
{\small \texttt{flaviana.iurlano@unige.it}}}
\begin{document}
\maketitle


\begin{abstract}
	We formulate the quasi-static elastic contact problem with Coulomb friction in a very general setting, with possible jumps in time for both the load and the solution.  Exploiting ideas originating in our recent paper~\cite{BI-ma3as}, we exhibit an optimal condition on the magnitude of the friction coefficient under which we prove the existence of an absolutely continuous solution for arbitrary absolutely continuous loads in the case of the most general 2D problem.  We provide examples showing that, when the condition is violated, \emph{spontaneous jumps} in time of the solution may occur, even when the load varies absolutely continuously in time.  We argue that these spontaneous jumps in time of the solution in the quasi-static problem reveal a transition of the process from a quasi-static nature to a dynamic nature, interpreted as the mathematical signature of the onset of friction-induced vibrations in the elastodynamic contact problem with dry friction.
\end{abstract}

\section{Introduction}

It is common experience to observe vibratory responses in elastic contact problems with dry friction, while the load varies arbitrarily slowly in time. An everyday life example is the squeal of a car brake. However, there are also situations in which a car brake remains perfectly silent, and observations of non-vibratory responses in elastic contact problems with dry friction are not rare.

The explanation of the origin of this phenomenon and its prediction require the analysis of the mathematical problem raised by the elasticity equations coupled with contact conditions and dry friction at the boundary.  Ideally, the analysis should be carried out in the framework of elastodynamics.  Unfortunately, this is completely out of reach in the current state of knowledge, as no existence of solutions for the frictionless elastodynamic contact problem (the so-called dynamic Signorini problem) has ever been proved yet.  The only viable alternative is to analyze the \emph{quasi-static} elastic contact problem with dry friction, in the hope that transitions from slow responses to fast responses could be detected in this simplified framework.  This quirky hope stems from the following observation.  Dry friction laws, such as the ubiquitous Coulomb law, are \emph{rate-independent} (invariant under time rescaling). Hence, the quasi-static elastic contact problem with dry friction falls back into the class of rate-independent processes, such as perfect elasto-plasticity~\cite{DalMasoPlast} or brittle fracture~\cite{FrancfortMarigo}.  An extensive theory for rate-independent processes is available~\cite{FirstMielke}.  In the case of brittle fracture, finite jumps (in time) in the solution produced by infinitesimal changes in the load have been evidenced in the quasi-static theory.  This was interpreted~\cite{DalMasoLarsen} as the fact that unstable propagation makes brittle fracture inherently a dynamic process, which should be analyzed in the framework of elastodynamics.  Besides, this is evidenced by the sound emitted by crack propagation in glass.  As the dynamic nature of brittle fracture is detected by the emergence of jumping-in-time solutions in the quasi-static theory, it is natural to expect that vibratory responses in elastic contact problems with dry friction with slowly varying loads might similarly give rise to such discontinuities in the rate-independent quasi-static framework, as a mathematical signature of the dynamic nature of the process.  If it were the case, transitions from continuous to jumping solutions in the quasi-static elastic contact problem with dry friction and continuous loads would be of upmost interest.  They could serve as indicators of transitions from non-vibratory responses to vibratory responses in the elastodynamic contact problem with dry friction. That would pave the way towards a quantitative theory capable of predicting the onset of friction-induced vibrations such as brake squeal.

This paper substantiates the above-mentioned hope in the following way. We will provide a general formulation of the quasi-static elastic contact problem with Coulomb friction (Section~\ref{sec:quasiStaticProblem}), so that we will be able to state precisely what is meant by a continuous and a jumping (in time) solution. To the best of our knowledge, this is the first formulation of the problem that accounts for the possibility of jumps in the solution.  We will provide a sufficient condition on the magnitude of the friction coefficient, depending only on the geometry and the elastic moduli, ensuring that there exists a continuous-in-time solution for arbitrary continuous-in-time load, in the 2D setting (Section~\ref{sec:quasiStaticProblem} for the statement of the main result and Section~\ref{sec:ProofExistence} for its proof). We will also provide examples of jumping-in-time solutions despite continuous-in-time load, including a case where a jump is unavoidable  (while the aforementioned condition is violated, Section~\ref{sec:Examples}).  Finally, we will provide reasonable evidence of the optimality of our condition ensuring the existence of continuous solutions.

The analysis of the quasi-static elastic contact problem with Coulomb friction is a long-standing problem in the mathematical community. It has remained essentially unsolved since its first statement in the early seventies \cite{DuvautLions}. The reason is that the analysis involves products of weakly converging sequences, which lie beyond the reach of standard compactness or compensated compactness methods. In previous works, this issue was circumvented by relying on Nirenberg's shift technique to obtain the lacking compactness, at the price of nonphysical smallness conditions on the friction coefficient, both in the case of the time-discretized \cite{EckJarusek} and time-continuous \cite{Andersson} problems. The nonphysical condition on the friction coefficient has made impossible any analysis about transitions between qualitatively different solutions, such as continuous and jumping solutions, in the case of the time-continuous problem. In contrast, the present paper overcomes the lacking compactness difficulty by adapting a refined compactness argument, less restrictive than classical results and free of unphysical assumptions, originally introduced in our earlier work \cite{BI-ma3as} in view of solving the time-discretized problem. The passage to the continuous-time limit is achieved under a friction condition that is both optimized and physically meaningful, and is carried out in the generalized framework of right-continuous functions of bounded variation. This broader functional setting enables us to exhibit a sharp counterexample: absolutely continuous loads may give rise to non-absolutely continuous solutions when the optimal condition is violated. This demonstrates not only the necessity of our assumptions but also the tightness of the theory developed here. On a broader scale, this problem falls into the class of rate-independent processes with threshold depending on the current solution. It can even be viewed as the purest occurrence of such a problem. The specific techniques developed in this paper have the robustness to hopefully apply to a wider class of problems of that type.

We also argue that the analysis in this paper provides a first rigorous mathematical answer to a long-standing issue in the physics and engineering communities: is the basic Coulomb friction law able to account for the onset of friction-induced vibrations or is it necessary to consider some more sophisticated friction laws (for example, softening friction laws involving so-called `static' and `dynamic' friction coefficients)?  The analysis of the simplest two-degrees-of-freedom elastic quasi-static system (extensively discussed in Section~\ref{sec:discrete} of this paper) had evidenced the existence of a critical value of the friction coefficient expressed in terms of the elastic properties of the system. This critical value was associated \cite{Klarbring}, on the one hand, with a loss of uniqueness of solution for the incremental problem, and on the other hand, with the possible loss of either existence or uniqueness of solution for the rate problem. The loss of existence of solution for the rate problem had been informally related to the possible occurrence of jumps-in-time of the quasi-static system \cite{Klarbring}, while, to the best of our knowledge, no precise definition of what would mean a jumping-in-time solution has ever appeared. This paper places those results known for the simplest two-degrees-of-freedom quasi-static system in the proper frame and extends them to the case of a continuum. 

\subsection{Formulation of the quasi-static elastic contact problem with Coulomb friction}

The Coulomb law of dry friction is an empirical law.  It was proposed at the end of the 18th century by French physicist Charles-Augustin de Coulomb on the basis of his experimental analysis of the motion of sliding plates subjected to a precisely controlled overall applied force \(\mathbf{t}\).  It relates the sliding velocity \(\dot{\mathbf{u}}_{\rm t}\in \mathbb{R}^{N}\) to the tangential force \(\mathbf{t}_{\rm t}\in\mathbb{R}^{N}\) and the normal force \(t_{\rm n}\leq 0\).  It reads as:
\begin{align*}
	\dot{\mathbf{u}}_{\rm t} = \mathbf{0} & \quad\implies\quad |\mathbf{t}_{\rm t}| \leq -f\,t_{\rm n},\\
	\dot{\mathbf{u}}_{\rm t} \neq \mathbf{0} & \quad\implies\quad\mathbf{t}_{\rm t} = f\,t_{\rm n} \,\frac{\dot{\mathbf{u}}_{\rm t}}{|\dot{\mathbf{u}}_{\rm t}|},
\end{align*}
where \(|\cdot|\) stands for the Euclidean norm and \(f\geq 0\) is a given dimensionless, material-dependent friction coefficient.  The Coulomb law of dry friction is equivalent to the following weak formulation:
\[
	\forall\mathbf{v}\in\mathbb{R}^N,\qquad \mathbf{t}_{\rm t}\cdot\bigl(\mathbf{v}-\dot{\mathbf{u}}_{\rm t}\bigr) - f\,t_{\rm n}\bigl(|\mathbf{v}| - |\dot{\mathbf{u}}_{\rm t}|\bigr) \geq 0,
\]
where `\(\cdot\)' stands for the usual Euclidean scalar product.  In the particular case \(f=0\) (no friction), the Coulomb friction law reduces to the frictionless condition \(\mathbf{t}_{\rm t}=\mathbf{0}\).

Our aim in this paper is to study the coupling of the Coulomb law of dry friction with linear elasticity.  Consider a smooth bounded open subset \(\Omega\) of \(\mathbb{R}^N\) (\(N=2\) or \(N=3\)), whose boundary is the union of three disjoint parts \(\partial\Omega = \overline{\Gamma}_{U}\cup\overline{\Gamma}_{T}\cup\overline{\Gamma}_{C}\).  The domain \(\Omega\) is the reference configuration of a linearly elastic body.  We will prescribe respectively Dirichlet conditions on \(\Gamma_U\), Neumann conditions on \(\Gamma_T\) and contact conditions on \(\Gamma_C\).  The displacement is denoted by \(\mathbf{u}:\Omega \rightarrow\mathbb{R}^N\), the (linearized) strain by \(\boldsymbol{\varepsilon}(\mathbf{u}):=(\nabla\mathbf{u}+{}^t\nabla\mathbf{u})/2\), the stress by \(\boldsymbol{\sigma}(\mathbf{u})=\boldsymbol{\Lambda}\,\boldsymbol{\varepsilon}(\mathbf{u})\), where the elastic modulus tensor \(\boldsymbol{\Lambda}\) is assumed to satisfy the usual symmetry condition:
\begin{equation}
	\label{eq:reqLamda1}
	\forall\boldsymbol{\varepsilon}_1,\boldsymbol{\varepsilon}_2\in \mathbb{M}^{N\times N}_\text{sym},\quad \forall x\in \Omega,\qquad \boldsymbol{\varepsilon}_1:\boldsymbol{\Lambda}(x)\,\boldsymbol{\varepsilon}_2 = \boldsymbol{\varepsilon}_2:\boldsymbol{\Lambda}(x)\,\boldsymbol{\varepsilon}_1,
\end{equation}
and the strong ellipticity condition:
\begin{equation}
	\label{eq:reqLambda2}
	\exists \alpha>0,\quad \forall \boldsymbol{\varepsilon}\in \mathbb{M}^{N\times N}_\text{sym},\quad \forall x\in \Omega,\qquad \boldsymbol{\varepsilon}:\boldsymbol{\Lambda}(x)\,\boldsymbol{\varepsilon} \geq \alpha\, \boldsymbol{\varepsilon}:\boldsymbol{\varepsilon}.  
\end{equation}
Above, `:' stands for the Frobenius scalar product between real matrices $A:B:=\sum_{i,j}A_{ij}B_{ij}$.  The outward unit normal to \(\Omega\) will be denoted by \(\mathbf{n}\) and any vector \(\mathbf{v}:\partial\Omega \rightarrow \mathbb{R}^N\) will be split into normal and tangential parts: \(\mathbf{v} = v_{\rm n}\mathbf{n} + \mathbf{v}_{\rm t}\) where the scalar product \(\mathbf{n}\cdot \mathbf{v}_{\rm t}=0\) vanishes.  The load consists in a given volume force \(\mathbf{F}:\Omega \rightarrow\mathbb{R}^N\) and a given surface force \(\mathbf{T}:\Gamma_T \rightarrow\mathbb{R}^N\) on \(\Gamma_T\), both possibly varying along time \(s\in [0,S]\).  The initial gap in the direction \(\mathbf{n}\) between the elastic body \(\Omega\) and a given rigid obstacle is represented by a function \(g:\Gamma_C \rightarrow \mathbb{R}\), while the given friction coefficient between the two objects is denoted by \(f:\Gamma_C \rightarrow \left[0,+\infty\right[\).  The formal quasi-static elastic contact problem with Coulomb friction consists in finding a displacement \(\mathbf{u}(s):\Omega\to\mathbb{R}^N\) defined for time \(s\in [0,S]\), satisfying a given initial condition and, defining the surface traction $\mathbf{t}:= \boldsymbol\sigma(\mathbf{u})\,\mathbf{n}=t_{\rm n}\mathbf{n}+\mathbf{t}_{\rm t}$:
\begin{equation}
	\label{eq:SignoriniCoulombCont}
	\left\{\quad
	\begin{aligned}
		& \text{div}\,\boldsymbol{\sigma}(\mathbf{u}) + \mathbf{F} = \mathbf{0}, \qquad & & \text{ in }\Omega,\\
		& \mathbf{u} = \mathbf{0},\qquad\qquad & & \text{ on }\Gamma_U,\\
		& \boldsymbol\sigma(\mathbf{u)\,\mathbf{n}} = \mathbf{T},\qquad & & \text{ on }\Gamma_T,\\
		& u_{\rm n} - g \leq 0, \qquad t_{\rm n} \leq 0, \qquad t_{\rm n}\,(u_{\rm n}-g) = 0, \qquad & & \text{ on }\Gamma_C,\\
		& \forall \mathbf{v}\in\mathbb{R}^N,\qquad \mathbf{t}_{\rm t}\cdot\bigl(\mathbf{v}-\dot{\mathbf{u}}_{\rm t}\bigr) - ft_{\rm n} \bigl(|\mathbf{v}|-|\dot{\mathbf{u}}_{\rm t}|\bigr) \geq 0, \qquad & & \text{ on }\Gamma_C,
	\end{aligned}
	\right.
\end{equation}
where the dot stands for the time derivative.

Clearly, system \eqref{eq:SignoriniCoulombCont} is invariant under any time monotone reparametrization, making it a rate-independent process. The preferred approach for analyzing a rate-independent process is to introduce a time discretization \(0=s_0<s_1<\cdots<s_i<\cdots<s_m=S\), and to consider the problem raised on one time step by replacing the velocity \(\dot{\mathbf{u}}\) by \((\mathbf{u}_i - \mathbf{u}_{i-1})/(s_i-s_{i-1})\).  Denoting by \(\mathbf{w}:=\mathbf{u}_{i-1}\) the displacement at the preceding time step, which is supposed to be given, the formal problem on one time step is now to find a (time independent) displacement \(\mathbf{u}:\Omega \rightarrow\mathbb{R}^N\) satisfying:
\begin{subequations}
	\label{eq:SignoriniCoulombDiscret}
	\begin{empheq}[left={\empheqlbrace\quad}]{align}
		& \text{div}\,\boldsymbol{\sigma}(\mathbf{u}) + \mathbf{F} = \mathbf{0}, \qquad & & \text{ in }\Omega,\\
		& \mathbf{u} = \mathbf{0},\qquad\qquad & & \text{ on }\Gamma_U,\\
		& \boldsymbol\sigma(\mathbf{u})\,\mathbf{n} = \mathbf{T},\qquad & & \text{ on }\Gamma_T,\\
		& u_{\rm n} - g \leq 0, \qquad t_{\rm n} \leq 0, \qquad t_{\rm n}\,(u_{\rm n}-g) = 0, \qquad & & \text{ on }\Gamma_C,\label{eq:Signorini}\\
		& \forall \mathbf{v}\in\mathbb{R}^N,\qquad \mathbf{t}_{\rm t}\cdot\bigl(\mathbf{v}-\mathbf{u}_{\rm t}\bigr) - ft_{\rm n} \bigl(|\mathbf{v}-\mathbf{w}_{\rm t}|-|\mathbf{u}_{\rm t}-\mathbf{w}_{\rm t}|\bigr) \geq 0, \qquad & & \text{ on }\Gamma_C.
		\label{eq:Frict}
	\end{empheq}
\end{subequations}

In the sequel, we will refer to problem~\eqref{eq:SignoriniCoulombCont} as the \emph{quasi-static} problem, and to problem~\eqref{eq:SignoriniCoulombDiscret} as the \emph{incremental} problem.

\subsection{Historical background}

The mathematical analysis of contact problems with Coulomb friction in linear elasticity was started by Duvaut and Lions~\cite{DuvautLions} in the early seventies. They considered only the incremental problem~\eqref{eq:SignoriniCoulombDiscret}.  They observed that when \(-ft_n\) is replaced with a given \(\tau\geq 0\) in line~\eqref{eq:Frict} of problem~\eqref{eq:SignoriniCoulombDiscret} (incremental Coulomb law), then one gets the following minimum problem:
\[
\mathbf{u} = \argmin_{\substack{\mathbf{v},\\
\mathbf{v}=\mathbf{0},\text{ on }\Gamma_U,\\
v_{\rm n}\leq g,\text{ on }\Gamma_C,\\
}} \frac{1}{2}\int_{\Omega} \boldsymbol{\varepsilon}(\mathbf{v}):\boldsymbol{\Lambda}\boldsymbol{\varepsilon}(\mathbf{v}) \,{\rm d}x - \int_{\Omega} \mathbf{F}\cdot\mathbf{v} \,{\rm d}x - \int_{\Gamma_T}\mathbf{T}\cdot\mathbf{v} + \int_{\Gamma_C}\tau\,\bigl|\mathbf{v}_{\rm t}-\mathbf{w}_{\rm t}\bigr|,
\]
which can be uniquely solved by the direct method of the Calculus of Variations, under appropriate regularity assumptions on the data.  This remark suggests a fixed point strategy applied to the mapping \(\tau\mapsto -ft_n\), in order to solve the incremental problem~\eqref{eq:SignoriniCoulombDiscret}.  But no sufficient regularity of that mapping could be obtained to run this approach.

The next progress came from Jarušek's PhD thesis~\cite{JarusekPhD1,JarusekPhD2} who developed an original idea of Nečas under his supervision.  He was able to run successfully the fixed-point strategy by applying Tikhonov's theorem, which is a version of Brouwer's theorem in a locally convex topological vector space such as a Hilbert space endowed with the weak topology.  The compactness needed to apply this theorem was obtained by requiring additional regularity on the data and proving corresponding regularity on the solution by use of local rectification together with Nirenberg's shift technique, an approach that relies on highly technical arguments.  This way, the first existence result for problem~\eqref{eq:SignoriniCoulombDiscret} was proved, provided the friction coefficient was small enough.  An essential limitation of the result was that the condition on the friction coefficient, coming from the shift technique, was hardly related to the physics of the problem.  An alternate approach based on a penalty approximation of the problem was designed later by Eck and Jarušek~\cite{EckJarusek}. Nirenberg's shift technique was also invoked there, still resulting in a nonphysical condition on the friction coefficient.  Shortly after, Andersson~\cite{Andersson} published the first existence result for the quasi-static problem~\eqref{eq:SignoriniCoulombCont}. His proof was based on the existence result of Eck and Jarušek for the incremental problem and therefore required the nonphysical condition on the friction coefficient stemming from Nirenberg's shift technique. In addition, to be able to pass to the limit in the time discretization, Andersson had to assume another condition on the friction coefficient. He therefore had two conditions on the friction coefficient, which, though formally similar, were completely unrelated one with the other. Since the condition of Eck and Jarušek was obviously not related to the physics, Andersson didn't discuss at all the conditions under which the existence of a solution for the quasi-static problem was proved, and, in particular, did not comment at all about his own new condition.  Although an undeniably valuable contribution, Andersson's result was obscured by the two conditions he introduced which were clearly not optimal and whose connection to the physics was not clear and not discussed.  As far as we are aware,  Andersson's article~\cite{Andersson} is the only single article ever published, that contains results on the mathematical analysis of the quasi-static elastic contact problem with Coulomb friction~\eqref{eq:SignoriniCoulombCont}.

It is absolutely striking that, despite the apparent simplicity of the equations of the problem, their ubiquity in applications and their potentiality in predicting the onset of vibrations in slow frictional sliding, the number of articles published in the last 50 years, containing results on the mathematical analysis of either problem~\eqref{eq:SignoriniCoulombDiscret} or problem~\eqref{eq:SignoriniCoulombCont} is so small (smaller than 10).  As is often the case with fundamental problems that resist analysis, numerous studies have explored various regularizations of the problem. However, these efforts have contributed little to advancing a deeper understanding of the core issue.

In a recently published article~\cite{BI-ma3as}, we have come up with a new strategy to prove the existence of solutions for the incremental problem~\eqref{eq:SignoriniCoulombDiscret}. This new strategy yielded in particular the existence of solution for the incremental problem in 2D isotropic elasticity without any restriction on the friction coefficient, definitely proving the nonphysical character of the Eck and Jarušek condition. In this paper, we are going to prove that this new strategy can be adapted to the quasi-static problem~\eqref{eq:SignoriniCoulombCont}, yielding the existence of solutions for the quasi-static problem in 2D, under a condition on the friction coefficient which is presumably optimal, supported by strong arguments. It will turn out that our condition on the friction coefficient is physically meaningful, as it is the condition preventing \emph{spontaneous jumps in time} in the solutions.  As the analysis in this paper strongly relies on the ideas and results of~\cite{BI-ma3as}, we now provide a compendium, for the convenience of the reader.

\subsection{General overview on previous analysis of the incremental problem}

\label{sec:overviewIncremental}
Our new strategy in~\cite{BI-ma3as} for analyzing the incremental problem~\eqref{eq:SignoriniCoulombDiscret} was quite natural.  It simply consists in splitting the incremental problem into the following two subproblems:
\begin{itemize}
	\item Subproblem 1. Given \(t_{\rm n}\leq 0\) on \(\Gamma_C\), find a displacement \(\mathbf{u}:\Omega\to\mathbb{R}^N\) satisfying:
		\[
		\mathbf{u} = \argmin_{\substack{\mathbf{v},\\
		\mathbf{v}=\mathbf{0},\text{ on }\Gamma_U,\\
		}} \frac{1}{2}\int_{\Omega} \boldsymbol{\varepsilon}(\mathbf{v}):\boldsymbol{\Lambda}\boldsymbol{\varepsilon}(\mathbf{v}) \,{\rm d}x - \int_{\Omega} \mathbf{F}\cdot\mathbf{v} \,{\rm d}x - \int_{\Gamma_T}\mathbf{T}\cdot\mathbf{v} - \int_{\Gamma_C}t_{\rm n}v_{\rm n}-\int_{\Gamma_C}ft_{\rm n}\,\bigl|\mathbf{v}_{\rm t}-\mathbf{w}_{\rm t}\bigr|;
		\]
	\item Subproblem 2. Find \(t_{\rm n}\leq 0\) on \(\Gamma_C\), such that:
		\[
			\forall t\leq0, \qquad \int_{\Gamma_C}\bigl(u_{\rm n}(t_{\rm n})-g\bigr)\bigl(t-t_{\rm n}\bigr) \geq 0.
		\]
\end{itemize}
Under usual regularity assumptions on the data (to be stated precisely in the sequel), Subproblem~1 is a standard coercive, continuous, strictly convex minimum problem, having a unique solution \(\mathbf{u}(t_{\rm n})\). As \(t_{\rm n}\leq 0\) has been picked up arbitrarily, the minimizer \(\mathbf{u}(t_{\rm n})\) satisfies every line in problem~\eqref{eq:SignoriniCoulombDiscret}, except for the contact conditions~\eqref{eq:Signorini}. Subproblem~2 is precisely the variational inequality that encodes the contact condition. Hence, we view the  incremental problem~\eqref{eq:SignoriniCoulombDiscret} as a variational inequality associated with a nonlinear operator that encodes elasticity and friction.  The benefit of doing so is that one can then rely on Brezis' contribution~\cite{BrezisIeq} to the theory of variational inequalities, aiming at identifying the most general class of operator for which one can solve the variational inequality for a Galerkin approximation of the operator, using Brouwer's theorem and pass to the limit in the variational inequality. Brezis calls the corresponding operators \emph{pseudomonotone}.  Hence, solvability of the incremental problem~\eqref{eq:SignoriniCoulombDiscret} is reduced to proving that the operator \(u_n(t_n)\) associated with Subproblem 1 is pseudomonotone and coercive, which ensures solvability of Subproblem~2 by Brezis' result. 

Coercivity is readily proved in arbitrary dimension.  Pseudomonotonicity is more challenging.  As it splits the condition to check into four subconditions, we undertook to prove that the operator \(u_n(t_n)\) belongs to a subclass of pseudomonotone operators, called the class of Leray-Lions operators.  The first three subconditions are readily checked in arbitrary dimension. The fourth subcondition requires showing that one can pass to the limit in a product of two weakly converging sequences. In~\cite{BI-ma3as}, we were only able to do so in 2D, by establishing a new compensated-compactness-like fine property of the elastic Neumann-to-Dirichlet operator in 2D (see Theorem~\ref{theo:BIcompensatedCompactness} in this article). This new fine property has yielded unconditional existence of solutions for the incremental problem~\eqref{eq:SignoriniCoulombDiscret} in 2D isotropic (either homogeneous or heterogeneous) elasticity, with arbitrarily large friction coefficient.  In the case of anisotropic elasticity, a condition on the magnitude of the friction coefficient shows up. This condition has a local character and takes the form:
\begin{equation}
	\label{eq:ConditionFrictionAnisotropic}
	\text{for a.a. }x\in\Gamma_C,\qquad f(x) \leq f_{\rm max}(\Lambda(x)),
\end{equation}
where \(f_{\rm max}(\Lambda)\) is a maximal admissible value depending only on the elastic moduli \(\Lambda\) at the corresponding point \(x\in\Gamma_C\).  Again, this maximal admissible value is infinite in the case of isotropy, so that \cite{BI-ma3as} was the first article providing existence results for the incremental problem~\eqref{eq:SignoriniCoulombDiscret} with arbitrarily large friction coefficient.  Reasonable evidence was provided in~\cite{BI-ma3as} that condition~\eqref{eq:ConditionFrictionAnisotropic} is optimal, in the sense that it is necessary and sufficient to solve the steady sliding frictional contact problem for the anisotropic homogeneous half-space.  Hence, this new strategy has yielded complete understanding of the incremental problem~\eqref{eq:SignoriniCoulombDiscret} in 2D.  The checking of the fourth subcondition in higher dimension is still a work in progress and raises deep issues in concentration-compactness.

\subsection{General overview of our strategy to analyze the quasi-static problem}

\label{sec:overviewQuasiStatic}
Our strategy to analyze the quasi-static problem~\eqref{eq:SignoriniCoulombCont} is very similar to the one we used for the incremental problem~\eqref{eq:SignoriniCoulombDiscret}.  It consists in splitting the quasi-static problem into the following two subproblems:
\begin{itemize}
	\item Subproblem 1. Given \(t_{\rm n}\leq 0\) on \([0,S]\times\Gamma_C\), find a displacement \(\mathbf{u}:[0,S]\times\Omega\to\mathbb{R}^N\) satisfying, for all test function \(\mathbf{v}:\Omega\to\mathbb{R}^N\), vanishing on \(\Gamma_U\), and all \(s\in [0,S]\):
	\begin{multline*}
		\int_{\Omega} \boldsymbol{\varepsilon}(\mathbf{u}(s)):\boldsymbol{\Lambda}\boldsymbol{\varepsilon}\bigl(\mathbf{v}-\mathbf{\dot{u}}(s)\bigr) \geq \int_{\Omega} \mathbf{F}(s)\cdot\bigl(\mathbf{v}-\mathbf{\dot{u}}(s)\bigr) + \int_{\Gamma_T}\mathbf{T}(s)\cdot\bigl(\mathbf{v}-\mathbf{\dot{u}}(s)\bigr) +\mbox{}\\\mbox{}+ \int_{\Gamma_C}t_{\rm n}(s)\,\Bigl[v_{\rm n}-\dot{u}_{\rm n}(s)\Bigr] +\int_{\Gamma_C} ft_{\rm n}(s)\,\Bigl[\bigl|\mathbf{v}_{\rm t}\bigr|-|\mathbf{\dot{u}}_{\rm t}(s)\bigr|\Bigr];
	\end{multline*}
	\item Subproblem 2. Find \(t_{\rm n}\leq 0\) on \(\Gamma_C\), such that:
		\[
			\forall s\in [0,S]\quad\forall t\leq0, \qquad \int_{\Gamma_C}\Bigl(u_{\rm n}(s)-g\Bigr)\Bigl(t-t_{\rm n}(s)\Bigr) \geq 0.
		\]
\end{itemize}
Under usual regularity assumptions on the data (to be stated precisely in the sequel), Subproblem~1 is a classical evolution problem introduced and studied in the seventies by Moreau~\cite{Moreau}, under the name of \emph{sweeping process}. A sweeping process can roughly be viewed as a succession of strictly convex minimum problems, parametrized by time \(s\). 
Given an initial condition, it can be uniquely solved. Very interestingly for our purpose, this evolution inequality can be given a sense, even when data \(\mathbf{F}(s)\), \(\mathbf{T}(s)\), \(t_{\rm n}(s)\) and solution \(\mathbf{u}(s)\) exhibit jumps in time. In the case of the sweeping process, absolutely continuous data always yield absolutely continuous solutions.  Hence, as was the case for the incremental problem, Subproblem~1 yields an operator \(u_n(t_n)\) acting on \(t_n\) viewed as a function of time and space.  However, this operator lacks both the coercivity and the pseudomonotonicity properties needed to apply Brezis' result in order to solve Subproblem~2.  Hence, one has to rely on the solvability of the incremental problem and to pass to the limit in the time discretization.  This approach fits harmoniously with the sweeping process theory, whose solvability is precisely proved on the basis of such time stepping approximations. Some difficulties will arise on the way, though. First, we will have to adapt some pseudomonotonicity arguments, and this step will rely crucially on the new fine property of the elastic Neumann-to-Dirichlet operator that we proved in \cite{BI-ma3as}.  Second, the necessary estimates to pass to the limit in the time discretization will require to assume a condition on the friction coefficient.  We will take great care to write the minimal condition.  Under this condition, we then will be able to prove the existence of a solution for the quasi-static problem~\eqref{eq:SignoriniCoulombCont}, which is absolutely continuous in time, provided that the load is also absolutely continuous in time.
Detailed analysis of a spatially discrete counterpart of the quasi-static problem (for all \(s\in [0,S]\), \(\mathbf{u}(s)\in\mathbb{R}^2\)) will then yield reasonable evidence that the condition on the friction coefficient is optimal.  In particular, when the condition does not hold, we will provide an example of a case where the solution \emph{must jump}, while the load varies continuously in time.

\section{Setting of the problem and main results}

\subsection{Geometry and traces}

We are given a bounded open set \(\Omega\subset \mathbb{R}^N\), connected and of class \(C^{1,1}\), which entails in particular that the outward unit normal \(\mathbf{n}\) belongs to \(W^{1,\infty}(\partial\Omega;\mathbb{R}^N)\).  We consider three nonintersecting open subsets \(\Gamma_U\), \(\Gamma_T\) and \(\Gamma_C\) of the boundary \(\partial \Omega\) of \(\Omega\), such that \(\partial\Omega=\overline{\Gamma}_U\cup\overline{\Gamma}_T\cup\overline{\Gamma}_C\) (where \(\overline{\Gamma}\) stands for the closure of \(\Gamma\) in \(\partial\Omega\)).  For the sake of simplicity, we make the following additional hypotheses.
\begin{itemize}
	\item The subset \(\Gamma_U\) has positive surface measure \(\mathcal{H}^{N-1}\).  As a homogeneous Dirichlet condition will be prescribed on \(\Gamma_U\), this hypothesis is made for convenience, to obtain some coercivity.  It is not essential and could be dropped at the price of additional complexity of the presentation.
	\item The subset \(\Gamma_C\subset \partial\Omega\) is of class \(C^{0,1}\) (that is, Lipschitz) and \(\text{dist}(\Gamma_U,\Gamma_C)>0\).  Again, this last hypothesis is made for convenience and is probably not essential.  It ensures that the space of traces on \(\Gamma_C\) of functions in \(H^1(\Omega)\) that vanish on \(\Gamma_U\):
	\[
	\Bigl\{ u_{|\Gamma_C} \bigm| u\in H^1(\Omega) \; \text{ and }\;u=0\;\text{ on }\Gamma_U\Bigr\}
	\]
	is exactly \(H^{1/2}(\Gamma_C)\).  Indeed, the space of traces of functions in \(H^1(\Omega)\) is \(H^{1/2}(\partial\Omega)\) and restrictions to \(\Gamma_C\) of functions in \(H^{1/2}(\partial\Omega)\) are in \(H^{1/2}(\Gamma_C)\).  Vice versa, when \(\text{dist}(\Gamma_U,\Gamma_C)>0\), it is always possible to find an extension to \(\partial\Omega\) of a given function in \(H^{1/2}(\Gamma_C)\), so that this extension vanishes identically in \(\Gamma_U\).  When this hypothesis is not made, the above trace space can be strictly smaller than \(H^{1/2}(\Gamma_C)\) (see \cite[Section~11]{LionsMagenes}).
\end{itemize}

Considering an arbitrary open \(C^1\) manifold \(\omega\) with dimension \(d\), embedded in an Euclidean space, the spaces \(L^2(\omega)\), \(L^\infty(\omega)\), \(H^1(\omega)\), \(W^{1,\infty}(\omega)\) and \(H^{1/2}(\omega)\) are defined in terms of their respective norms:
\begin{align}
	\|v\|_{L^2(\omega)} & := \left(\int_{\omega} |v|^2\,{\rm d}\mathcal{H}^{d}\right)^{1/2},\nonumber\\
	\|v\|_{L^\infty(\omega)} & := \esup_{x\in\omega} |v(x)|,\nonumber\\
	\|v\|_{H^1(\omega)} & := \left(\|v\|_{L^2(\omega)}^2 + \|\nabla v\|_{L^2(\omega)}^2\right)^{1/2},\nonumber\\
	\|v\|_{W^{1,\infty}(\omega)} & := \max\left\{\|v\|_{L^\infty(\omega)},\|\nabla v\|_{L^\infty(\omega)}\right\},\nonumber\\
	\|v\|_{H^{1/2}(\omega)} & := \left(\int_{\omega} |v|^2\,{\rm d}\mathcal{H}^{d}+\int_{\omega\times\omega}\frac{|v(x)-v(y)|^2}{|x-y|^{d+1}}\,{\rm d}\mathcal{H}^{d}(x)\, {\rm d}\mathcal{H}^{d}(y)\right)^{1/2}.\label{eq:defH12}
\end{align}
From the definition of these norms, we have the inequalities:
\begin{equation}
	\label{eq:multLischitz}
	\Bigl\| |v|\Bigr\|_{H^{1/2}(\omega)} \leq \|v\|_{H^{1/2}(\omega)},\qquad \|fv\|_{H^{1/2}(\omega)} \leq \sqrt{2\max\{1,G(\omega)\}}\;\|f\|_{W^{1,\infty}(\omega)}\|v\|_{H^{1/2}(\omega)},
\end{equation}
where \(G(\omega)\) is the following constant (which is finite whenever \(\omega\) is bounded) depending only on the geometry of \(\omega\):
\begin{equation}
	\label{eq:geometryConstant}
		G(\omega):=\sup_{x\in\omega}\int_{\omega}\frac{1}{|x-y|^{d-1}}\,{\rm d}\mathcal{H}^{d}(y).
\end{equation}
We define the space \(H^{-1/2}(\omega)\) as the dual space of \(H^{1/2}(\omega)\) and we systematically endow it with the dual norm. Recalling that the space of \(C^\infty\) functions with compact support is dense in \(H^{1/2}(\omega)\), the space \(H^{-1/2}(\omega)\) is a space of distributions on \(\omega\). Here and henceforth $\langle\cdot,\cdot\rangle$ will denote the duality pairing. For \(w\in H^{1/2}(\Gamma_C)\) and \(\tau\in H^{-1/2}(\Gamma_C)\cap \mathscr{M}(\Gamma_C)\), we will also use the notation \(\langle w,\tau\rangle=\int_{\Gamma_C}w\tau\) which stems from the identification of \(L^2\) with its dual space.
 
In the sequel, it will be convenient to consider the closed subspace of \(H^1(\Omega;\mathbb{R}^N)\) defined by:
\begin{equation}
	\label{eq:defV}
		V:=\Bigl\{\mathbf{u}\in H^1(\Omega;\mathbb{R}^N) \bigm| \mathbf{u}_{|\Gamma_U}=\mathbf{0}\Bigr\}.
\end{equation}
Given \(\boldsymbol\Lambda\in L^\infty(\Omega)\) satisfying requirements~\eqref{eq:reqLamda1} and \eqref{eq:reqLambda2}, by Korn's inequality there exist two positive real constants \(c_0\) and \(C_0\) such that, for all \(\mathbf{u}\in V\):
\begin{equation}
	\label{eq:Korn}
		c_0^2\bigl\|\mathbf{u}\bigr\|_{H^1}^2 \leq \int_{\Omega} \boldsymbol{\varepsilon}(\mathbf{u}):\boldsymbol{\Lambda}\boldsymbol{\varepsilon}(\mathbf{u}) \,{\rm d}x \leq C_0^2\bigl\|\mathbf{u}\bigr\|_{H^1}^2,
\end{equation}
so that the inner product:
\begin{equation}
	\label{eq:InnerProduct}
	a(\mathbf{u},\mathbf{v}):=\int_{\Omega} \boldsymbol{\varepsilon}(\mathbf{u}):\boldsymbol{\Lambda}\boldsymbol{\varepsilon}(\mathbf{v}) \,{\rm d}x
\end{equation}
induces a norm on \(V\) that is equivalent to that of \(H^1(\Omega;\mathbb{R}^N)\).
The two following linear mappings are classically well-defined and continuous:
\begin{align*}
	\text{\rm tr} & \left\{\begin{array}{rcl}
		H^1(\Omega;\mathbb{R}^N) & \to & H^{1/2}(\Gamma_C\cup\Gamma_T;\mathbb{R}^N)\\[1.0ex]
		\mathbf{u} & \mapsto & \mathbf{u}_{|\Gamma_C\cup\Gamma_T}\end{array}\right.\\[3.0ex] \text{\rm ext} & \left\{\begin{array}{rcl}
	H^{1/2}(\Gamma_C;\mathbb{R}^N) & \to & H^{1}(\Omega;\mathbb{R}^N)\\[1.0ex]
	\mathbf{u} & \mapsto & \displaystyle\argmin_{\substack{\mathbf{v}\in V,\\
		\mathbf{v}=\mathbf{u},\text{ on }\Gamma_C,\\
		}} \int_{\Omega} \boldsymbol{\varepsilon}(\mathbf{v}):\boldsymbol{\Lambda}\boldsymbol{\varepsilon}(\mathbf{v}) \,{\rm d}x\end{array}\right.
\end{align*}
Considering \(\mathbf{F}\in L^2(\Omega;\mathbb{R}^N)\), \(\mathbf{T}\in L^2(\Gamma_T;\mathbb{R}^N)\), \(\mathbf{u}\in V\) such that \(\text{div}\,\boldsymbol{\sigma}(\mathbf{u})+\mathbf{F}=\mathbf{0}\) in \(\Omega\) and \(\boldsymbol\sigma(\mathbf{u})\,\mathbf{n}=\mathbf{T}\) on \(\Gamma_T\), Stoke's theorem yields, for all \(\mathbf{v}\in H^1(\Omega;\mathbb{R}^N)\) such that \(\mathbf{v}_{|\Gamma_U}=\mathbf{0}\):
\begin{equation}
	\label{eq:Stokes}
		\Bigl\langle \boldsymbol{\sigma}(\mathbf{u})\mathbf{n},\mathbf{v}_{|\Gamma_C}\Bigr\rangle = \int_{\Omega} \boldsymbol{\varepsilon}(\mathbf{u}):\boldsymbol{\Lambda}\boldsymbol{\varepsilon}(\mathbf{v}) \,{\rm d}x - \int_{\Omega} \mathbf{F}\cdot\mathbf{v} \,{\rm d}x - \int_{\Gamma_T}\mathbf{T}\cdot\mathbf{v},
\end{equation}
which enables to define \(\mathbf{t}:=\boldsymbol{\sigma}(\mathbf{u})\mathbf{n}\) as an element of \(H^{-1/2}(\Gamma_C;\mathbb{R}^N)\) (the dual space of \(H^{1/2}(\Gamma_C;\mathbb{R}^N)\)) with:
\begin{equation}
	\label{eq:tBound}
	\| \mathbf{t}\|_{H^{-1/2}} \leq \max\bigl\{1,C_0,\|\text{tr}\|\bigr\}\|\text{ext}\|\Bigl\{\bigl(a(\mathbf{u},\mathbf{u})\bigr)^{1/2} + \|\mathbf{F}\|_{L^2} + \|\mathbf{T}\|_{L^2}\Bigr\},
\end{equation}
Using the above assumed \(C^{1,1}\) regularity of \(\Omega\), we can also define the normal part \(t_{\rm n}\) of \(\mathbf{t}\) as an element of \(H^{-1/2}(\Gamma_C)\) by the formula:
\[
	\forall v\in H^{1/2}(\Gamma_C),\qquad \langle t_{\rm n},v\rangle := \bigl\langle \mathbf{t},v\,\mathbf{n}\bigr\rangle.
\]
The tangential part will be defined by \(\mathbf{t}_{\rm t}:=\mathbf{t}-t_{\rm n}\mathbf{n}\) and belongs to \(H^{-1/2}(\Gamma_C;\mathbb{R}^{N})\).

Putting everything together, we have proved the following lemma, which will prove useful in the sequel.
\begin{lemm}
	\label{lemm:condAndersson}
	Let \(\mathbf{F}\in L^2(\Omega;\mathbb{R}^N)\), \(\mathbf{T}\in L^2(\Gamma_T;\mathbb{R}^N)\), \(\mathbf{u}\in V\setminus\{\mathbf{0}\}\) such that \(\text{\rm div}\,\boldsymbol{\sigma}(\mathbf{u})+\mathbf{F}=\mathbf{0}\) in \(\Omega\) and \(\boldsymbol\sigma(\mathbf{u})\,\mathbf{n}=\mathbf{T}\) on \(\Gamma_T\).  Let also \(f\in W^{1,\infty}(\Gamma_C)\) and \(\mathbf{v}\in V\setminus\{\mathbf{0}\}\).  Then:
	\begin{multline*}
		\frac{\Bigl|\bigl\langle f\,\boldsymbol{\sigma}(\mathbf{u})\mathbf{n}\;,\;|\mathbf{v}_{\rm t}|\mathbf{n} \bigr\rangle_{\Gamma_C}\Bigr|}{\bigl(a(\mathbf{v},\mathbf{v})\bigr)^{1/2}\Bigl\{\bigl(a(\mathbf{u},\mathbf{u})\bigr)^{1/2} + \|\mathbf{F}\|_{L^2} + \|\mathbf{T}\|_{L^2}\Bigr\}}\\ \mbox{}\leq \frac{\max\bigl\{1,C_0,\|\text{\rm tr}\|\bigr\}\|\text{\rm ext}\|\|\text{\rm tr}\|}{c_0}\sqrt{2\max\{1,G(\Gamma_C)\}}\;\|f\|_{W^{1,\infty}},
	\end{multline*}
	where the geometric constant \(G(\Gamma_C)\) is defined in~\eqref{eq:geometryConstant} and \(c_0\) is the constant in Korn's inequality~\eqref{eq:Korn}.
\end{lemm}

\subsection{Rigorous results for the incremental problem}

\begin{defi}[Incremental problem]
	\label{def:IncrementalProblem}
	Consider the arbitrary data \(f\in W^{1,\infty}(\Gamma_C;\left[0,+\infty\right[)\), \(\mathbf{F}\in L^2(\Omega;\mathbb{R}^N)\), \(\mathbf{T}\in L^2(\Gamma_T;\mathbb{R}^N)\), \(g\in H^{1/2}(\Gamma_C)\) and \(\mathbf{w}_{\rm t}\in H^{1/2}(\Gamma_C);\mathbb{R}^N)\), with \(\mathbf{w}_{\rm t}\cdot\mathbf{n}=0\). 
The incremental problem is that of finding \(\mathbf{u}\in V\) such that:
\begin{align}
	& \bullet \quad\forall\mathbf{v}\in V,\qquad \int_{\Omega} \boldsymbol{\varepsilon}(\mathbf{u}):\boldsymbol{\Lambda}\boldsymbol{\varepsilon}(\mathbf{v}-\mathbf{u}) \geq \int_{\Omega} \mathbf{F}\cdot(\mathbf{v}-\mathbf{u}) + \int_{\Gamma_T}\mathbf{T}\cdot(\mathbf{v}-\mathbf{u}) \nonumber\\ & \hspace*{7cm}\mbox{} + \int_{\Gamma_C}t_{\rm n}\,(v_{\rm n}-u_{\rm n})+\int_{\Gamma_C}ft_{\rm n}\Bigl[\bigl|\mathbf{v}_{\rm t}-\mathbf{w}_{\rm t}\bigr|-\bigl|\mathbf{u}_{\rm t}-\mathbf{w}_{\rm t}\bigr|\Bigr],
	\label{eq:optimalityIncrementalPb}\\
	& \bullet \quad t_{\rm n}\leq 0\;\text{ on }\Gamma_C,\quad \text{ and }\quad	\forall t\in H^{-1/2}(\Gamma_C),\text{ such that }\:t\leq 0,\quad \int_{\Gamma_C}\bigl(u_{\rm n}-g\bigr)\bigl(t-t_{\rm n}\bigr)\geq 0,\label{eq:SignoriniIncrementalPb}
\end{align}
where \(t_{\rm n}\) stands for the normal component of \(\mathbf{t}:=\boldsymbol{\sigma}(\mathbf{u})\mathbf{n}\) on \(\Gamma_C\).
\end{defi}

Note that whenever \(\tau=t_{\rm n}\leq 0\) is given, then the first statement in Definition~\ref{def:IncrementalProblem} is equivalent (necessary and sufficient optimality condition) to the strictly convex minimum problem:
\[
\mathbf{u} = \argmin_{\substack{\mathbf{v}\in V\\
	}} \frac{1}{2}\int_{\Omega} \boldsymbol{\varepsilon}(\mathbf{v}):\boldsymbol{\Lambda}\boldsymbol{\varepsilon}(\mathbf{v}) \,{\rm d}x - \int_{\Omega} \mathbf{F}\cdot\mathbf{v} \,{\rm d}x - \int_{\Gamma_T}\mathbf{T}\cdot\mathbf{v} - \int_{\Gamma_C}t_{\rm n}v_{\rm n}-\int_{\Gamma_C}f\tau\,\bigl|\mathbf{v}_{\rm t}-\mathbf{w}_{\rm t}\bigr|.
\]

In a recently published article~\cite{BI-ma3as}, we proved the following existence result for the incremental problem in 2D, which is the first result of this type with arbitrarily large friction coefficient.

\begin{theo}[{\cite[Theorem 2.13]{BI-ma3as}}]
	\label{theo:BIincremental}
	We suppose that \(N=2\) and \(\boldsymbol{\Lambda}\in W^{1,\infty}(\Omega)\) satisfying condition~\eqref{eq:ConditionFrictionAnisotropic} (no condition if \(\boldsymbol{\Lambda}\) is isotropic at each \(x\in \Omega\)).  Let \(f\in W^{1,\infty}(\Gamma_C;\left[0,+\infty\right[)\), \(\mathbf{F}\in L^2(\Omega;\mathbb{R}^2)\), \(\mathbf{T}\in L^2(\Gamma_T;\mathbb{R}^2)\), \(g\in H^{1/2}(\Gamma_C)\) and \(w_{\rm t}\in H^{1/2}(\Gamma_C)\) be given data.  Then, there exists a solution \(\mathbf{u}\in V\) of the incremental problem formulated in Definition~\ref{def:IncrementalProblem}.
\end{theo}

No existence result for large friction coefficient is available in higher dimension yet. The only available existence result is that of Eck and Jarušek~\cite{EckJarusek}, which relies on the nonphysical smallness condition on the friction coefficient that comes from Nirenberg's shift technique.  Regarding uniqueness, examples of multiple solutions have been exhibited~\cite{Hild} for finite friction coefficient.  But, no such examples have ever been produced for arbitrarily small friction coefficient, while no uniqueness result is available for small friction coefficient, either.

\subsection{Bounded variation in time and sweeping processes}

Let \(S\geq 0\) be a given time horizon and \(H\) a Hilbert space.   For a function \(\mathbf{u}:[0,S]\to H\), we define the total variation of \(\mathbf{u}\) in time as:
\begin{equation}
	\label{eq:defVar}
		\text{var}(\mathbf{u};[0,S]) := \sup\sum_{i=1}^m \|\mathbf{u}(s_i)-\mathbf{u}(s_{i-1})\|_H,
\end{equation}
where the supremum is taken over all finite subdivisions \(s_0<s_1<\cdots<s_m\) of the interval \([0,S]\).  The space of functions of bounded variation in time from \([0,S]\) to \(H\) is defined by: 
\[
	\textsl{BV}([0,S];H) := \Bigl\{\mathbf{u}:[0,S]\to H \bigm| \text{var}(\mathbf{u};[0,S])<+\infty\Bigr\}.
\]
Each function \(\mathbf{u}\in\textsl{BV}([0,S];H)\) has a right limit \(\mathbf{u}(s+)\) at each \(s\in \left[0,S\right[\).  By convention, we set \(\mathbf{u}(S+):=\mathbf{u}(S)\).  Likewise, there exists a left limit \(\mathbf{u}(s-)\) at each \(s\in \left]0,S\right]\), and we set \(\mathbf{u}(0-):=\mathbf{u}(0)\). The set of \(s\in [0,S]\) where \(\mathbf{u}\) has a jump (\(\mathbf{u}(s+)\neq \mathbf{u}(s-)\)) is at most countable. Having a left and right limit at each \(s\in [0,S]\), the function \(\mathbf{u}\) is classically the uniform limit on \([0,S]\) of a sequence of piecewise constant functions. This makes any \(\mathbf{u}\in\textsl{BV}([0,S];H)\) \emph{universally integrable}, that is, integrable with respect to any measure on \([0,S]\).  The space of \(H\)-valued measures with bounded total variation on \([0,S]\) is denoted by \(\mathscr{M}([0,S];H)\).  For any \(\mathbf{u}\in\textsl{BV}([0,S];H)\), its distributional derivative \(\dot{\mathbf{u}}\) belongs to \(\mathscr{M}([0,S];H)\) (\(\dot{\mathbf{u}}\) is sometimes called the Stieljes measure of \(\mathbf{u}\)) and we have:
\[
	\int_{\left]s_1,s_2\right[} \dot{\mathbf{u}} = \mathbf{u}(s_2-)-\mathbf{u}(s_1+),\qquad \int_{\left[s_1,s_2\right]} \dot{\mathbf{u}} = \mathbf{u}(s_2+)-\mathbf{u}(s_1-).
\]
By definition, any measure \(\mathbf{v}\in\mathscr{M}([0,S];H)\) has bounded total variation measure \(\|\mathbf{v}\|\in\mathscr{M}([0,S];\left[0,+\infty\right[)\). As any Hilbert space has the so-called \emph{Radon-Nikodym property} (see, e.g., \cite[Section VII.7, p.218]{diestelUhl}), any measure \(\mathbf{v}\in\mathscr{M}([0,S];H)\) admits a density \(\bar{\mathbf{v}}\in L^1([0,S],\|\mathbf{v}\|;H)\) with respect to its total variation measure \(\|\mathbf{v}\|\), which will be called in the sequel the Radon-Nikodym derivative \(\bar{\mathbf{v}}={\rm d}\mathbf{v}/{\rm d}\|\mathbf{v}\|\). Finally \cite[Propositions 11.1 and 12.1]{MoreauBV} gives: 
\begin{equation}
	\label{eq:bvDissipation}
	\dot{\wideparen{\|\mathbf{u}(s)\|_H^2}} = \bigl(\mathbf{u}(s+)+\mathbf{u}(s-)\bigr)\cdot\dot{\mathbf{u}} \leq 2\mathbf{u}(s+)\cdot\dot{\mathbf{u}},
\end{equation}
where `$\cdot$' stands for the scalar product of $H$ and \(\mathbf{u}(s+)\cdot\dot{\mathbf{u}}\) is a shorthand notation for \(\mathbf{u}(s+)\cdot\bar{\mathbf{v}}\,\|\dot{\mathbf{u}}\|\), with \(\bar{\mathbf{v}}:={\rm d}\dot{\mathbf{u}}/{\rm d}\|\dot{\mathbf{u}}\|\). In the sequel, we will focus on right-continuous functions with bounded variation and adopt the notation:
\[
	\textsl{BV}^+([0,S];H) := \Bigl\{\mathbf{u}\in\textsl{BV}([0,S]; H) \bigm| \forall s\in [0,S],\;\mathbf{u}(s)=\mathbf{u}(s+)\Bigr\}.
\]
A function \(\mathbf{u}:[0,S]\to H\) is said to be \emph{absolutely continuous} if, for all \(\epsilon>0\), there exists \(\eta>0\) such that, for all finite subdivision \(s_0<s_1<\cdots<s_m\):
\[
	\sum_{i=1}^m \bigl|s_i-s_{i-1}\bigr|<\eta, \qquad \implies\qquad\sum_{i=1}^m \bigl\|\mathbf{u}(s_i)-\mathbf{u}(s_{i-1})\bigr\|_H < \epsilon.
\]
Any absolutely continuous function has bounded variation, admits a classical derivative \(\dot{\mathbf{u}}(s)\), at almost all \(s\in [0,S]\). The (almost everywhere defined) function \(\dot{\mathbf{u}}\) is Lebesgue integrable, and satisfies~\cite[Theorem 1.124]{BarBuPrecupanu}:
\[
	\forall s_1,s_2\in [0,S],\qquad \mathbf{u}(s_2)-\mathbf{u}(s_1) = \int_{s_1}^{s_2} \dot{\mathbf{u}}(s)\,{\rm d}s.
\]
A function \(\mathbf{u}\in\textsl{BV}([0,S];H)\) is absolutely continuous (resp. right-continuous, resp. left-continuous) on \([0,S]\) if and only if the real-valued function \(s\to \text{var}(\mathbf{u};[0,s])\) is absolutely continuous (resp. right-continuous, resp. left-continuous) on \([0,S]\) \cite[Corollary 4.4]{MoreauBV}.

We now turn to state elementary facts about sweeping processes. The motivation of doing so is the following.
Our formulation of the incremental problem in Definition~\ref{def:IncrementalProblem} is a convex minimization problem (when \(t_{\rm n}\) is fixed), coupled with a variational inequality that encodes the contact condition.  This structure has played a crucial role in the analysis of the incremental problem in~\cite{BI-ma3as}. It turns out that the evolution inequality in Subproblem~1 in Section~\ref{sec:overviewQuasiStatic} (when \(t_{\rm n}\) is fixed), also carries a structure: that of a sweeping process which is essentially a convex minimization problem parametrized by time \(s\). The notation and results recalled in the remainder of this section are not strictly needed in themselves, but they serve to illustrate the underlying structure that will play an essential role later in yielding a natural formulation of the quasi-static contact problem with Coulomb friction within the framework of right-continuous functions of bounded variation. This, in turn, will open up the possibility of accounting, for the first time, for jumping solutions in such a problem.

Denoting by \(\mathscr{C}(H)\) the set of nonempty closed convex subsets of \(H\), we consider a function \(C:[0,S]\to \mathscr{C}(H)\). In the terminology of Moreau~\cite{Moreau}, \(C(s)\) is a \emph{moving convex set}  in the Hilbert space \(H\).  Considering the subdifferential \(\partial I_{C(s)}\) of the indicatrix function \(I_{C(s)}\) of \(C(s)\) (that takes the value \(0\) on \(C(s)\) and \(+\infty\) outside \(C(s)\)), and given an arbitrary initial condition \(\mathbf{u}_0\in C(0)\), the sweeping process associated with the moving convex set \(C(s)\) is the formal evolution problem of finding a function \(\mathbf{u}:[0,S]\to H\) such that:
\begin{equation}
	\label{eq:sweepingProcess}
	\mathbf{u}(0)=\mathbf{u}_0,\qquad\text{and, for a.a. } s\in [0,S],\quad -\dot{\mathbf{u}}(s)\in \partial I_{C(s)}\bigl(\mathbf{u}(s)\bigr),
\end{equation}
which amounts to requiring that the velocity \(\dot{\mathbf{u}}(s)\) is in the cone of the inward normal vectors to the convex set \(C(s)\) at \(\mathbf{u}(s)\). Roughly speaking, the moving point \(s \mapsto \mathbf{u}(s)\) remains at rest as long as it lies in the interior of C(s); when it is reached by the boundary of the moving set, it can only move in an inward normal direction, as if pushed by the boundary. It is sufficient that a nonempty closed convex moving set \(C:[0,S]\to \mathscr{C}(H)\) is absolutely continuous when \(\mathscr{C}(H)\) is endowed with the Hausdorff distance, to ensure that there exists a unique absolutely continuous function \(\mathbf{u}:[0,S]\to H\) satisfying~\eqref{eq:sweepingProcess} \cite{Moreau}.  The evolution problem~\eqref{eq:sweepingProcess} can be extended to the case where \(C(s)\) is not necessarily continuous, but is only assumed to have bounded variation.  Such is the case when \(C(s)\) is a piecewise constant function.  In that case, the solution of the sweeping process will be constructed as a piecewise constant function, in a way that the new position of the point \(\mathbf{u}_i\), corresponding to the time interval \(\left[s_i,s_{i+1}\right[\), is obtained by projecting the old position \(\mathbf{u}_{i-1}\), corresponding to the time interval \(\left[s_{i-1},s_{i}\right[\), on the new convex set \(C_i\), using the algorithm:
\[
	\mathbf{u}_i = \text{Proj}\,(\mathbf{u}_{i-1},C_i),\qquad \Longleftrightarrow\qquad  -(\mathbf{u}_i-\mathbf{u}_{i-1})\in \partial I_{C_i}(\mathbf{u}_i).
\]
This is equivalent to require that the Radon-Nikodym derivative \(\mathbf{v}:={\rm d}\dot{\mathbf{u}}/{\rm d}\|\dot{\mathbf{u}}\|\) of the measure \(\dot{\mathbf{u}}\in\mathscr{M}([0,S];H)\) with respect to its total variation measure satisfies:
\begin{equation}
	\label{eq:sweepingProcessWithJumps}
		\text{for $\|\dot{\mathbf{u}}\|$-a.a. } s\in [0,S],\quad -\mathbf{v}(s)\in \partial I_{C(s)}\bigl(\mathbf{u}(s)\bigr),
\end{equation}
provided that \emph{$C(s)$ and \(\mathbf{u}(s)\) are right-continuous functions}. Hence, as noted by Moreau~\cite{Moreau}, the right-continuous functions with bounded variation are the natural framework to study sweeping processes with jumps in time. Besides, Moreau proved in \cite{Moreau} that if \(C:[0,S]\to \mathscr{C}(H)\) is right-continuous with bounded variation (\(\mathscr{C}(H)\) being endowed with the Hausdorff distance) and \(\mathbf{u}_0\in C(0)\), then there exists a unique right-continuous function \(\mathbf{u}:[0,S]\to H\) with bounded variation satisfying both the initial condition and~\eqref{eq:sweepingProcessWithJumps}.

Now we check that Subproblem 1 in Section~\ref{sec:overviewQuasiStatic} can be interpreted as a sweeping process.

Given \(\mathbf{F}:[0,S]\to L^2(\Omega;\mathbb{R}^N)\), \(\mathbf{T}:[0,S]\to L^2(\Gamma_T;\mathbb{R}^N)\) and \(\tau:[0,S]\to H^{-1/2}(\Gamma_C)\), such that, for all \(s\), \(\tau(s)\leq 0\) on \(\Gamma_C\), the set:
\begin{multline*}
	C(s):=\Biggl\{\mathbf{u}\in V \Bigm| \forall \mathbf{v}\in V,\\
	\int_{\Omega} \boldsymbol{\varepsilon}(\mathbf{u}):\boldsymbol{\Lambda}\boldsymbol{\varepsilon}(\mathbf{v})\:{\rm d}x \geq \int_{\Omega} \mathbf{F}(s)\cdot\mathbf{v}\:{\rm d}x + \int_{\Gamma_T}\mathbf{T}(s)\cdot\mathbf{v} + \int_{\Gamma_C}\tau(s)\,v_{\rm n}+\int_{\Gamma_C}f\tau(s)\,|\mathbf{v}_{\rm t}|\Biggr\}.
\end{multline*}
is a moving convex set (in the terminology of Moreau~\cite{Moreau}) in the Hilbert space \(V\) (defined by formula~\eqref{eq:defV}). Recalling the notation $\mathbf{t}:= \boldsymbol\sigma(\mathbf{u})\,\mathbf{n}=t_{\rm n}\mathbf{n}+\mathbf{t}_{\rm t}$, \(C(s)\) is the set of those \(\mathbf{u}\in H^1(\Omega;\mathbb{R}^N)\) satisfying:
\[
\left\{\quad
	\begin{aligned}
		& \text{div}\,\boldsymbol{\sigma}(\mathbf{u}) + \mathbf{F}(s) = \mathbf{0}, \qquad & & \text{ in }\Omega,\\
		& \mathbf{u} = \mathbf{0},\qquad\qquad & & \text{ on }\Gamma_U,\\
		& \boldsymbol\sigma(\mathbf{u})\,\mathbf{n} = \mathbf{T}(s),\qquad & & \text{ on }\Gamma_T,\\
		& t_{\rm n} = \tau(s), \qquad & & \text{ on }\Gamma_C,\\
		& \forall \mathbf{v}\in H^{1/2}(\Gamma_C;\mathbb{R}^N),\qquad \bigl\langle\mathbf{t}_{\rm t}(s),\mathbf{v}_{\rm t}\bigr\rangle - \bigl\langle f\tau(s) , |\mathbf{v}_{\rm t}|\bigr\rangle \geq 0, \qquad & & \text{ on }\Gamma_C.
	\end{aligned}
\right.
\]
It is obviously a nonempty closed convex subset of \(V\), for all \(s\in [0,S]\).  In the sequel, the Hilbert space \(V\) will be endowed with the inner product~\eqref{eq:InnerProduct}.  Considering the support function \(S_{C(s)}\) of the convex set \(C(s)\):
\[
	S_{C(s)}(\mathbf{v}) = \int_{\Omega} \mathbf{F}(s)\cdot\mathbf{v}\:{\rm d}x + \int_{\Gamma_T}\mathbf{T}(s)\cdot\mathbf{v} + \int_{\Gamma_C}\tau(s)\,v_{\rm n}-\int_{\Gamma_C}f\tau(s)\,|\mathbf{v}_{\rm t}|,
\]
we have the classic equivalence:
\[
	-\dot{\mathbf{u}}(s)\in \partial I_{C(s)}\bigl(\mathbf{u}(s)\bigr) \qquad \Longleftrightarrow\qquad \mathbf{u}(s)\in \partial S_{C(s)}\bigl(-\dot{\mathbf{u}}(s)\bigr),
\]
coming from the fact that the support function \(S_C\) is the Legendre-Fenchel conjugate of the indicatrix function \(I_C\).

\begin{defi}
	\label{defi:intMeasure}
	Let \(\mathbf{u}\in \textsl{BV}([0,S];H^1(\Omega;\mathbb{R}^N))\) and \(\mathbf{v}\in \mathscr{M}([0,S];H^1(\Omega;\mathbb{R}^N))\). Denoting by \(\bar{\mathbf{v}}={\rm d}\mathbf{v}/{\rm d}\|\mathbf{v}\|\in L^1([0,S],\|\mathbf{v}\|;H^1(\Omega;\mathbb{R}^N))\) the Radon-Nikodym derivative of \(\mathbf{v}\) with respect to its total variation measure, we will systematically use the shorthand notation:
	\[
		\int_{[0,S]}\int_{\Omega} \boldsymbol{\varepsilon}(\mathbf{u}(s)):\boldsymbol{\Lambda}\boldsymbol{\varepsilon}\bigl(\mathbf{v}\bigr) := \int_{[0,S]}\int_{\Omega} \boldsymbol{\varepsilon}(\mathbf{u}(s)):\boldsymbol{\Lambda}\boldsymbol{\varepsilon}\bigl(\bar{\mathbf{v}}(s)\bigr)\,{\rm d}x\,{\rm d}\|\mathbf{v}\|.
	\]
	Likewise, with \(\tau\in \textsl{BV}([0,S];H^{-1/2}(\Gamma_C))\) such that \(\tau(s)\geq 0\), we will use the shorthand notation:
	\[
		\int_{[0,S]}\int_{\Gamma_C} f\tau(s)\,\bigl|\mathbf{v}_{\rm t}\bigr| := \int_{[0,S]}\int_{\Gamma_C} f\tau(s)\,\bigl|\bar{\mathbf{v}}_{\rm t}(s)\bigr|\,{\rm d}x\,{\rm d}\|\mathbf{v}\|,
	\]
	where, for \(\mathbf{w}\in H^{1,2}(\Gamma_C,\mathbb{R}^N)\), the function \(|\mathbf{w}|\in H^{1/2}(\Gamma_C,[0,\infty))\) is defined as \(|\mathbf{w}|(x):=|\mathbf{w}(x)|\).
\end{defi}

Straightforward application of the theory of sweeping processes~\cite{Moreau} then yields the following result.

\begin{theo}[Moreau]
	\label{theo:ExistenceSweeping}
	Let \(\mathbf{u}_0\in C(0)\), and functions \(\mathbf{F}\in \textsl{BV}^+([0,S];L^2(\Omega;\mathbb{R}^N))\), \(\mathbf{T}\in \textsl{BV}^+([0,S];L^2(\Gamma_T;\mathbb{R}^N))\), \(f\in W^{1,\infty}(\Gamma_C;\left[0,+\infty\right[)\) and \(\tau\in \textsl{BV}^+([0,S];H^{-1/2}(\Gamma_C))\) with \(\tau(s)\leq 0\), for all \(s\in [0,S]\), be given.  Then, there exists a unique function \(\mathbf{u}\in \textsl{BV}^+([0,S];V)\) such that \(\mathbf{u}(0)=\mathbf{u}_0\) and, for all \(\mathbf{v}\in \mathscr{M}([0,S];V)\):
	\begin{multline*}
		\int_{[0,S]}\int_{\Omega} \boldsymbol{\varepsilon}(\mathbf{u}(s)):\boldsymbol{\Lambda}\boldsymbol{\varepsilon}\bigl(\mathbf{v}-\mathbf{\dot{u}}\bigr) \geq \int_{[0,S]}\int_{\Omega} \mathbf{F}(s)\cdot\bigl(\mathbf{v}-\mathbf{\dot{u}}\bigr) + \int_{[0,S]}\int_{\Gamma_T}\mathbf{T}(s)\cdot\bigl(\mathbf{v}-\mathbf{\dot{u}}\bigr) +\mbox{}\\\mbox{}+ \int_{[0,S]}\int_{\Gamma_C}\tau(s)\,\Bigl[v_{\rm n}-\dot{u}_{\rm n}\Bigr] +\int_{[0,S]}\int_{\Gamma_C} f\tau(s)\,\Bigl[\bigl|\mathbf{v}_{\rm t}\bigr|-|\mathbf{\dot{u}}_{\rm t}\bigr|\Bigr].
	\end{multline*}
	In the case where the functions \(\mathbf{F}\), \(\mathbf{T}\) and \(\tau\) are absolutely continuous in time, the solution \(\mathbf{u}\) is also absolutely continuous in time and satisfies, for all \(\mathbf{v}\in V\) and almost all \(s\in [0,S]\):
	\begin{multline*}
		\int_{\Omega} \boldsymbol{\varepsilon}(\mathbf{u}(s)):\boldsymbol{\Lambda}\boldsymbol{\varepsilon}\bigl(\mathbf{v}-\mathbf{\dot{u}}(s)\bigr) \geq \int_{\Omega} \mathbf{F}(s)\cdot\bigl(\mathbf{v}-\mathbf{\dot{u}}(s)\bigr) + \int_{\Gamma_T}\mathbf{T}(s)\cdot\bigl(\mathbf{v}-\mathbf{\dot{u}}(s)\bigr) +\mbox{}\\\mbox{}+ \int_{\Gamma_C}\tau(s)\,\Bigl[v_{\rm n}-\dot{u}_{\rm n}(s)\Bigr] + \int_{\Gamma_C} f\tau(s)\,\Bigl[\bigl|\mathbf{v}_{\rm t}\bigr|-|\mathbf{\dot{u}}_{\rm t}(s)\bigr|\Bigr].
	\end{multline*}
\end{theo}

Theorem~\ref{theo:ExistenceSweeping} will not be used in the sequel, but it is nonetheless useful, as it reveals that Subproblem~1 in Section~\ref{sec:overviewQuasiStatic} exhibits the structure of a sweeping process. This structural insight is crucial for formulating Subproblem~1 correctly and precisely within the framework of right-continuous functions with bounded variation. Adopting this setting is essential, since the solution to the quasi-static problem may exhibit spontaneous jumps, as will be shown later on.

\subsection{The quasi-static problem and statement of our solvability result}

\label{sec:quasiStaticProblem}
Using the content of the previous subsection, we are now able to define precisely what is meant by a \emph{possibly jumping-in-time} solution of the quasi-static problem~\eqref{eq:SignoriniCoulombCont}. 

First, we have to define the admissible initial conditions.  Considering the nonempty closed convex subset \(K_g\) of the Hilbert space \(V\) (see formula~\eqref{eq:defV}) defined by:
\[
K_g := \Bigl\{\mathbf{u}\in V \bigm| u_{\rm n}\leq g\text{ on }\Gamma_C\Bigr\},
\]
\(\mathbf{F}\in L^2(\Omega;\mathbb{R}^N)\), \(\mathbf{T}\in L^2(\Gamma_T;\mathbb{R}^N)\) and \(\mathbf{t}_{\rm t}\in H^{-1/2}(\Gamma_C;\mathbb{R}^N)\), such that \(\mathbf{t}_{\rm t}\cdot\mathbf{n}=0\) on \(\Gamma_C\) (\(\mathbf{t}_{\rm t}\) is a given tangent traction on \(\Gamma_C\)), we denote by \(\mathbf{u}=\mathscr{I}\bigl(\mathbf{F},\mathbf{T},\mathbf{t}_{\rm t}\bigr)\) the unique solution \(\mathbf{u}\in K_g\) of the variational inequality:
\[
	\forall \mathbf{v}\in K_g,\qquad\int_{\Omega} \boldsymbol{\varepsilon}(\mathbf{u}):\boldsymbol{\Lambda}\boldsymbol{\varepsilon}(\mathbf{v}-\mathbf{u}) \,{\rm d}x\geq \int_{\Omega} \mathbf{F}\cdot\bigl(\mathbf{v}-\mathbf{u}\bigr) + \int_{\Gamma_T}\mathbf{T}\cdot\bigl(\mathbf{v}-\mathbf{u}\bigr) + \bigl\langle\mathbf{t}_{\rm t},\mathbf{v}_{\rm t}-\mathbf{u}_{\rm t}\bigr\rangle.
\]
Then, we set:
\begin{multline}
	\mathscr{A}(\mathbf{F},\mathbf{T},f) := \biggl\{\mathbf{u}=\mathscr{I}\bigl(\mathbf{F},\mathbf{T},\mathbf{t}_{\rm t}\bigr) \Bigm| \mathbf{t}_{\rm t}\in H^{-1/2}(\Gamma_C;\mathbb{R}^N) \;\ \text{such that:}\\ \mathbf{t}_{\rm t}\cdot\mathbf{n}=0\quad \text{and}\quad \forall\mathbf{v}\in H^{1/2}(\Gamma_C;\mathbb{R}^N),\quad \bigl\langle \mathbf{t}_{\rm t},\mathbf{v}_{\rm t}\bigr\rangle-\bigl\langle f\,\boldsymbol{\sigma}(\mathbf{u})\mathbf{n},|\mathbf{v}_{\rm t}|\mathbf{n}\bigr\rangle \geq 0\biggr\},\label{eq:defCf}
\end{multline}
which is nonempty as it contains \(\mathscr{I}(\mathbf{F},\mathbf{T},\mathbf{0})\).  Equivalently, \(\mathscr{A}(\mathbf{F},\mathbf{T},f)\) is the set of those \(\mathbf{u}\in V\) such that there exists \(t_{\rm n}\in H^{-1/2}(\Gamma_C)\), \(t_{\rm n}\leq 0\) such that:
\begin{align}
	& \bullet\quad 	\forall \mathbf{v}\in V,\qquad \int_{\Omega} \boldsymbol{\varepsilon}(\mathbf{u}):\boldsymbol{\Lambda}\boldsymbol{\varepsilon}(\mathbf{v}) \geq \int_{\Omega} \mathbf{F}\cdot\mathbf{v} + \int_{\Gamma_T}\mathbf{T}\cdot\mathbf{v} + \int_{\Gamma_C}t_{\rm n}\,v_{\rm n}+ft_{\rm n}\,|\mathbf{v}_{\rm t}|,\label{eq:admissible1}\\
	& \bullet\quad 	 	\forall t\in H^{-1/2}(\Gamma_C), \;\text{ with }\;t\leq 0,\qquad \int_{\Gamma_C}\bigl(u_{{\rm n}}-g\bigr)\bigl(t-t_{\rm n}\bigr)\geq 0.\label{eq:admissible2}
\end{align}
Taking \(\mathbf{v}=\lambda\mathbf{v}\) (\(\lambda>0\)) in~\eqref{eq:optimalityIncrementalPb}, dividing by \(\lambda\) and taking the limit as \(\lambda\to +\infty\), we see that all the solutions of the incremental problem for arbitrary \(\mathbf{w}_{\rm t}\in H^{1/2}(\Gamma_C;\mathbb{R}^N)\) (\(\mathbf{w}_{\rm t}\cdot\mathbf{n}=0\)) are in \(\mathscr{A}(\mathbf{F},\mathbf{T},f)\).  

Recalling the notation $\mathbf{t}:= \boldsymbol\sigma(\mathbf{u})\,\mathbf{n}=t_{\rm n}\mathbf{n}+\mathbf{t}_{\rm t}$, \(\mathscr{A}(\mathbf{F},\mathbf{T},f)\) is also the set of those \(\mathbf{u}\in H^1(\Omega;\mathbb{R}^N)\) satisfying:
\[
\left\{\quad
	\begin{aligned}
		& \text{div}\,\boldsymbol{\sigma}(\mathbf{u}) + \mathbf{F} = \mathbf{0}, \qquad & & \text{ in }\Omega,\\
		& \mathbf{u} = \mathbf{0},\qquad\qquad & & \text{ on }\Gamma_U,\\
		& \boldsymbol\sigma(\mathbf{u})\,\mathbf{n} = \mathbf{T},\qquad & & \text{ on }\Gamma_T,\\
		& u_{\rm n}\leq g,\quad t_{\rm n} \leq 0,\quad (u_{\rm n}-g)t_{\rm n}=0, \qquad & & \text{ on }\Gamma_C,\\
		& \forall \mathbf{v}\in H^{1/2}(\Gamma_C;\mathbb{R}^N),\qquad \bigl\langle\mathbf{t}_{\rm t},\mathbf{v}_{\rm t}\bigr\rangle - \bigl\langle ft_{\rm n} , |\mathbf{v}_{\rm t}|\bigr\rangle \geq 0, \qquad & & \text{ on }\Gamma_C.
	\end{aligned}
\right.
\]
\mbox{}

\begin{defi}[Quasi-static problem]
	\label{def:quasi-static}
	Consider arbitrary \(\mathbf{F}\in \textsl{BV}^+([0,S];L^2(\Omega;\mathbb{R}^N))\), \(\mathbf{T}\in \textsl{BV}^+([0,S];L^2(\Gamma_T;\mathbb{R}^N))\), \(f\in W^{1,\infty}(\Gamma_C;\left[0,+\infty\right[)\) and an initial condition \(\mathbf{u}_0\in \mathscr{A}(\mathbf{F}(0),\mathbf{T}(0),f)\). We will say that \(\mathbf{u}\in \textsl{BV}^+([0,S];V)\) is a solution of the quasi-static problem, if \(\mathbf{u}(0)=\mathbf{u}_0\) and, there exists a function \(t_n\in \textsl{BV}^+([0,S];H^{-1/2}(\Gamma_C))\) with \(t_{\rm n}(s)\leq 0\), for all \(s\in [0,S]\), and such that:
		\begin{multline*}\bullet\;\forall \mathbf{v}\in \mathscr{M}([0,S];V),\qquad \int_{[0,S]}\int_{\Omega} \boldsymbol{\varepsilon}(\mathbf{u}(s)):\boldsymbol{\Lambda}\boldsymbol{\varepsilon}\bigl(\mathbf{v}-\mathbf{\dot{u}}\bigr) \geq \int_{[0,S]}\int_{\Omega} \mathbf{F}(s)\cdot\bigl(\mathbf{v}-\mathbf{\dot{u}}\bigr) +\mbox{}\\\mbox{}+ \int_{[0,S]}\int_{\Gamma_T}\mathbf{T}(s)\cdot\bigl(\mathbf{v}-\mathbf{\dot{u}}\bigr) + \int_{[0,S]}\int_{\Gamma_C}t_{\rm n}(s)\,\Bigl[v_{\rm n}-\dot{u}_{\rm n}\Bigr] + \int_{[0,S]}\int_{\Gamma_C}ft_{\rm n}(s)\,\Bigl[\bigl|\mathbf{v}_{\rm t}\bigr|-|\mathbf{\dot{u}}_{\rm t}\bigr|\Bigr],
	\end{multline*}	
	\(\displaystyle\hspace*{0.45cm}\bullet\;\:\forall s\in [0,S],\quad\forall t\in H^{-1/2}(\Gamma_C), \;\text{ with }\;t\leq 0,\qquad \int_{\Gamma_C}\bigl(u_{\rm n}(s)-g\bigr)\bigl(t-t_{\rm n}(s)\bigr)\geq 0,\)

	\medskip
	\noindent{}where the precise meaning of the integrals in the first condition is given by Definition~\ref{defi:intMeasure}.
\end{defi}

Note that the first condition entails that \(t_{\rm n}(s)\) is the normal part of the surface traction \(\boldsymbol\sigma(\mathbf{u}(s))\mathbf{n}\) on \(\Gamma_C\), for all \(s\in [0,S]\).  

As functions \(\mathbf{u}\in \textsl{BV}^+([0,S];V)\) may contain jumps in time \(\mathbf{u}(s+)\neq \mathbf{u}(s-)\), this definition clearly accounts for the possibility of jumps in time of the solution \(\mathbf{u}\) of the quasi-static problem.  As far as we are aware, Definition~\ref{def:quasi-static} is the first definition to explicitly account for a possibly jumping solution of the quasi-static problem~\eqref{eq:SignoriniCoulombCont}.  As will be shown in the sequel, there may exist solutions \(\mathbf{u}\) in the sense of Definition~\ref{def:quasi-static} of the quasi-static problem that exhibit jumps in time, while the data \(\mathbf{F}\), \(\mathbf{T}\) are absolutely continuous (or even Lipschitz-continuous) in time.  This fact, although already informally stated~\cite{Klarbring} in the engineering literature, seems to be formalized here for the first time.

We now state the minimal condition on the friction coefficient \(f\) for which we are going to prove the existence of an absolutely continuous solution of the quasi-static problem in the sense of Definition~\ref{def:quasi-static}, for any absolutely continuous data \(\mathbf{F}\), \(\mathbf{T}\).  The only ground on which it is based is that it is the minimal condition required by our proof of existence, as it will be shown in Section~\ref{sec:ProofExistence}. However, we will construct a simple example in Section~\ref{sec:Examples} of a system where a solution of the quasi-static problem \emph{must jump in time}, while the data \(\mathbf{F}\), \(\mathbf{T}\) are Lipschitz-continuous in time and whenever the aforementioned condition is broken.

\begin{defi}
	\label{def:conditionC}
	Given \(f\in W^{1,\infty}(\Gamma_C;\left[0,+\infty\right[)\) and the data \(\mathbf{F}\in \textsl{BV}^+([0,S];L^2(\Omega;\mathbb{R}^N))\), \(\mathbf{T}\in \textsl{BV}^+([0,S];L^2(\Gamma_T;\mathbb{R}^N))\), we set:
	\begin{align*}
		\bar{\mathscr{A}}(\mathbf{F},\mathbf{T},f) & := \bigcup_{s\in [0,S]}\mathscr{A}(\mathbf{F}(s),\mathbf{T}(s),f),\\	
		\tilde{\mathscr{A}}(f) & := \bigcup_{\substack{\mathbf{F}\in L^2(\Omega;\mathbb{R}^N)),\\
		\mathbf{T}\in L^2(\Gamma_T;\mathbb{R}^N)),\\
		}}\mathscr{A}(\mathbf{F},\mathbf{T},f).
	\end{align*}
	A friction coefficient \(f\in W^{1,\infty}(\Gamma_C;\left[0,+\infty\right[)\) is said to fulfil condition~$\mathscr{C}$ for the data \(\mathbf{F}\) and \(\mathbf{T}\), if:
		\[
			\sup_{\substack{\mathbf{u},\mathbf{v}\in \bar{\mathscr{A}}(\mathbf{F},\mathbf{T},f),\\
		\mathbf{u}\neq \mathbf{v},\\
		}}\frac{\bigl\langle f\,\boldsymbol{\sigma}(\mathbf{u}-\mathbf{v})\mathbf{n}\;,\;|\mathbf{u}_{\rm t}-\mathbf{v}_{\rm t}|\mathbf{n} \bigr\rangle_{\Gamma_C}}{\|\mathbf{u}-\mathbf{v}\|_{a}\Bigl\{\|\mathbf{u}-\mathbf{v}\|_{a} + \bigl\|\text{\rm div}\,\boldsymbol{\sigma}(\mathbf{u}-\mathbf{v})\bigr\|_{L^2(\Omega)} + \bigl\|\boldsymbol{\sigma}(\mathbf{u}-\mathbf{v})\mathbf{n}\bigr\|_{L^2(\Gamma_T)}\Bigr\}}<1,
		\]
		where \(\|\mathbf{u}-\mathbf{v}\|_{a}:=\bigl(a(\mathbf{u}-\mathbf{v},\mathbf{u}-\mathbf{v})\bigr)^{1/2}\) is a shorthand notation.

		A friction coefficient \(f\in W^{1,\infty}(\Gamma_C;\left[0,+\infty\right[)\) is said to fulfil uniform condition~$\mathscr{C}$ if it fulfills condition~$\mathscr{C}$ for all \(\mathbf{F}\in L^2(\Omega;\mathbb{R}^N)\), \(\mathbf{T}\in L^2(\Gamma_T;\mathbb{R}^N)\), that is:
		\[
			\sup_{\substack{\mathbf{u},\mathbf{v}\in \tilde{\mathscr{A}}(f),\\
		\mathbf{u}\neq \mathbf{v},\\
		}}\frac{\bigl\langle f\,\boldsymbol{\sigma}(\mathbf{u}-\mathbf{v})\mathbf{n}\;,\;|\mathbf{u}_{\rm t}-\mathbf{v}_{\rm t}|\mathbf{n} \bigr\rangle_{\Gamma_C}}{\|\mathbf{u}-\mathbf{v}\|_{a}\Bigl\{\|\mathbf{u}-\mathbf{v}\|_{a} + \bigl\|\text{\rm div}\,\boldsymbol{\sigma}(\mathbf{u}-\mathbf{v})\bigr\|_{L^2(\Omega)} + \bigl\|\boldsymbol{\sigma}(\mathbf{u}-\mathbf{v})\mathbf{n}\bigr\|_{L^2(\Gamma_T)}\Bigr\}}<1,
		\]
\end{defi}

It is easy to see that if \(f\) fulfills condition~$\mathscr{C}$ for the data \(\mathbf{F}\) and \(\mathbf{T}\) (respectively uniform condition~$\mathscr{C}$), then the same is true for \(\alpha f\), for all \(\alpha\in [0,1]\).  Uniform condition~$\mathscr{C}$ is a condition on the friction coefficient \(f\), the elastic moduli \(\boldsymbol{\Lambda}\) and the geometry \(\Omega\), \(\Gamma_U\), \(\Gamma_T\), \(\Gamma_C\), \(g\) only.  By Lemma~\ref{lemm:condAndersson}, every friction coefficient \(f\in W^{1,\infty}(\Gamma_C;\left[0,+\infty\right[)\) such that:
\begin{equation}
	\label{eq:conAndersson}
		\|f\|_{W^{1,\infty}(\Gamma_C)}<\frac{\sqrt{2}c_0}{2\max\{1,C_0\}\|\text{\rm ext}\|\|\text{\rm tr}\|\sqrt{\max\{1,G(\Gamma_C)\}}}
\end{equation}
fulfills uniform condition~$\mathscr{C}$.  Hence, there is always a range of infinitely many friction coefficients satisfying uniform condition~$\mathscr{C}$.  Condition~\eqref{eq:conAndersson} first appeared in~\cite{Andersson} in close association with the nonphysical Eck-Jarušek condition coming from Nirenberg's shift technique.

Condition~$\mathscr{C}$ is designed to be the minimal condition ensuring the following corollary.

\begin{prop}  
	\label{prop:estimateConditionC}
	Let \(\mathbf{F}\in \textsl{BV}^+([0,S];L^2(\Omega;\mathbb{R}^N))\), \(\mathbf{T}\in \textsl{BV}^+([0,S];L^2(\Gamma_T;\mathbb{R}^N))\).   Let \(f\in W^{1,\infty}(\Gamma_C;\left[0,+\infty\right[)\) fulfilling condition~$\mathscr{C}$ for the data \(\mathbf{F}\) and \(\mathbf{T}\) (respectively uniform condition $\mathscr{C}$). Let \(s_1\in [0,S]\) and \(\mathbf{u}_1\in \mathscr{A}(\mathbf{F}(s_1),\mathbf{T}(s_1),f)\). Let
	\(s_2\in [0,S]\) and \(\mathbf{u}_2\in V\) be a solution of the incremental problem (see Definition~\ref{def:IncrementalProblem}) for the data \(\mathbf{F}(s_2)\), \(\mathbf{T}(s_2)\) and \(\mathbf{w}_{\rm t}=\mathbf{u}_{1,{\rm t}}\).  Then,
	\[
		\bigl\| \mathbf{u}_2-\mathbf{u}_1\bigr\|_{H^1} \leq C\Bigl(\bigl\|\mathbf{F}(s_2)-\mathbf{F}(s_1)\bigr\|_{L^2(\Omega)} + \bigl\|\mathbf{T}(s_2)-\mathbf{T}(s_1)\bigr\|_{L^2(\Gamma_T)}\Bigr),
	\]
	for some \emph{positive} constant \(C\) depending only on \(f\), \(\boldsymbol{\Lambda}\), \(\Omega\), \(\Gamma_U\), \(\Gamma_T\), \(\Gamma_C\), \(g\), and also possibly on \(\mathbf{F}\), \(\mathbf{T}\).
\end{prop}

\noindent\textbf{Proof.} Chosing \(\mathbf{v}=\mathbf{u}_2-\mathbf{u}_1\) in the variational inequality~\eqref{eq:admissible1} for \(\mathbf{u}_1\) and \(\mathbf{v}=\mathbf{u}_1\) in the optimality condition~\eqref{eq:optimalityIncrementalPb} for \(\mathbf{u}_2\), and taking the sum of the corresponding two inequalities, we obtain:
\begin{multline*}
	\int_{\Omega} \boldsymbol{\varepsilon}(\mathbf{u}_2-\mathbf{u}_1):\boldsymbol{\Lambda}\boldsymbol{\varepsilon}(\mathbf{u}_2-\mathbf{u}_1) - \int_{\Gamma_C}f\bigl(t_{2,{\rm n}}-t_{1,{\rm n}}\bigr)\bigl|\mathbf{u}_{2,{\rm t}}-\mathbf{u}_{1,{\rm t}}\bigr|\leq \mbox{}\\ \int_{\Omega} \bigl(\mathbf{F}(s_2)-\mathbf{F}(s_1)\bigr)\cdot\bigl(\mathbf{u}_2-\mathbf{u}_1\bigr) + \int_{\Gamma_T}\bigl(\mathbf{T}(s_2)-\mathbf{T}(s_1)\bigr)\cdot\bigl(\mathbf{u}_2-\mathbf{u}_1\bigr) + \int_{\Gamma_C}\bigl(t_{2,{\rm n}}-t_{1,{\rm n}}\bigr)\bigl(u_{2,{\rm n}}-u_{1,{\rm n}}\bigr).
\end{multline*}
Choosing \(t=t_{2,{\rm n}}\) in the variational inequality~\eqref{eq:admissible2} for \(\mathbf{u}_1\) and \(t=t_{1,{\rm n}}\) in that~\eqref{eq:SignoriniIncrementalPb} for \(\mathbf{u}_2\), and taking the sum of the corresponding two inequalities, we obtain:
\[
	\int_{\Gamma_C}\bigl(t_{2,{\rm n}}-t_{1,{\rm n}}\bigr)\bigl(u_{2,{\rm n}}-u_{1,{\rm n}}\bigr)\leq 0,
\]
hence:
\begin{multline}
	\int_{\Omega} \boldsymbol{\varepsilon}(\mathbf{u}_2-\mathbf{u}_1):\boldsymbol{\Lambda}\boldsymbol{\varepsilon}(\mathbf{u}_2-\mathbf{u}_1) - \int_{\Gamma_C}f\bigl(t_{2,{\rm n}}-t_{1,{\rm n}}\bigr)\,\bigl|\mathbf{u}_{2,{\rm t}}-\mathbf{u}_{1,{\rm t}}\bigr| \leq \mbox{}\\ \max\bigl\{1,\|\text{\rm tr}\|\bigr\}\bigl\| \mathbf{u}_2-\mathbf{u}_1\bigr\|_{H^1}\Bigl(\bigl\|\mathbf{F}(s_2)-\mathbf{F}(s_1)\bigr\|_{L^2(\Omega)} + \bigl\|\mathbf{T}(s_2)-\mathbf{T}(s_1)\bigr\|_{L^2(\Gamma_T)}\Bigr).\label{eq:claimPart1}
\end{multline}
Now, condition~$\mathscr{C}$ ensures that:
\begin{multline}
	\bigl\| \mathbf{u}_2-\mathbf{u}_1\bigr\|_{a}^2 - \Bigl\langle f\,\bigl(t_{2,{\rm n}}-t_{1,{\rm n}}\bigr)\;,\;|\mathbf{u}_{2,{\rm t}}-\mathbf{u}_{1,{\rm t}}|\Bigr\rangle_{\Gamma_C}\geq \alpha\bigl\| \mathbf{u}_2-\mathbf{u}_1\bigr\|_{a}^2\\\mbox{}-\beta\bigl\| \mathbf{u}_2-\mathbf{u}_1\bigr\|_{a}\Bigl(\bigl\|\mathbf{F}(s_2)-\mathbf{F}(s_1)\bigr\|_{L^2(\Omega)} + \bigl\|\mathbf{T}(s_2)-\mathbf{T}(s_1)\bigr\|_{L^2(\Gamma_T)}\Bigr),\label{eq:claimPart2}
\end{multline} 
where we denote the supremum in condition~$\mathscr{C}$ by \(1-\alpha\), so that \(\alpha>0\).  It is now sufficient to gather inequalities~\eqref{eq:claimPart1} and~\eqref{eq:claimPart2} using~\eqref{eq:Korn} to get the claim. \qed

\bigskip
It's worth expressing the conclusion of Proposition~\ref{prop:estimateConditionC} in words.    Condition~$\mathscr{C}$ ensures that if two consecutive values \(\mathbf{F}(s_1)\) and \(\mathbf{F}(s_2)\) of the load are close in a time stepping discretization of the quasi-static problem, then the corresponding successive values \(\mathbf{u}_1\) and \(\mathbf{u}_2\) of the solution obtained by solving the incremental problem are also close.  Hence, condition~$\mathscr{C}$ is meant to ensure that if the load has no jumps in time, then the solution of the quasi-static problem has no jumps in time either, and it is designed as being the minimal condition ensuring this property.  Note that Proposition~\ref{prop:estimateConditionC} yields no information on the uniqueness of solution for the incremental problem.

We are now ready to state our existence theorem for the quasi-static problem as formulated in Definition~\ref{def:quasi-static}, under condition~$\mathscr{C}$.  

\begin{theo} 
	\label{theo:mainExistenceTheorem}
	We suppose \(N=2\) and \(\boldsymbol{\Lambda}\in W^{1,\infty}(\Omega)\) satisfying condition~\eqref{eq:ConditionFrictionAnisotropic} (in particular, no condition if \(\boldsymbol{\Lambda}\) is isotropic at each \(x\in \Omega\)). 
	Consider \(\mathbf{F}\in \textsl{BV}^+([0,S];L^2(\Omega;\mathbb{R}^2))\) and \(\mathbf{T}\in \textsl{BV}^+([0,S];L^2(\Gamma_T;\mathbb{R}^2))\). Let \(f\in W^{1,\infty}(\Gamma_C;\left[0,+\infty\right[)\) fulfilling condition~$\mathscr{C}$ for the data \(\mathbf{F}\) and \(\mathbf{T}\).  Let \(\mathbf{u}_0\in \mathscr{A}(\mathbf{F}(0),\mathbf{T}(0),f)\) be an admissible initial condition.  Then, there exists a solution \(\mathbf{u}\in \textsl{BV}^+([0,S];V)\) of the quasi-static problem~\eqref{eq:SignoriniCoulombCont}, in the sense of Definition~\ref{def:quasi-static}.  In the case where \(\mathbf{F}\), \(\mathbf{T}\) are absolutely continuous in time, the solution \(\mathbf{u}\) is also absolutely continuous in time and satisfies:
	\begin{multline*}\bullet\;\text{for a.a. }s\in [0,S],\quad\forall \mathbf{v}\in V,\qquad \int_{\Omega} \boldsymbol{\varepsilon}(\mathbf{u}(s)):\boldsymbol{\Lambda}\boldsymbol{\varepsilon}\bigl(\mathbf{v}-\mathbf{\dot{u}}(s)\bigr)\,{\rm d}x \geq \int_{\Omega} \mathbf{F}(s)\cdot\bigl(\mathbf{v}-\mathbf{\dot{u}}(s)\bigr)\,{\rm d}x +\mbox{}\\\mbox{}+ \int_{\Gamma_T}\mathbf{T}(s)\cdot\bigl(\mathbf{v}-\mathbf{\dot{u}}(s)\bigr) + \int_{\Gamma_C}t_{\rm n}(s)\,\Bigl[v_{\rm n}-\dot{u}_{\rm n}(s)\Bigr] + \int_{\Gamma_C}ft_{\rm n}(s)\,\Bigl[\bigl|\mathbf{v}_{\rm t}\bigr|-|\mathbf{\dot{u}}_{\rm t}(s)\bigr|\Bigr],
	\end{multline*}	
	\(\displaystyle\hspace*{0.5cm}\bullet\;\:\forall s\in [0,S],\quad\forall t\in H^{-1/2}(\Gamma_C), \;\text{ with }\;t\leq 0,\qquad \int_{\Gamma_C}\bigl(u_{\rm n}(s)-g\bigr)\bigl(t-t_{\rm n}(s)\bigr)\geq 0.\)
\end{theo}

The detailed proof of this theorem is postponed to Section~\ref{sec:ProofExistence}.  It relies on our proof of existence of solution for the 2D incremental problem in~\cite{BI-ma3as}, the compensated-compactness-like fine property of the 2D elastic Neumann-to-Dirichlet operator proved in~\cite{BI-ma3as}, and some pseudomonotonicity-like arguments. 

Before focusing on the details of the proof, we are going to discuss examples and present, in particular, a striking example showing that when condition~$\mathscr{C}$ is violated, absolutely continuous loads can induce solutions having spontaneous jumps in time, and that there may even not exist continuous solutions in that case.  It is likely, though, that there always exists a quasi-static solution with bounded variation, whenever the incremental problem is solvable, even if condition~$\mathscr{C}$ is not satisfied. This will not be proved in this paper.  On one hand, jumping solutions have little physical relevance in the quasi-static setting but, on the other hand, having an unconditional existence result of a solution with bounded variation would be a prerequisite to classify those absolutely continuous loads that induce jumping solutions and those that do not, in the case of large friction coefficients that violate condition~$\mathscr{C}$.  Our result only states that under condition~$\mathscr{C}$, there always exists an absolutely continuous solution for any absolutely continuous load.

\section{Various examples}
\label{sec:Examples}

\subsection{Case of the spatially discrete system in $\mathbb{R}^2$}
\label{sec:discrete}

We consider a point particle in $\mathbb{R}^2$ evolving in a quadratic elastic potential under the action of a prescribed force \(\mathbf{F}(s)=(F_{\rm n}(s),F_{\rm t}(s))\).  The quadratic elastic potential is defined by a positive definite symmetric stiffness matrix:
\[
	\mathbf{K}=\begin{pmatrix}
		k_{\rm nn} & k_{\rm nt}\\
		k_{\rm nt} & k_{\rm tt}
	\end{pmatrix}.
\]
Possibly reversing the tangential direction, if necessary, we can assume that \(k_{\rm nt}\geq 0\) without restricting the generality. The location of the particle is denoted by \(\mathbf{u}(s)=(u_{\rm n}(s),u_{\rm t}(s))\).  The particle is in contact with a rigid wall at \(u_{\rm n}=0\) and the contact is frictional, with a friction coefficient \(f\geq 0\).  The contact force is denoted by \(\mathbf{t}(s)=(t_{\rm n}(s),t_{\rm t}(s))\).

A given initial condition \(\mathbf{u}_0\in\mathbb{R}^2\) is said to be admissible with respect to \(\mathbf{F}\in\mathbb{R}^2\) if:
	\begin{itemize}
		\item \(\mathbf{K}\mathbf{u}_0=\mathbf{F}+\mathbf{t}\),
		\item \(u_{0{\rm n}}\leq 0,\quad t_{\rm n}\leq 0, \quad u_{0{\rm n}}\,t_{\rm n}=0\),
		\item \(|t_{\rm t}|\leq -ft_{\rm n}.\)
	\end{itemize}

\begin{defi}[Discrete quasi-static problem]
	\label{def:quasi-static-discrete}
	Let \(\mathbf{F}\in\textsl{BV}^+([0,S];\mathbb{R}^2)\) be a given external force, and \(\mathbf{u}_0\in\mathbb{R}^2\) some initial condition supposed to be admissible with respect to \(\mathbf{F}(0)\).  We will say that \(\mathbf{u},\mathbf{t}\in \textsl{BV}^+([0,S];\mathbb{R}^2)\) is a solution of the quasi-static problem, if:
	\begin{itemize}
		\item \(\mathbf{u}(0)=\mathbf{u}_0\),
		\item \(\forall s\in [0,S],\quad \mathbf{K}\mathbf{u}(s)=\mathbf{F}(s)+\mathbf{t}(s)\),
		\item \(\forall s\in [0,S],\quad u_{\rm n}(s)\leq 0,\quad t_{\rm n}(s)\leq 0, \quad u_{\rm n}(s)\,t_{\rm n}(s)=0\),
		\item \(\forall s\in [0,S],\quad |t_{\rm t}(s)|\leq -ft_{\rm n}(s),\quad \text{ and }\quad t_{\rm t}\dot{u}_{\rm t}=ft_{\rm n}|\dot{u}_{\rm t}|\).
	\end{itemize}
\end{defi}

Above \(|\dot{u}_{\rm t}|\) stands for the total variation measure of the measure \(\dot{u}_{\rm t}\).

\begin{defi}[Discrete incremental problem]
	\label{def:incremental-discrete}
	Let \(\mathbf{F}\in\mathbb{R}^2\) be a given external force, and \(w_t\in\mathbb{R}\).  We will say that \(\mathbf{u},\mathbf{t}\in \mathbb{R}^2\) is a solution of the incremental problem, if:
	\begin{itemize}
		\item \(\mathbf{K}\mathbf{u}=\mathbf{F}+\mathbf{t}\),
		\item \(u_{\rm n}\leq 0,\quad t_{\rm n}\leq 0, \quad u_{\rm n}t_{\rm n}=0\),
		\item \(|t_{\rm t}|\leq -ft_{\rm n},\quad \text{ and }\quad t_{\rm t}\bigl(u_{\rm t}-w_{\rm t}\bigr)=ft_{\rm n}|u_{\rm t}-w_{\rm t}|\).
	\end{itemize}
\end{defi}

Above, \(|u_{\rm t}-w_{\rm t}|\) now stands for the absolute value of the real number \(u_{\rm t}-w_{\rm t}\).

It will sometimes be convenient to view the discrete incremental problem as a fixed point problem. Set:
\[
K:= \Bigl\{\mathbf{u}\in\mathbb{R}^2\bigm| u_{\rm n}\leq 0\Bigr\}.
\]
Pick an arbitrary \(\sigma\geq 0\) and consider the minimum problem:
\begin{equation}
	\label{eq:discreteMinProb}
		\mathbf{u}=\argmin_{\mathbf{v}\in K}\Bigl\{\frac{1}{2}\mathbf{v}\cdot\mathbf{K}\mathbf{v}-\mathbf{F}\cdot\mathbf{v} + \sigma\bigl|v_{\rm t}-w_{\rm t}\bigr|\Bigr\},
\end{equation}
whose unique minimizer is characterized either by the optimality inequality:
\begin{equation}
	\label{eq:discreteOptIneq}
		\forall \mathbf{v}\in K,\qquad \mathbf{u}\cdot\mathbf{K}(\mathbf{v} - \mathbf{u}) \geq \mathbf{F}\cdot(\mathbf{v} - \mathbf{u}) - \sigma\Bigl(\bigl|v_{\rm t}-w_{\rm t}\bigr| - \bigl|u_{\rm t}-w_{\rm t}\bigr|\Bigr),
\end{equation}
or by the conditions:
\begin{itemize}
	\item \(\mathbf{K}\mathbf{u}=\mathbf{F}+\mathbf{t}\),
	\item \(u_{\rm n}\leq 0,\quad t_{\rm n}\leq 0, \quad u_{\rm n}t_{\rm n}=0\),
	\item \(|t_{\rm t}|\leq \sigma,\quad \text{ and }\quad t_{\rm t}\bigl(u_{\rm t}-w_{\rm t}\bigr)=-\sigma|u_{\rm t}-w_{\rm t}|\).
\end{itemize}
Based on this uniquely solvable minimum problem, one can define:
\begin{equation}
	\label{eq:defP}
	P\left\{\begin{array}{rcl}\left[0,+\infty\right[ & \to & \left[0,+\infty\right[\\
		\sigma & \mapsto & -ft_{\rm n}\end{array}\right.
\end{equation}
so that any fixed point of \(P\) is a solution of the discrete incremental problem, and reciprocally.

\begin{prop}
	The discrete incremental problem has always a solution.
\end{prop}

\noindent\textbf{Proof.} Taking \(\sigma_1,\sigma_2\in \left[0,+\infty\right[\) and denoting by \(\mathbf{u}_1,\mathbf{u}_2\in K\) the corresponding minimizers~\eqref{eq:discreteMinProb}, the sum of the corresponding optimality inequalities~\eqref{eq:discreteOptIneq} yields:
\begin{equation}
	\label{eq:discreteOptIneqSum}
		(\mathbf{u}_1-\mathbf{u}_2)\cdot\mathbf{K}(\mathbf{u}_1-\mathbf{u}_2) \leq |\sigma_1-\sigma_2|\bigl|u_{1{\rm t}}-u_{2{\rm t}}\bigr|,
\end{equation}
which shows that the mapping \(P\) defined by~\eqref{eq:defP} is Lipschitz-continuous.

Denote:
\[
	\langle x\rangle^+:=\max\{x,0\}\quad\text{ and }\quad \langle x\rangle^-:=\max\{-x,0\},
\]
and pick:
\[
	\sigma\geq \Sigma_1:=\max\Biggl\{\biggl|k_{\rm tt}w_{\rm t}-F_{\rm t}-\frac{k_{\rm nt}}{k_{\rm nn}}\bigl\langle k_{\rm nt}w_{\rm t}-F_{\rm n}\bigr\rangle^+\biggr|\;,\;f\bigl\langle k_{\rm nt}w_{\rm t}-F_{\rm n}\bigr\rangle^-\Biggr\}.
\]
The corresponding unique minimizer \(\mathbf{u}\) of the minimum problem~\eqref{eq:discreteMinProb} is:
\begin{align*}
	u_{\rm n} &= -\frac{1}{k_{\rm nn}}\bigl\langle k_{\rm nt}w_{\rm t}-F_{\rm n}\bigr\rangle^+, \qquad & u_{\rm t} &= w_{\rm t},\\
	t_{\rm n} &= -\bigl\langle k_{\rm nt}w_{\rm t}-F_{\rm n}\bigr\rangle^-,\qquad & t_{\rm t} &= k_{\rm tt}w_{\rm t}-F_{\rm t}-\frac{k_{\rm nt}}{k_{\rm nn}}\bigl\langle k_{\rm nt}w_{\rm t}-F_{\rm n}\bigr\rangle^+,
\end{align*}
because it satisfies \(|t_{\rm t}|\leq \Sigma_1 \leq \sigma\). We have \(P(\sigma)\leq \Sigma_1\).  Setting:
\[
	\Sigma_2:=\max_{\sigma\in[0,\Sigma_1]}P(\sigma),\qquad \Sigma:=\max\{\Sigma_1,\Sigma_2\},
\]
we see that \(P:[0,\Sigma]\to [0,\Sigma]\) is a continuous function, which has at least one fixed point by Brouwer's theorem. \qed

\bigskip
In the discrete setting, the solution of the incremental problem is unique, for small friction coefficient. This has been known for a long time, although it is difficult to find a reference for this fact.   Examples of multiple solutions for the discrete incremental problem with large friction coefficient seem to have been mentioned for the first time in~\cite{Janovsky}.

\begin{prop}
	\label{prop:discreteIncrementalProblem}
	If \(fk_{\rm nt}/k_{\rm tt}<1\), then the mapping \(P\) defined by~\eqref{eq:defP} is a strict contraction (has Lipschitz-modulus strictly smaller than \(1\)) and the solution of the discrete incremental problem is therefore unique. In that case, the solution mapping \(\mathbf{F}\mapsto\mathbf{u}\) is Lipschitz-continuous.  If \(f=k_{\rm tt}/k_{\rm nt}\), then one can find an external force \(\mathbf{F}\) such that the discrete incremental problem has a continuum of solutions.
\end{prop}

\noindent\textbf{Proof.} 

\noindent\textbf{Step 1.} \textit{Uniqueness.}

Assume the condition \(fk_{\rm nt}/k_{\rm tt}<1\).  Consider \(\sigma_1,\sigma_2\in \left[0,+\infty\right[\) and denote by \(\mathbf{u}_1,\mathbf{u}_2\in K\) the corresponding minimizers~\eqref{eq:discreteMinProb}.  Setting \(\mathbf{t}_i:=\mathbf{K}\mathbf{u}_i-\mathbf{F}\), the contraction property to prove reads as:
\begin{equation}
	\label{eq:contraction}
		\bigl|t_{1{\rm n}}-t_{2{\rm n}}\bigr| \leq \frac{k_{\rm nt}}{k_{\rm tt}}\bigl|\sigma_1-\sigma_2\bigr|.
\end{equation}
Since \(t_{i{\rm n}}u_{i{\rm n}}=0\), it is sufficient to prove it for the following three cases.

\noindent\textbf{Case 1.} \textit{\(t_{1{\rm n}}=t_{2{\rm n}}=0\).} The expected inequality~\eqref{eq:contraction} is trivially satisfied.

\noindent\textbf{Case 2.} \textit{\(u_{1{\rm n}}=u_{2{\rm n}}=0\).} In that case, inequality~\eqref{eq:discreteOptIneqSum} reduces to:
\[
	k_{\rm tt}\bigl(u_{1{\rm t}}-u_{2{\rm t}}\bigr)^2 \leq |\sigma_1-\sigma_2|\bigl|u_{1{\rm t}}-u_{2{\rm t}}\bigr|,
\]
which, by \(t_{1{\rm n}}-t_{2{\rm n}}=k_{\rm nt}(u_{1{\rm t}}-u_{2{\rm t}})\), readily yields the expected inequality~\eqref{eq:contraction}.

\noindent\textbf{Case 3.} \textit{\(u_{1{\rm n}}=0\), \(t_{2{\rm n}}=0\).} We recall that \(t_{1{\rm n}}\leq 0\) and \(u_{2{\rm n}}\leq 0\), so that:
\[
	\underbrace{t_{1{\rm n}}-t_{2{\rm n}}}_{\mbox{}\leq 0} = \underbrace{k_{\rm nn}\bigl(u_{1{\rm n}}-u_{2{\rm n}}\bigr)}_{\mbox{}\geq 0} + \underbrace{k_{\rm nt}\bigl(u_{1{\rm t}}-u_{2{\rm t}}\bigr)}_{\mbox{}\leq 0},
\]
and:
\begin{equation}
	\label{eq:eqOne}
		|t_{1{\rm n}}-t_{2{\rm n}}| =  k_{\rm nt}\bigl|u_{1{\rm t}}-u_{2{\rm t}}\bigr| - k_{\rm nn}|u_{2{\rm n}}|.
\end{equation}
As the friction law can be equivalently rewritten as::
\[
	\forall v\in\mathbb{R},\qquad t_{i,{\rm t}}\bigl(v-u_{i{\rm t}}\bigr) + \sigma_i\bigr[|v-w_{\rm t}|-|u_{i{\rm t}}-w_{\rm t}|\bigr]\geq 0,
\]
we have:
\[
\bigl(t_{1{\rm t}}-t_{2{\rm t}}\bigr)\bigl(u_{1{\rm t}}-u_{2{\rm t}}\bigr) \leq |\sigma_1-\sigma_2|\bigl|u_{1{\rm t}}-u_{2{\rm t}}\bigr|,
\]
so that:
\[
	k_{\rm nt}\bigl(u_{1{\rm n}}-u_{2{\rm n}}\bigr)\bigl(u_{1{\rm t}}-u_{2{\rm t}}\bigr)+k_{\rm tt}\bigl(u_{1{\rm t}}-u_{2{\rm t}}\bigr)^2 \leq |\sigma_1-\sigma_2|\bigl|u_{1{\rm t}}-u_{2{\rm t}}\bigr|,
\]
which yields:
\begin{equation}
	\label{eq:eqTwo}
		-k_{\rm nt}|u_{2{\rm n}}|+k_{\rm tt}\bigl|u_{1{\rm t}}-u_{2{\rm t}}\bigr| =: S \leq |\sigma_1-\sigma_2|.
\end{equation}
Inverting the linear system formed by equations~\eqref{eq:eqOne} and~\eqref{eq:eqTwo}, we obtain:
\[
	|u_{2{\rm n}}| = \frac{-k_{\rm tt}|t_{1{\rm n}}-t_{2{\rm n}}|+k_{\rm nt}S}{k_{\rm nn}k_{\rm tt}-k_{\rm nt}^2}\geq 0,
\]
which yields~\eqref{eq:contraction} also in this case, as the determinant \(k_{\rm nn}k_{\rm tt}-k_{\rm nt}^2\) is positive.

\medskip
\noindent\textbf{Step 2.} \textit{Lipschitz-continuity of the solution mapping.}

Taking \(\sigma_1,\sigma_2\in \left[0,+\infty\right[\), \(\mathbf{F}_1,\mathbf{F}_2\in \mathbb{R}^2\) and denoting by \(\mathbf{u}_i:=\mathbf{u}(\mathbf{F}_i,\sigma_i)\in K\) the corresponding minimizers~\eqref{eq:discreteMinProb}, the sum of the corresponding optimality inequalities~\eqref{eq:discreteOptIneq} yields:
\begin{equation}
	\label{eq:optimalityLipschitz}
		(\mathbf{u}_1-\mathbf{u}_2)\cdot\mathbf{K}(\mathbf{u}_1-\mathbf{u}_2) \leq \bigl(\mathbf{F}_1-\mathbf{F}_2\bigr)\cdot\bigl(\mathbf{u}_1-\mathbf{u}_2\bigr) + |\sigma_1-\sigma_2|\bigl|u_{1{\rm t}}-u_{2{\rm t}}\bigr|,
\end{equation}
which shows that the mapping \((\mathbf{F},\sigma)\mapsto\mathbf{u}\) is Lipschitz-continuous.  Assume \(fk_{\rm nt}/k_{\rm tt}<1\).  By Step~1, we know that there exists a unique \(\sigma(\mathbf{F})\in\left[0,+\infty\right[\) such that:
\begin{equation}
	\label{eq:fixedPoint}
		P(\sigma) = -f\Bigl[k_{\rm nn}u_{\rm n}(\mathbf{F},\sigma)+k_{\rm nt}u_{\rm t}(\mathbf{F},\sigma)-F_{\rm n}\Bigr] = \sigma.
\end{equation}
Equation~\eqref{eq:contraction} yields:
\[
	\forall\mathbf{F}\in\mathbb{R}^2,\qquad	\biggl|\Bigl[k_{\rm nn}u_{\rm n}(\mathbf{F},\sigma_1)+k_{\rm nt}u_{\rm t}(\mathbf{F},\sigma_1)\Bigr] - \Bigl[k_{\rm nn}u_{\rm n}(\mathbf{F},\sigma_2)+k_{\rm nt}u_{\rm t}(\mathbf{F},\sigma_2)\Bigr]\biggr| \leq \frac{k_{\rm nt}}{k_{\rm tt}}\bigl|\sigma_1-\sigma_2\bigr|,
\]
and equation~\eqref{eq:optimalityLipschitz} yields:
\[
	\forall\sigma\in\left[0,+\infty\right[,\qquad\biggl|\Bigl[k_{\rm nn}u_{\rm n}(\mathbf{F}_1,\sigma)+k_{\rm nt}u_{\rm t}(\mathbf{F}_1,\sigma)\Bigr] - \Bigl[k_{\rm nn}u_{\rm n}(\mathbf{F}_2,\sigma)+k_{\rm nt}u_{\rm t}(\mathbf{F}_2,\sigma)\Bigr]\biggr| \leq C\bigl|\mathbf{F}_1-\mathbf{F}_2\bigr|,
\]
for some constant \(C>0\) depending only on \(\mathbf{K}\).  Bringing these last two estimates in formula~\eqref{eq:fixedPoint} shows the Lipschitz-continuity of the function \(\sigma(\mathbf{F})\):
\[
	\bigl|\sigma(\mathbf{F}_1)-\sigma(\mathbf{F}_2)\bigr|\leq  \frac{f(C+1)}{1-fk_{\rm nt}/k_{\rm tt}}\bigl|\mathbf{F}_1-\mathbf{F}_2\bigr|.
\]
The Lipschitz-continuity of the solution mapping \(\mathbf{F}\mapsto\mathbf{u}\bigl(\mathbf{F},\sigma(\mathbf{F})\bigr)\) now follows from~\eqref{eq:optimalityLipschitz}.

\medskip
\noindent\textbf{Step 3.} \textit{Example of multiple solutions.}

We now assume that \(f=k_{\rm tt}/k_{\rm nt}\), and we are going to construct a continuum of solutions for the discrete incremental problem with \(w_{\rm t}=0\).  We will look for solutions which all satisfy \(t_{\rm t}=ft_{\rm n}\). As:
\[
	u_{{\rm n}}= \overbrace{\frac{\bigl(k_{\rm tt}-fk_{\rm nt}\bigr)t_{{\rm n}}}{k_{\rm nn}k_{\rm tt}-k_{\rm nt}^2}}^{\mbox{}=0} + \frac{k_{\rm tt}F_{\rm n}-k_{\rm nt}F_{\rm t}}{k_{\rm nn}k_{\rm tt}-k_{\rm nt}^2}
\]
we will consider only external forces such that \(F_{\rm n}=k_{\rm nt}F_{\rm t}/k_{\rm tt}\) to ensure that \(u_{\rm n}=0\).  We now need to ensure that \(t_{\rm n}\leq 0\) and \(u_{\rm t}\geq 0\) in:
\[
u_{{\rm t}} = \frac{\bigl(fk_{\rm nn}-k_{\rm nt}\bigr)t_{{\rm n}}}{k_{\rm nn}k_{\rm tt}-k_{\rm nt}^2} + \frac{k_{\rm nn}F_{\rm t}-k_{\rm nt}F_{\rm n}}{k_{\rm nn}k_{\rm tt}-k_{\rm nt}^2},
\]
which is always possible, provided that \(F_{\rm t}>0\).

Finally, take \(F_{\rm t}>0\) arbitrarily, set \(F_{\rm n}=k_{\rm nt}F_{\rm t}/k_{\rm tt}\) with \(f=k_{\rm tt}/k_{\rm nt}\).  Take \(t_{\rm n}\) arbitrarily in the interval \(\bigl[-\frac{k_{\rm nt}F_{\rm t}}{k_{\rm tt}},0\bigr]\) and set \(t_{\rm t}=ft_{\rm n}\), \(\mathbf{u}:=\mathbf{K}^{-1}(\mathbf{F}+\mathbf{t})\).  Then, one has:
\begin{align*}
	t_{{\rm n}}&\leq 0, & u_{{\rm n}}&=0,\\
	t_{{\rm t}}&=ft_{{\rm n}}, & u_{{\rm t}}&\geq 0, \qquad\implies\qquad t_{{\rm t}}u_{{\rm t}}=ft_{{\rm n}}|u_{{\rm t}}|,
\end{align*}
hence a continuum of solutions of the discrete incremental problem.  \qed

\bigskip
We now turn to the analysis of the discrete quasi-static problem.  We first exhibit an example of a \emph{jumping} solution.

\begin{prop}
	\label{prop:jumping-solution}
	Let \(f\geq k_{\rm tt}/k_{\rm nt}\) and \(\mathbf{u}_{0}=\mathbf{0}\).  Pick an arbitrary \(R>0\) and define the piecewise affine Lipschitz-continuous external force \(\mathbf{F}:[0,2]\to\mathbb{R}^2\) by:
\[
	\mathbf{F}(s):=\left|\begin{array}{ll}
		\Bigl(sR/f,sR\Bigr) & \text{ if } s\in [0,1],\\[2.0ex]
		\Bigl(R/f+(s-1),R+(f+1)(s-1)\Bigr) & \text{ if } s\in [1,2].
	\end{array}\right.
\]
Then, the functions \(\mathbf{u},\mathbf{t}\in \textsl{BV}^+([0,2];\mathbb{R}^2)\) defined by:
\begin{align*}
	\mathbf{u}(s) & :=\left|\begin{array}{ll}\mathbf{0}, & \text{ if } s\in [0,1[,\\[2.0ex]\displaystyle
		\biggl(\frac{\overbrace{R(k_{\rm tt}/f-k_{\rm nt})+\bigl[k_{\rm tt}-k_{\rm nt}(f+1)\bigr](s-1)}^{\mbox{}\leq 0}}{k_{\rm nn}k_{\rm tt}-k_{\rm nt}^2}, & \\[2.0ex]\displaystyle
		\hspace*{2cm}\frac{R(k_{\rm nn}-k_{\rm nt}/f)+\bigl[k_{\rm nn}(f+1)-k_{\rm nt}\bigr](s-1)}{k_{\rm nn}k_{\rm tt}-k_{\rm nt}^2}\biggr) & \text{ if } s\in [1,2]
	\end{array}\right.\\
	\mathbf{t}(s) & :=\left|\begin{array}{ll}-\mathbf{F}(s)=-\Bigl(sR/f,sR\Bigr) & \text{ if } s\in [0,1[,\\[2.0ex]
		\mathbf{0}, & \text{ if } s\in [1,2]
	\end{array}\right.
\end{align*}
is a solution of the discrete quasi-static problem~\eqref{def:quasi-static-discrete} with a jump at \(s=1\).
\end{prop}

The proof of Proposition~\ref{prop:jumping-solution} is left to the reader. Interestingly, the jump in the above solution is unavoidable in the following sense.

\begin{prop} 
	\label{prop:jumping-discrete}
	Let \(f\geq k_{\rm tt}/k_{\rm nt}\), \(\mathbf{u}_{0}=\mathbf{0}\) and \(\mathbf{F}\) be the Lipschitz-continuous function defined in Proposition~\ref{prop:jumping-solution}.  Then, the solution \(\mathbf{u}\) defined in Proposition~\ref{prop:jumping-solution} is the only solution of the discrete quasi-static problem~\eqref{def:quasi-static-discrete} that equals \(\mathbf{0}\) for \(s\in [0,1[\).
\end{prop}

\noindent\textbf{Proof.}  At \(s=1\), \(\mathbf{t}(1-)=(-R/f,-R)\) with \(R>0\).  Assume that \(\mathbf{u}(s)\) is continuous at \(s=1\), then \(\mathbf{t}(s)\) is continuous at \(s=1\) as well.  As a result, \(t_{\rm n}(s)\) is negative on a right neighborhood of \(s=1\), and \(u_{\rm n}(s)=0\) on that right neighborhood of \(s=1\).  The continuous function \(t_{\rm t}(s)\) must also be negative on a right neighborhood of \(s=1\) and the Coulomb law implies that the restriction of the measure \(\dot{u}_{\rm t}\) to that right neighborhood must be nonnegative, which entails that the function \(u_{\rm t}(s)\) is nondecreasing and therefore nonnegative on that right neighborhood, as it vanishes at \(s=1\). But:
\begin{align*}
	t_{\rm n}(s) &= -F_{\rm n}(s) + k_{\rm nt}u_{\rm t}(s),\\
	t_{\rm t}(s) &= -F_{\rm t}(s) + k_{\rm tt}u_{\rm t}(s),\\
	t_{\rm t}(s) &\geq ft_{\rm n}(s),
\end{align*}
which yields:
\[
	\underbrace{(fk_{\rm nt}-k_{\rm tt})}_{\mbox{}\geq 0}u_{\rm t}(s) \leq fF_{\rm n}(s) - F_{\rm t}(s)=-(s-1),
\]
and therefore a contradiction. \qed

\bigskip
\noindent\textbf{Remark.} We could also take the left limit of the solution \(\mathbf{u}\) at \(s=1\) as initial condition for the discrete quasi-static problem for \(s\in [1,2]\) associated with external force \(\mathbf{F}\) defined for \(s\geq 1\). That way, Proposition~\ref{prop:jumping-solution} would still hold, showing that the solution \emph{must jump} at initial time \(s=1\), and therefore that the discrete quasi-static problem has no continuous solution, while the external force is a \(C^\infty\) function of time.

\bigskip
Interestingly, this result was informally discovered in the engineering literature, without any proper definition of what a jumping solution to the quasi-static problem would mean.  In~\cite{Klarbring}, Klarbring studies the \emph{rate} problem associated with the discrete quasi-static problem: given a configuration of the system and a force rate \(\dot{\mathbf{F}}\), find the velocity \(\dot{\mathbf{u}}\) of the system.  He is able to find examples where the rate problem has possibly multiple solutions and also examples where the rate problem has no solution at all, which he calls `examples of non-existence of solution for the quasi-static problem'.  Our jumping solution corresponds to a case where the rate problem has no solution.

We now turn to the analysis of the counterpart of Condition~$\mathscr{C}$, in the case of the discrete quasi-static problem. We cannot transpose directly Definition~\ref{def:conditionC} to the discrete setting, due to a basic difference between the continuum and the discrete problem. In the case of the continuum problem, both the external force \((\mathbf{F},\mathbf{T})\) and the contact force \(\mathbf{t}\) are uniquely determined from the knowledge of the current solution \(\mathbf{u}\). This is no longer the case in the discrete setting, as \(\mathbf{t}+\mathbf{F}=\mathbf{K}\mathbf{u}\).  To find the counterpart of Condition~$\mathscr{C}$ in the discrete setting, we will therefore turn to what Condition~$\mathscr{C}$ has achieved in the proof of existence of a solution for the quasi-static problem: it was designed as the minimal condition to ensure that the increment of the solution \(\mathbf{u}\) on one time step is uniformly controlled by the increment of the external force \(\mathbf{F}\) (Proposition~\ref{prop:estimateConditionC}).  In the discrete setting, this means that we need to ensure that the ratio of the increment of the solution on one time step and the increment of the external force is uniformly bounded.  The reader will check that under this condition, the demonstration of Theorem~\ref{theo:mainExistenceTheorem} can readily be adapted to the discrete quasi-static problem, yielding existence of an absolutely continuous solution \(\mathbf{u}\) whenever the external force \(\mathbf{F}\) is absolutely continuous.

\begin{prop}
	For the discrete quasi-static problem~\eqref{def:quasi-static-discrete}, uniform condition~$\mathscr{C}$ reduces to \(fk_{\rm nt}/k_{\rm tt}<1\), in the sense that, under this condition, if \(\mathbf{u_{i-1}}\) is the solution associated with force \(\mathbf{F}_{i-1}\) at time \(s_{i-1}\), in a time stepping approach of the quasi-static problem, then solution \(\mathbf{u}_i\) associated with the external force \(\mathbf{F}_i\) at time \(s_{i}\), as provided by the incremental problem, satisfies the uniform estimate:
	\[
	\bigl|\mathbf{u}_i-\mathbf{u}_{i-1}\bigr| \leq C\bigl|\mathbf{F}_i-\mathbf{F}_{i-1}\bigr|,
	\]
	for some constant \(C>0\) depending only on \(f\) and \(\mathbf{K}\), and where \(|\cdot|\) stands for the Euclidean norm.  Under this condition, the discrete quasi-static problem has a solution which is absolutely continuous whenever the external force \(\mathbf{F}\) is absolutely continuous.
\end{prop}

\noindent\textbf{Proof.}  Set \(w_{\rm t}:=u_{i-1,{\rm t}}\). First, assume \(fk_{\rm nt}/k_{\rm tt}<1\).  Then, note that \(\mathbf{u}_{i-1}\) is the unique solution, provided by Proposition~\ref{prop:discreteIncrementalProblem}, of the incremental problem associated with \(w_{\rm t}\) and external force \(\mathbf{F}_{i-1}\).  As \(\mathbf{u_{i}}\) is the unique solution of the incremental problem associated with \(w_{\rm t}\) and external force \(\mathbf{F}_{i}\), the claim follows from the Lispchitz-continuity of the solution mapping \(\mathbf{F}\mapsto\mathbf{u}\) provided by Proposition~\ref{prop:discreteIncrementalProblem}.  If \(fk_{\rm nt}/k_{\rm tt}\geq 1\), Proposition~\ref{prop:discreteIncrementalProblem} shows that we may have \(\mathbf{u}_i\neq \mathbf{u}_{i-1}\), while \(\mathbf{F}_i= \mathbf{F}_{i-1}\).  \qed

\bigskip
\noindent\textbf{Synthesis and discussion.} We have proved the following results.
\begin{itemize}
	\item The discrete incremental problem is always solvable, thanks to Brouwer's theorem.
	\item If \(fk_{\rm nt}/k_{\rm tt}<1\), the discrete incremental problem has a unique solution.
	\item If \(f=k_{\rm tt}/k_{\rm nt}\), the discrete incremental problem may have a continuum of solutions.
	\item If \(fk_{\rm nt}/k_{\rm tt}<1\), the discrete quasi-static problem has a solution which is absolutely continuous whenever the external force \(\mathbf{F}\) is absolutely continuous.
	\item If \(fk_{\rm nt}/k_{\rm tt}\geq 1\), the discrete quasi-static problem may have a solution with jumps, while the external force is absolutely continuous.
	\item If \(fk_{\rm nt}/k_{\rm tt}\geq 1\), the discrete quasi-static problem with absolutely continuous external force may have no continuous solution.
\end{itemize}
Incidentally, some examples were provided in~\cite{bibiEngSci} showing that the discrete quasi-static problem may have multiple solutions for \(C^\infty\) external forces and arbitrarily small friction coefficient \(f\). The results about the discrete incremental problem were already known, those about the discrete quasi-static problem seem to be new, as no proper definition of a jumping solution seems to have ever appeared in the literature. Discontinuous solutions to the discrete quasi-static problem were mentioned in~\cite{Gastaldi}, but this reference misses a clear formulation of the evolution problem in the bounded variation setting.  It was also known~\cite{Klarbring} that the rate-problem (finding the velocity of the particle, given the current configuration and the rate of the external force) associated with the discrete quasi-static problem may have no solution, and the author informally relates this to possible discontinuities in the solution of the quasi-static problem.

This has lead the engineering community to consider the value \(f=k_{\rm tt}/k_{\rm nt}\) as a critical value, connected to some transition in the qualitative response of the system, and to vibrations induced by friction. The community has therefore sought to define such a critical value for the friction coefficient \(f\) in more general settings, and, in particular, for the continuum quasi-static problem.  As this critical value was associated to a loss of uniqueness in the discrete setting, it was believed that the continuum quasi-static problem would also lose uniqueness for some critical value of the friction coefficient \(f\).  But, repeated attempts to study uniqueness in the continuum incremental problem were never able to deliver either a proof of uniqueness for the incremental problem with small friction coefficient, or an example of multiple solutions for arbitrary small friction coefficient.

Also, in the case of the continuum incremental problem, the only known existence results~\cite{BookJiri} available required a small friction coefficient \(f\), suggesting that the transition in the response of the system for a continuum would be in terms of the existence of a solution, rather than in terms of uniqueness.  The situation has therefore remained unclear for a long time, until our recent article~\cite{BI-ma3as} appeared, which has shown that the continuum incremental problem is always solvable, even for large friction coefficient \(f\).  This has clarified the situation and has paved the way to the present article.  In this article, we argue that the transition in the qualitative response of the system in the case of a continuum should be analyzed in terms of the existence (or not) 	of an absolutely continuous solution to the quasi-static problem with absolutely continuous external forces.

\subsection{Examples for a continuum}

In this section, we briefly show how the examples of solutions for the discrete incremental and quasi-static problems displayed in Sections~\ref{sec:discrete} can be transposed to examples of solutions in the continuum setting.  This transposition was first performed by P. Hild~\cite{Hild} to exhibit an example of a continuum of solutions for the continuum incremental problem. We briefly review his idea, and use it to exhibit an example of a jumping solution for the continuum quasi-static problem, under absolutely continuous load.

We consider the triangle \((ABC)\) in \(\mathbb{R}^2\) with vertices \(A=(-1,0)\), \(B=(0,0)\) and \(C=(x_C,y_C<0)\), filled with an isotropic homogeneous linear elastic material with Young's modulus \(E\) and Poisson's ratio \(\nu=0\).  The triangle is clamped on the edge \(BC\) (\(\Gamma_U=BC\)) and is loaded by a vanishing external volume force \(\mathbf{F}=(0,0)\) and a constant external traction \(\mathbf{T}=(T_x,T_y)\in\mathbb{R}^2\) on the edge \(AC\) (\(\Gamma_T=AC\)).  Finally, it obeys frictional contact conditions on the edge \(AB\) (\(\Gamma_C=AB\)) with respect to the rigid obstacle defined by the half-plane \(y\geq 0\) (\(g=0\)). The friction coefficient \(f\) is supposed to be constant on \(\Gamma_C\).

We will consider only displacements \(\mathbf{u}=(u_x(x,y),u_y(x,y))\) that are linearly interpolated between the zero displacement \(\mathbf{u}^B=\mathbf{u}^C=\mathbf{0}\) of vertices \(B\) and \(C\) and the displacement \(\mathbf{u}^A=(u_x^A,u_y^A)\) of vertex \(A\):
\begin{equation}
	\label{eq:linearDisplacement}
	\mathbf{u}(x,y) = \Biggl(\frac{x-x_B}{x_A-x_B} - \biggl(\frac{x_C-x_B}{x_A-x_B}\biggr)\biggl(\frac{y-y_B}{y_C-y_B}\biggr)\Biggr)\Biggl(1 - \biggl(\frac{x_C-x_B}{x_A-x_B}\biggr)\biggl(\frac{y_A-y_B}{y_C-y_B}\biggr)\Biggr)^{-1}\mathbf{u}^A.
\end{equation}
The (linearized) strain \(\boldsymbol{\varepsilon}(\mathbf{u}):=(\nabla\mathbf{u}+{}^t\nabla\mathbf{u})/2\) is therefore a constant matrix whose entries are linear functions of \(u_x^A\) and \(u_y^A\). The elastic energy takes the form:
\[
	\frac{|ABC|}{2}\boldsymbol{\varepsilon}(\mathbf{u}):\boldsymbol{\Lambda}\,\boldsymbol{\varepsilon}(\mathbf{u}) = \frac{1}{2}\mathbf{u}^A\cdot\tilde{\mathbf{K}}\mathbf{u}^A,
\]
for some constant definite positive symmetric matrix \(\tilde{\mathbf{K}}\) depending only on the geometry of the triangle and the elastic properties \(E\), \(\nu\) of the material.  As the stress  \(\boldsymbol{\sigma}(\mathbf{u})=\boldsymbol{\Lambda}\,\boldsymbol{\varepsilon}(\mathbf{u})\), is constant on the triangle, it is automatically divergence-free, and the associated traction \(\mathbf{t}=\boldsymbol{\sigma}(\mathbf{u})\cdot\mathbf{n}\) is constant on each edge of the triangle.  Setting:
\[
	\tilde{\mathbf{F}} := \frac{|AC|^2}{2}\mathbf{T},
\]
it is readily checked that a displacement field \(\mathbf{u}(x,y)\) of the form~\eqref{eq:linearDisplacement} solves the incremental problem for the triangle \(ABC\) if and only if \(\mathbf{u}^A\) solves the discrete incremental problem~\eqref{def:incremental-discrete} with stiffness matrix \(\tilde{\mathbf{K}}\), external force \(\tilde{\mathbf{F}}\) and friction coefficient \(f\).  Hence, the continuum of solutions for the discrete incremental problem, displayed in Proposition~\ref{prop:discreteIncrementalProblem}, can be transposed to a continuum of solutions for the continuum incremental problem.  Likewise, considering a time-varying \(\mathbf{T}\in \textsl{BV}^+([0,S];\mathbb{R}^2)\), a displacement \(\mathbf{u}(s)\) of the form~\eqref{eq:linearDisplacement}, for each \(s\in [0,S]\), solves the quasi-static problem, as formulated in Definition~\ref{def:quasi-static}, if and only if \(\mathbf{u}^A(s)\) solves the discrete quasi-static problem~\eqref{def:quasi-static-discrete} with stiffness matrix \(\tilde{\mathbf{K}}\), external force \(\tilde{\mathbf{F}}(s)\) and friction coefficient \(f\).  Hence, the jumping solution for the discrete quasi-static problem, displayed in Proposition~\ref{prop:jumping-discrete}, can be transposed to a jumping solution for the continuum quasi-static problem, with Lipschitz-continuous (in time) load.  We have therefore proved:

\begin{prop}
	In the case of a continuum, the incremental problem may have a continuum of solutions, and the quasi-static problem may have a jumping solution, even in the case of an absolutely continuous (in time) load.
\end{prop}

\section{Proof of Theorem~\ref{theo:mainExistenceTheorem}}
\label{sec:ProofExistence}

We can now turn to the proof of Theorem~\ref{theo:mainExistenceTheorem}, our main result.  The proof relies on the solvability of the incremental problem as proved in Theorem~\ref{theo:BIincremental} and on a mix of ideas from Moreau's sweeping process~\cite{Moreau} and from Brezis' pseudomonotonicity~\cite{BrezisIeq}.  The application of pseudomonotonicity to friction problems is very far from being straightforward and requires to prove a compensated-compactness-like result.  In the case of the two-dimensional problem, we were able to obtain the needed pseudomonotonicity to solve the incremental problem, from a new fine property of the elastic Neumann-to-Dirichlet operator, which we proved in a previous work~\cite{BI-ma3as}.  It will also play a crucial role in our existence proof for the quasistatic problem and will therefore currently restrict it to the two-dimensional case.  The extension to the three-dimensional case is still a work in progress, and raises deep issues in concentration-compactness.

In all the sequel, the analysis will be restricted to the two-dimensional case \(N=2\).  We first recall the statement of our theorem from~\cite{BI-ma3as}, for easy further reference.  We consider a direct orthonormal basis \((\mathbf{n},\boldsymbol{\tau})\) on \(\Gamma_C\), where \(\mathbf{n}\) still stands for the outward unit normal.  Given arbitrary \(t\in H^{-1/2}(\Gamma_C)\), we consider the two well-posed linear elastic problems \(\mathscr{P}_{\rm n}(t)\), \(\mathscr{P}_{\rm t}(t)\) of finding the displacement \(\mathbf{u}\in H^1(\Omega;\mathbb{R}^2)\) such that:
\[
	\mathscr{P}_{\rm n}(t)\left\{\quad
\begin{array}{ll}
	\text{div}\,\boldsymbol{\sigma}(\mathbf{u}) = \mathbf{0}, & \text{ in }\Omega,\\[0.4ex]
		\mathbf{u} = \mathbf{0}, & \text{ on }\Gamma_U,\\[0.4ex]
		\boldsymbol\sigma(\mathbf{u})\,\mathbf{n} = \mathbf{0}, & \text{ on }\Gamma_T,\\[0.4ex]
		\boldsymbol\sigma(\mathbf{u})\,\mathbf{n}=t\mathbf{n}, & \text{ on }\Gamma_C,
\end{array}
\right.\qquad
\mathscr{P}_{\rm t}(t)\left\{\quad
\begin{array}{ll}
		\text{div}\,\boldsymbol{\sigma}(\mathbf{u}) = \mathbf{0}, & \text{ in }\Omega,\\[0.4ex]
		\mathbf{u} = \mathbf{0}, & \text{ on }\Gamma_U,\\[0.4ex]
		\boldsymbol\sigma(\mathbf{u})\,\mathbf{n} = \mathbf{0}, & \text{ on }\Gamma_T,\\[0.4ex]
		\boldsymbol\sigma(\mathbf{u})\,\mathbf{n}=t\boldsymbol{\tau}, & \text{ on }\Gamma_C,
\end{array}
\right.
\]

\begin{defi}
	\label{defi:linearOperators}
	Let \(L_{\rm n},L_{\rm t}:H^{-1/2}(\Gamma_C)\to H^{1/2}(\Gamma_C)\) be the two linear operators that map \(t\in H^{-1/2}(\Gamma_C)\) to the normal component \(u_{\rm n}\in H^{1/2}(\Gamma_C)\) of the solution \(\mathbf{u}\in H^1(\Omega;\mathbb{R}^2)\) of the problem \(\mathscr{P}_{\rm n}(t)\) (respectively, \(\mathscr{P}_{\rm t}(t)\)).
\end{defi}

\begin{prop}
	\label{prop:bilinearMapping}
	The bilinear mapping \(t_1,t_2\mapsto \langle t_1,L_{\rm n}(t_2)\rangle\) is a scalar product on \(H^{-1/2}(\Gamma_C)\), that induces a norm on \(H^{-1/2}(\Gamma_C)\) that is equivalent to the norm of \(H^{-1/2}(\Gamma_C)\).
\end{prop}

\noindent\textbf{Proof.}  Denote by \(\mathbf{u}(t)\) the solution of the problem \(\mathscr{P}_{\rm n}(t)\). Since:
\[
	\bigl\langle t_1,L_{\rm n}(t_2)\bigr\rangle = \int_{\Omega} \boldsymbol{\varepsilon}\bigl(\mathbf{u}(t_1)\bigr):\boldsymbol{\Lambda}\boldsymbol{\varepsilon}\bigl(\mathbf{u}(t_2)\bigr),
\]
it is obvious that the bilinear mapping \(t_1,t_2\mapsto \langle t_1,L_{\rm n}(t_2)\rangle\) is a scalar product on \(H^{-1/2}(\Gamma_C)\). Combining inequality~\eqref{eq:tBound} with the open mapping theorem now yields the claim.\qed

\bigskip

We focus first on the isotropic case. The following compensated-compactness-like result was first proved in~\cite{BI-ma3as}.

\begin{theo}[{\cite[Theorem 2.9]{BI-ma3as}}]
	\label{theo:BIcompensatedCompactness}
	We suppose that \(N=2\) and \(\boldsymbol{\Lambda}\in W^{1,\infty}(\Omega)\) is \emph{isotropic} at each \(x\in \Omega\).  Let \(\tau^k\), \(t_{\rm t}^k\) be two weakly converging sequences in \(H^{-1/2}(\Gamma_C)\). We denote their weak limits by \(\tau\), \(t_{\rm t}\), respectively.  We assume that:
	\[
		\forall k\in\mathbb{N},\quad\forall v\in H^{1/2}(\Gamma_C),\qquad \bigl\langle t_{\rm t}^k,v\bigr\rangle \leq \bigl\langle \tau^k,|v|\bigr\rangle,
	\]	
	which means that \(\tau^k\), \(t_{\rm t}^k\) are measures such that \(|t_{\rm t}^k|\leq \tau^k\).  Then:
	\[
		\lim_{k\to+\infty}\Bigl\langle \tau^k,L_{\rm t}\bigl(t_{\rm t}^k\bigr)\Bigr\rangle = \Bigl\langle \tau,L_{\rm t}\bigl(t_{\rm t}\bigr)\Bigr\rangle.
	\]
\end{theo}	

The above theorem was crucial in the proof of the existence of a solution to the incremental problem in dimension \(N=2\) (Theorem~\ref{theo:BIincremental}).  We recall (formula~\eqref{eq:defV}) that \(V\) denotes the closed subspace of \(H^1(\Omega;\mathbb{R}^2)\) defined by the Dirichlet boundary condition on \(\Gamma_U\).

We are given arbitrary \(\mathbf{F}\in \textsl{BV}^+([0,S];L^2(\Omega;\mathbb{R}^2))\), \(\mathbf{T}\in \textsl{BV}^+([0,S];L^2(\Gamma_T;\mathbb{R}^2))\) and an arbitrary initial condition \(\mathbf{u}_0\in V\), supposed to be compatible with \(\mathbf{F}(0)\), \(\mathbf{T}(0)\) in the sense that:
\begin{align}
		& \bullet \quad\forall \mathbf{v}\in V,\qquad \int_{\Omega} \boldsymbol{\varepsilon}(\mathbf{u}_0):\boldsymbol{\Lambda}\boldsymbol{\varepsilon}(\mathbf{v}) \geq \int_{\Omega} \mathbf{F}(0)\cdot\mathbf{v} + \int_{\Gamma_T}\mathbf{T}(0)\cdot\mathbf{v} + \int_{\Gamma_C}t_{0,{\rm n}}v_{\rm n}+\int_{\Gamma_C}ft_{0,{\rm n}}|v_{\rm t}|,\label{eq:compatibilityCondition}
		\\
		& \bullet \quad t_{0,{\rm n}}\leq 0,\;\text{ on }\Gamma_C,\quad \text{ and }\quad	\forall t\in H^{-1/2}(\Gamma_C),\text{ such that }\,t\leq 0,\quad \int_{\Gamma_C}\bigl(u_{0,{\rm n}}-g\bigr)\bigl(t-t_{0,{\rm n}}\bigr)\geq 0,\nonumber
\end{align}
where \(t_{0,{\rm n}}\) stands for the normal component of \(\boldsymbol{\sigma}(\mathbf{u}_0)\mathbf{n}\).

To prove the existence of a solution, we are going to construct a sequence of piecewise constant functions \(\mathbf{u}_+^m\in \textsl{BV}^+([0,S];V)\), by solving inductively finitely many incremental problems. Then, we will prove that the sequence \(\mathbf{u}_+^m\) converges uniformly on \([0,S]\) towards a solution \(\mathbf{u}\in \textsl{BV}^+([0,S];V)\) of the evolution inequality in Theorem~\ref{theo:mainExistenceTheorem}.

We first construct an appropriate subdivision \(0=s_0^m<s_1^m<\cdots<s_{n_m}^m=S\) of the time interval \([0,S]\).  It is designed in such a way that the corresponding piecewise constant right continuous approximations of the loads \(\mathbf{F}\), \(\mathbf{T}\) converge to \(\mathbf{F}\), \(\mathbf{T}\) uniformly (see Proposition~\ref{prop:convData} below).  Recalling definition~\eqref{eq:defVar} of the variation of a function, we define the variation \(v(s)\) of the load over the time interval \([0,s]\) as:
\begin{equation}
	\label{eq:defVarLoad}
		v(s):=\text{var}\bigl(\mathbf{F};[0,s]\bigr) + \text{var}\bigl(\mathbf{T};[0,s]\bigr),
\end{equation}
so that \(v(s)\) is a right-continuous nondecreasing bounded function of \(s\in [0,S]\). Picking \(m\in \mathbb{N}\), the \(s_i^m\) are defined inductively, starting from \(s_0^m=0\), by:
\[
	s_i^m := \sup\biggl\{s\in \bigl]s_{i-1}^m,S\bigr]\bigm| v(s)\leq v(s_{i-1}^m) + \frac{v(S)}{m+1}\biggr\},
\]
which is well-defined (the supremum is taken over a nonempty set), thanks to the right-continuity of \(v\). The subdivision is finite, with \(n_m\leq m+1\), since \(v(s_{i}^m)-v(s_{i-1}^m)\geq v(S)/(m+1)\).  We set:
\begin{align*}
	\mathbf{F}_+^m(s) & 	:= \left|\begin{array}{ll}\mathbf{F}_{i-1}^m:=\mathbf{F}(s_{i-1}^m),\quad & \text{if }s\in \left[s_{i-1}^m,s_i^m\right[,\\[1.0ex]
	\mathbf{F}_{n_m}^m:=\mathbf{F}(S),\quad &\text{if }s=S,\end{array}\right.\\[1.5ex]
	\mathbf{F}_-^m(s) & 	:= \left|\begin{array}{ll}
		\mathbf{F}_{0}^m:=\mathbf{F}(0),\quad & \text{if }s=0,\\[1.0ex]
	\mathbf{F}_{i}^m:=\mathbf{F}(s_{i}^m),\quad & \text{if }s\in \left]s_{i-1}^m,s_i^m\right],\end{array}\right.
\end{align*}
and \(\mathbf{T}_+^m\), \(\mathbf{T}_-^m\) are defined similarly.  We immediately obtain the following proposition.

\begin{prop}
	\label{prop:convData}
	The piecewise constant functions \(\mathbf{F}_+^m\), \(\mathbf{F}_-^m\), \(\mathbf{T}_+^m\), \(\mathbf{T}_-^m\) satisfy the following estimates, for all \(m\in \mathbb{N}\) and all \(s\in [0,S]\):
	\begin{align*}
		\bigl\|\mathbf{F}_+^m(s)-\mathbf{F}(s)\bigr\|_{L^2} + \bigl\|\mathbf{T}_+^m(s)-\mathbf{T}(s)\bigr\|_{L^2} & \leq \frac{v(S)}{m+1},\\
		\bigl\|\mathbf{F}_-^m(s)-\mathbf{F}(s)\bigr\|_{L^2} + \bigl\|\mathbf{T}_-^m(s)-\mathbf{T}(s)\bigr\|_{L^2} & 	\leq v\bigl(\sigma^m(s)\bigr)-v(s),\\[1.0ex]
		\text{\rm var}\bigl(\mathbf{F}_+^m;[0,S]\bigr) + \text{\rm var}\bigl(\mathbf{T}_+^m;[0,S]\bigr) & \leq v(S),
	\end{align*}
	where \(\sigma^m(s):=s_i^m\), for \(s\in \left]s_{i-1}^m,s_i^m\right]\), and \(\sigma^m(0):=0\).
	Therefore, the functions \(\mathbf{F}_+^m\) and \(\mathbf{T}_+^m\) converge uniformly on \([0,S]\) towards \(\mathbf{F}\) and \(\mathbf{T}\), respectively, as \(m\to\infty\).  Recalling that \(v\) is right-continuous, the functions \(\mathbf{F}_-^m\) and \(\mathbf{T}_-^m\) converge pointwisely in \([0,S]\) towards \(\mathbf{F}\) and \(\mathbf{T}\), respectively, as \(m\to\infty\), and are uniformly bounded on \([0,S]\). Therefore, for any finite positive measure \(\mu\) on \([0,S]\), the function \((\mathbf{F}_-^m,\mathbf{T}_-^m)\) converges strongly in \(L^2([0,S],\mu;L^2(\Omega;\mathbb{R}^2)\times L^2(\Gamma_T;\mathbb{R}^2))\) towards \((\mathbf{F},\mathbf{T})\), by the dominated convergence theorem.
\end{prop}

We now construct \((\mathbf{u}_i^m,t_{i,{\rm n}}^m)\in V\times H^{-1/2}(\Gamma_C)\) inductively (\(i\in\{0,\ldots,n_m\}\)), starting from \(\mathbf{u}_0^m:=\mathbf{u}_0\), \(t_{0,{\rm n}}^m:=t_{0,{\rm n}}\), by solving the following incremental problem on the basis of Theorem~\ref{theo:BIincremental}:
\begin{align}
	& \bullet \quad\forall \mathbf{v}\in V,\qquad
	\int_{\Omega} \boldsymbol{\varepsilon}(\mathbf{u}_i^m):\boldsymbol{\Lambda}\boldsymbol{\varepsilon}(\mathbf{v}-\mathbf{u}_i^m) \geq \int_{\Omega} \mathbf{F}_i^m\cdot(\mathbf{v}-\mathbf{u}_i^m) + \int_{\Gamma_T}\mathbf{T}_i^m\cdot(\mathbf{v}-\mathbf{u}_i^m) \nonumber\\ & \hspace*{5cm} \mbox{} + \int_{\Gamma_C}t_{i,{\rm n}}^m\,(v_{\rm n}-u_{i,{\rm n}}^m)+\int_{\Gamma_C}ft_{i,{\rm n}}^m\Bigl[\bigl|v_{\rm t}-u_{i-1,{\rm t}}^m\bigr|-\bigl|u_{i,{\rm t}}^m-u_{i-1,{\rm t}}^m\bigr|\Bigr].\label{eq:optimalityInequality}
	\\
	& \bullet \quad t_{i,{\rm n}}^m\leq 0,\;\text{ on }\Gamma_C,\quad \text{ and }\quad	\forall t\in H^{-1/2}(\Gamma_C),\text{ such that }\:t\leq 0,\quad \int_{\Gamma_C}\bigl(u_{i,{\rm n}}^m-g\bigr)\bigl(t-t_{i,{\rm n}}^m\bigr)\geq 0.\label{eq:varIneq}
\end{align}

We set:
\begin{align*}
	\mathbf{u}_+^m(s) & := \left|\begin{array}{ll}\mathbf{u}_{i-1}^m, & \text{if } s\in \left[s_{i-1}^m,s_i^m\right[,\\[1.0ex]
	\mathbf{u}_{n_m}^m, & \text{if } s=S,\end{array}\right.\qquad &
	t_{+,{\rm n}}^m(s) & := \left|\begin{array}{ll}t_{i-1}^m, & \text{if } s\in \left[s_{i-1,{\rm n}}^m,s_i^m\right[,\\[1.0ex]
	t_{n_m,{\rm n}}^m, & \text{if } s=S,\end{array}\right.\\[1.5ex]
	\mathbf{u}_-^m(s) & 	:= \left|\begin{array}{ll}
	\mathbf{u}_{0}^m=\mathbf{u}_{0},\quad & \text{if }s=0,\\[1.0ex]
\mathbf{u}_{i}^m,\quad & \text{if }s\in \left]s_{i-1}^m,s_i^m\right],\end{array}\right.\qquad &
	t_{-,{\rm n}}^m(s) & := \left|\begin{array}{ll}
	t_{0,{\rm n}}^m,\quad & \text{if }s=0,\\[1.0ex]
	t_{i,{\rm n}}^m,\quad & \text{if }s\in \left]s_{i-1}^m,s_i^m\right].\end{array}\right.
\end{align*}
These piecewise constant functions actually yield an exact solution for the quasi-static problem with piecewise constant load.

\bigskip
\begin{prop}
	\label{lem:approxSweeping}
	The functions \(\mathbf{u}_+^m\in \textsl{BV}^+([0,S];V)\), \(t_{+,{\rm n}}^m\in \textsl{BV}^+([0,S];H^{-1/2}(\Gamma_C))\) satisfy:
	\begin{multline*}
		\bullet\;\forall \mathbf{v}\in \mathscr{M}([0,S];V),\qquad \int_{[0,S]}\int_{\Omega} \boldsymbol{\varepsilon}(\mathbf{u}_+^m(s)):\boldsymbol{\Lambda}\boldsymbol{\varepsilon}\bigl(\mathbf{v}-\mathbf{\dot{u}}_+^m\bigr) \geq \int_{[0,S]}\int_{\Omega} \mathbf{F}_+^m(s)\cdot\bigl(\mathbf{v}-\mathbf{\dot{u}}_+^m\bigr) +\mbox{}\\\mbox{}+ \int_{[0,S]}\int_{\Gamma_T}\mathbf{T}_+^m(s)\cdot\bigl(\mathbf{v}-\mathbf{\dot{u}}_+^m\bigr) + \int_{[0,S]}\int_{\Gamma_C}t_{+,{\rm n}}^m(s)\,\Bigl[v_{\rm n}-\dot{u}_{+,{\rm n}}^m\Bigr]\\ \mbox{} + \int_{[0,S]}\int_{\Gamma_C}ft_{+,{\rm n}}^m(s)\,\Bigl[\bigl|v_{\rm t}\bigr|-|\dot{u}_{+,{\rm t}}^m\bigr|\Bigr],
	\end{multline*}
	\(\displaystyle\hspace*{0.45cm}\bullet\;\:\forall s\in [0,S],\;
	t_{+,{\rm n}}^m(s)\leq 0,\quad\forall t\in H^{-1/2}(\Gamma_C), \text{ with }\,t\leq 0,\quad \int_{\Gamma_C}\bigl(u_{+,{\rm n}}^m(s)-g\bigr)\bigl(t-t_{+,{\rm n}}^m(s)\bigr)\geq 0,\)	

	\bigskip
	\noindent{}where the precise meaning of the integrals in the first condition is given by Definition~\ref{defi:intMeasure}. Therefore, \(\mathbf{u}_+^m\) solves the quasi-static problem with data \(\mathbf{F}_+^m\), \(\mathbf{T}_+^m\) and initial condition \(\mathbf{u}_0\), in the sense of Definition~\ref{def:quasi-static}.
\end{prop}

\noindent\textbf{Proof.} Taking successively \(\mathbf{v}=\mathbf{u}_{i-1}^m\) and \(\mathbf{v}=2\mathbf{u}_{i}^m-\mathbf{u}_{i-1}^m\) in inequality~\eqref{eq:optimalityInequality} for \(\mathbf{u}_i^m\), we obtain:
\begin{multline*}
	\int_{\Omega} \boldsymbol{\varepsilon}(\mathbf{u}_i^m):\boldsymbol{\Lambda}\boldsymbol{\varepsilon}(\mathbf{u}_i^m-\mathbf{u}_{i-1}^m) = \int_{\Omega} \mathbf{F}_i^m\cdot(\mathbf{u}_i^m-\mathbf{u}_{i-1}^m) + \int_{\Gamma_T}\mathbf{T}_i^m\cdot(\mathbf{u}_i^m-\mathbf{u}_{i-1}^m) \\\mbox{} + \int_{\Gamma_C}t_{i,{\rm n}}^m\,(u_{i,{\rm n}}^m-u_{i-1,{\rm n}}^m)+\int_{\Gamma_C}ft_{i,{\rm n}}^m\bigl|u_{i,{\rm t}}^m-u_{i-1,{\rm t}}^m\bigr|,
\end{multline*}
which yields:
\begin{multline}
	\int_{[0,S]}\int_{\Omega} \boldsymbol{\varepsilon}(\mathbf{u}_+^m):\boldsymbol{\Lambda}\boldsymbol{\varepsilon}(\dot{\mathbf{u}}_+^m) = \int_{[0,S]}\int_{\Omega} \mathbf{F}_+^m\cdot\dot{\mathbf{u}}_+^m + \int_{[0,S]}\int_{\Gamma_T}\mathbf{T}_+^m\cdot\dot{\mathbf{u}}_+^m \\ \mbox{} + \int_{[0,S]}\int_{\Gamma_C}\Bigl(t_{+,{\rm n}}^m\,\dot{u}_{+,{\rm n}}^m+ft_{+,{\rm n}}^m\bigl|\dot{u}_{+,{\rm t}}^m\bigr|\Bigr).\label{eq:secondTermIneqEvol}
\end{multline}

Taking \(\lambda>0\), using \(\lambda \mathbf{v}\) as test function in the inequality~\eqref{eq:optimalityInequality} for \(\mathbf{u}_i^m\), dividing by \(\lambda\) and taking the limit \(\lambda\to+\infty\), we obtain:
\[	
		\forall \mathbf{v}\in V,\qquad
		\int_{\Omega} \boldsymbol{\varepsilon}(\mathbf{u}_i^m):\boldsymbol{\Lambda}\boldsymbol{\varepsilon}(\mathbf{v}) \geq \int_{\Omega} \mathbf{F}_i^m\cdot\mathbf{v} + \int_{\Gamma_T}\mathbf{T}_i^m\cdot\mathbf{v}  + \int_{\Gamma_C}t_{i,{\rm n}}^m\,v_{\rm n}+\int_{\Gamma_C}ft_{i,{\rm n}}^m\bigl|v_{\rm t}\bigr|,
\]
which yields:
\begin{equation}
	\label{eq:firstTermIneqEvol}
		\int_{[0,S]}\int_{\Omega} \boldsymbol{\varepsilon}(\mathbf{u}_+^m):\boldsymbol{\Lambda}\boldsymbol{\varepsilon}(\mathbf{v}) \geq \int_{[0,S]}\int_{\Omega} \mathbf{F}_+^m\cdot\mathbf{v} + \int_{[0,S]}\int_{\Gamma_T}\mathbf{T}_+^m\cdot\mathbf{v}  + \int_{[0,S]}\int_{\Gamma_C}t_{+,{\rm n}}^m\,v_{\rm n}+\int_{[0,S]}\int_{\Gamma_C}ft_{+,{\rm n}}^m\bigl|v_{\rm t}\bigr|,
\end{equation}
for all \(\mathbf{v}\in \mathscr{M}([0,S];V)\), whose Radon-Nikodym derivative \(\bar{\mathbf{v}}:={\rm d}\mathbf{v}/{\rm d}\|\mathbf{v}\|\) of \(\mathbf{v}\) with respect to its total variation \(\|\mathbf{v}\|\) is piecewise constant. As piecewise constant functions are dense in \(L^1([0,S],\mu;V)\) for any finite positive measure \(\mu\) on \([0,S]\), we can extend the above inequality to all \(\mathbf{v}\in \mathscr{M}([0,S];V))\) by density. \qed

\bigskip
We now need uniform estimates for the \((\mathbf{u}_i^m,t_{i,{\rm n}}^m)\).  In that aim, we need to be able to control the increment in the solution by the increment in the load.  This is not unconditionally possible, as proved by the discrete system studied in Section~\ref{sec:discrete}.   Being able to have such a control exactly amounts to discarding spontaneous jumps in the solution when the load does not jump.   Condition~$\mathscr{C}$ in Definition~\ref{def:conditionC} was precisely designed as the minimal condition to ensure such a control.  Combined with condition~\eqref{eq:optimalityInequality}, it yields the following uniform estimates.

\begin{prop}
	\label{prop:uniformEstimatesSweepingProcess}
	Under the hypotheses of Theorem~\ref{theo:mainExistenceTheorem}, for all \(m\in \mathbb{N}\) and \(i\in\{1,2,\ldots,n_m\}\), one has the uniform estimates:
	\begin{align*}
		\bigl\|\Delta\mathbf{u}_i^m\bigl\|_{H^{1}} +\bigl\|\Delta t_{i,{\rm n}}^m\bigr\|_{H^{-1/2}} & \leq C\Bigl\{\bigl\|\Delta\mathbf{F}_i^m\bigl\|_{L^{2}}+\bigl\|\Delta\mathbf{T}_i^m\bigl\|_{L^{2}}\Bigr\},\\
		\bigl\|\mathbf{u}_i^m\bigl\|_{H^1} +\bigl\|t_{i,{\rm n}}^m\bigr\|_{H^{-1/2}} & \leq C\Bigl\{\bigl\|\mathbf{u}_0\bigl\|_{H^1} + \bigl\|\mathbf{F}(0)\bigl\|_{L^2}+\bigl\|\mathbf{T}(0)\bigl\|_{L^2} + \text{\rm var}\bigl(\mathbf{F};[0,S]\bigl)+\text{\rm var}\bigl(\mathbf{T};[0,S]\bigl)\Bigr\},
	\end{align*}
	for some constant \(C>0\) depending only on \(f\), \(\boldsymbol{\Lambda}\), \(\Omega\), \(\Gamma_U\), \(\Gamma_T\), \(\Gamma_C\).  In these estimates, \(\Delta\mathbf{u}_i^m:=\mathbf{u}_i^m-\mathbf{u}_{i-1}^m\) denotes the increment.
\end{prop}

\noindent\textbf{Proof.}  

\noindent\textbf{First estimate.} We first take \(\mathbf{v}=\mathbf{u}_{i-1}^m\) in inequality~\eqref{eq:optimalityInequality} for \(\mathbf{u}_i^m\):
\[
	\int_{\Omega} \boldsymbol{\varepsilon}(\mathbf{u}_i^m):\boldsymbol{\Lambda}\boldsymbol{\varepsilon}(\Delta\mathbf{u}_i^m) \leq \int_{\Omega} \mathbf{F}_i^m\cdot\Delta\mathbf{u}_i^m + \int_{\Gamma_T}\mathbf{T}_i^m\cdot\Delta\mathbf{u}_i^m+ \int_{\Gamma_C}t_{i,{\rm n}}^m\Delta u_{i,{\rm n}}^m+\int_{\Gamma_C}ft_{i,{\rm n}}^m\bigl|\Delta u_{i,{\rm t}}^m\bigr|,
\]
and then take \(\mathbf{v}=\mathbf{u}_i^m\) in inequality~\eqref{eq:optimalityInequality} for \(\mathbf{u}_{i-1}^m\) (or in the compatibility condition~\eqref{eq:compatibilityCondition} in the case \(i=1\)):
\begin{multline*}
	-\int_{\Omega} \boldsymbol{\varepsilon}(\mathbf{u}_{i-1}^m):\boldsymbol{\Lambda}\boldsymbol{\varepsilon}(\Delta\mathbf{u}_i^m) \leq -\int_{\Omega}\mathbf{F}_{i-1}^m\cdot\Delta\mathbf{u}_i^m-\int_{\Gamma_T}\mathbf{T}_{i-1}^m\cdot\Delta\mathbf{u}_i^m \\ \mbox{} - \int_{\Gamma_C}t_{i-1,{\rm n}}^m\Delta u_{i,{\rm n}}^m+\int_{\Gamma_C}f\bigl(-t_{i-1,{\rm n}}^m\bigr)\underbrace{\Bigl[\bigl|u_{i,{\rm t}}^m-u_{i-2,{\rm t}}^m\bigr|-\bigl|u_{i-1,{\rm t}}^m-u_{i-2,{\rm t}}^m\bigr|\Bigr]}_{\displaystyle\mbox{}\leq \bigl|\Delta u_{i,{\rm t}}^m\bigr|}.
\end{multline*}
The sum of the two inequalities yields:
\[
	\int_{\Omega} \boldsymbol{\varepsilon}(\Delta\mathbf{u}_i^m):\boldsymbol{\Lambda}\boldsymbol{\varepsilon}(\Delta\mathbf{u}_i^m) \leq \int_{\Omega} \Delta\mathbf{F}_i^m\cdot\Delta\mathbf{u}_i^m + \int_{\Gamma_T}\Delta\mathbf{T}_i^m\cdot\Delta\mathbf{u}_i^m+ \int_{\Gamma_C}\Delta t_{i,{\rm n}}^m\Delta u_{i,{\rm n}}^m+f\Delta t_{i,{\rm n}}^m\bigl|\Delta u_{i,{\rm t}}^m\bigr|.
\]
As above,
\[
	\int_{\Gamma_C}\Delta t_{i,{\rm n}}^m\Delta u_{i,{\rm n}}^m = \int_{\Gamma_C}\bigl(t_{i,{\rm n}}^m-t_{i-1,{\rm n}}^m\bigr)\Bigl[\bigl(u_{i,{\rm n}}^m-g\bigr)-\bigl(u_{i-1,{\rm n}}^m-g\bigr)\Bigr]\leq 0,
\]
so that:
\begin{equation}
	\label{eq:firstIneq}
	\int_{\Omega} \boldsymbol{\varepsilon}(\Delta\mathbf{u}_i^m):\boldsymbol{\Lambda}\boldsymbol{\varepsilon}(\Delta\mathbf{u}_i^m) - \int_{\Gamma_C}f\Delta t_{i,{\rm n}}^m\bigl|\Delta \mathbf{u}_{i,{\rm t}}^m\bigr|\leq \int_{\Omega} \Delta\mathbf{F}_i^m\cdot\Delta\mathbf{u}_i^m + \int_{\Gamma_T}\Delta\mathbf{T}_i^m\cdot\Delta\mathbf{u}_i^m.
\end{equation}
Now, condition~$\mathscr{C}$ of Definition~\ref{def:conditionC} ensures that:
\begin{multline}
	\int_{\Omega} \boldsymbol{\varepsilon}(\Delta\mathbf{u}_i^m):\boldsymbol{\Lambda}\boldsymbol{\varepsilon}(\Delta\mathbf{u}_i^m) - \int_{\Gamma_C}f\Delta t_{i,{\rm n}}^m\bigl|\Delta \mathbf{u}_{i,{\rm t}}^m\bigr|\geq \alpha\int_{\Omega} \boldsymbol{\varepsilon}(\Delta\mathbf{u}_i^m):\boldsymbol{\Lambda}\boldsymbol{\varepsilon}(\Delta\mathbf{u}_i^m)\\\mbox{}-\beta\Biggl(\int_{\Omega} \boldsymbol{\varepsilon}(\Delta\mathbf{u}_i^m):\boldsymbol{\Lambda}\boldsymbol{\varepsilon}(\Delta\mathbf{u}_i^m)\Biggr)^{1/2}\Bigl(\bigl\|\Delta\mathbf{F}_i^m\bigr\|_{L^2(\Omega)} + \bigl\|\Delta\mathbf{T}_i^m\bigr\|_{L^2(\Gamma_T)}\Bigr),\label{eq:secondIneq}
\end{multline} 
where we denote the supremum in condition~$\mathscr{C}$ by \(1-\alpha\), so that \(\alpha>0\). Gathering inequalities~\eqref{eq:firstIneq},~\eqref{eq:secondIneq} and \eqref{eq:tBound} is now sufficient to yield the first estimate of the proposition. Actually, Proposition~\ref{prop:estimateConditionC} yields directly the first estimate of the proposition. The calculation in the proof of Proposition~\ref{prop:estimateConditionC} has just been repeated here for the sake of clarity.

\medskip
\noindent\textbf{Second estimate.} The  triangle inequality, the first estimate and the definition of the variation yield:
\begin{align*}
	\bigl\|\mathbf{u}_i^m - \mathbf{u}_0\bigl\|_{H^1} +\bigl\|t_{i,{\rm n}}^m - t_{0,{\rm n}}\bigr\|_{H^{-1/2}} & \leq C\sum_{j=1}^i\Bigl\{\bigl\|\Delta\mathbf{F}_j^m\bigl\|_{L^{2}}+\bigl\|\Delta\mathbf{T}_j^m\bigl\|_{L^{2}}\Bigr\} \\ \mbox{} & \leq C\Bigl\{\text{\rm var}\bigl(\mathbf{F};[0,S]\bigl)+\text{\rm var}\bigl(\mathbf{T};[0,S]\bigl)\Bigr\}.
\end{align*}
The second estimate of the proposition now follows from estimate~\eqref{eq:tBound}.\qed

\bigskip

\begin{prop}
	\label{prop:weakConvegences}
	There exist subsequences of \(\mathbf{u}_+^m\), \(\mathbf{u}_-^m\), \(t_{+,{\rm n}}^m\), \(t_{-,{\rm n}}^m\) (not relabelled) such that, for all \(s\in [0,S]\):
	\begin{align*}
		& \mathbf{u}_+^m(s) \rightharpoonup \mathbf{u}(s) \quad\text{ weakly in } H^1(\Omega;\mathbb{R}^2),\\ & \mathbf{u}_-^m(s) \rightharpoonup \mathbf{u}(s) \quad\text{ weakly in } H^1(\Omega;\mathbb{R}^2),\\ & t_{+,{\rm n}}^m(s) \rightharpoonup t_{\rm n}(s) \quad\text{ weakly in } H^{-1/2}(\Gamma_C),\\
		& t_{-,{\rm n}}^m(s) \rightharpoonup t_{\rm n}(s) \quad\text{ weakly in } H^{-1/2}(\Gamma_C),
	\end{align*}
	where \(\mathbf{u}\in \textsl{BV}^+([0,S];V)\) and \(t_{\rm n}\in \textsl{BV}^+([0,S];H^{-1/2}(\Gamma_C))\).
\end{prop}

\noindent\textbf{Proof.} Combining Proposition~\ref{prop:uniformEstimatesSweepingProcess} and the Helly selection theorem \cite[Theorem 1.126]{BarBuPrecupanu}, there exist \(\mathbf{u}\in \textsl{BV}([0,S];V)\) and \(t_{\rm n}\in \textsl{BV}([0,S];H^{-1/2}(\Gamma_C))\) such that, for all \(s\in [0,S]\):
\begin{align*}
	& \mathbf{u}_+^m(s) \rightharpoonup \mathbf{u}(s) \quad\text{ weakly in } H^1(\Omega;\mathbb{R}^2),\\ & t_{+,{\rm n}}^m(s) \rightharpoonup t_{\rm n}(s) \quad\text{ weakly in } H^{-1/2}(\Gamma_C).
\end{align*}  
There remains only to prove that \(\mathbf{u}\) and \(t\) are right-continuous.  By the first estimate of Proposition~\ref{prop:uniformEstimatesSweepingProcess}, we have, for all \(s_1\leq s_2\in [0,S]\):
\begin{equation}
	\label{eq:rightContinuityCoulomb}
		\bigl\|\mathbf{u}_+^m(s_2)-\mathbf{u}_+^m(s_1)\bigr\|_{H^1} + \bigl\|t_{+,{\rm n}}^m(s_2)-t_{+,{\rm n}}^m(s_1)\bigr\|_{H^{-1/2}} \leq C\Bigl\{\text{var}\bigl(\mathbf{F};[s_1,s_2]\bigr) + \text{var}\bigl(\mathbf{T};[s_1,s_2]\bigr)\Bigr\}.
\end{equation}
Taking the limit \(m\to\infty\) and using the weak lower semicontinuity of the norms, the same inequality holds for the limits \(\mathbf{u}\), \(t_{\rm n}\). As the right-hand side tends to zero as \(s_2\to s_1+\), we deduce that \(\mathbf{u}(s)\) and \(t_{\rm n}(s)\) are right-continuous, for all \(s\in \left[0,S\right[\). 

Similarly, there exist \(\tilde{\mathbf{u}}\in \textsl{BV}([0,S];V)\) and \(\tilde{t}_{\rm n}\in \textsl{BV}([0,S];H^{-1/2}(\Gamma_C))\) such that, for all \(s\in [0,S]\):
\begin{align*}
	& \mathbf{u}_-^m(s) \rightharpoonup \tilde{\mathbf{u}}(s) \quad\text{ weakly in } H^1(\Omega;\mathbb{R}^2),\\ & t_{-,{\rm n}}^m(s) \rightharpoonup \tilde{t}_{\rm n}(s) \quad\text{ weakly in } H^{-1/2}(\Gamma_C).
\end{align*}  
As estimate~\eqref{eq:rightContinuityCoulomb} also holds for \(\mathbf{u}_-^m\) and \(t_{-,{\rm n}}^m\), we deduce that \(\tilde{\mathbf{u}}(s)\) and \(\tilde{t}_{\rm n}(s)\) are right-continuous, for all \(s\in \left[0,S\right[\). 
Let \(Q\) be the countable set of all the \(s_i^m\) with \(m\in \mathbb{N}\) and \(i\in\{0,1,\ldots,n_m\}\).  We have:
\[
	\forall m\in \mathbb{N},\quad \forall s\in [0,S]\setminus Q,\qquad \mathbf{u}_+^m(s) = \mathbf{u}_-^m(s),
\]
which entails:
\[
	\forall s\in [0,S]\setminus Q,\qquad \mathbf{u}(s) = \tilde{\mathbf{u}}(s).
\]
As \(\mathbf{u}\) and \(\tilde{\mathbf{u}}\) are both right-continuous, and \([0,S]\setminus Q\) is dense in \([0,S]\), we deduce that they coincide on the whole interval \([0,S]\).  The same reasoning applies to \(t_{\rm n}\) and \(\tilde{t}_{\rm n}\). \qed

\bigskip
The following Theorem is the cornerstone in the proof of Theorem~\ref{theo:mainExistenceTheorem}. It relies on our compensated-compactness-like Theorem~\ref{theo:BIcompensatedCompactness} combined with pseudomonotonicity-like arguments.

\begin{theo}
	\label{theo:pseudomonotonicity}
	We suppose that \(N=2\) and \(\boldsymbol{\Lambda}\in W^{1,\infty}(\Omega)\) is \emph{isotropic} at each \(x\in \Omega\). For all \(s\in [0,S]\), the sequence \((t_{+,{\rm n}}^m(s),t_{-,{\rm n}}^m(s))\) converges strongly in \(H^{-1/2}(\Gamma_C)\times H^{-1/2}(\Gamma_C)\) towards \((t_{\rm n}(s),t_{\rm n}(s))\), as \(m\to\infty\).  For any finite positive measure \(\mu\) on \([0,S]\), the function \((t_{+,{\rm n}}^m,t_{-,{\rm n}}^m)\) converges strongly in \(L^2([0,S],\mu;H^{-1/2}(\Gamma_C)\times H^{-1/2}(\Gamma_C))\) towards \((t_{\rm n},t_{\rm n})\).
\end{theo}

\noindent\textbf{Proof.} By definition of \(\mathbf{u}_+^m\), we have for all \(s\in [0,S]\):
\[
t_{+,{\rm n}}^m(s)\leq 0\;\text{ on }\Gamma_C,\quad \text{ and }\quad	\forall t\in H^{-1/2}(\Gamma_C),\text{ such that }\:t\leq 0,\quad \int_{\Gamma_C}\bigl(u_{+,\rm n}^m(s)-g\bigr)\bigl(t-t_{+,{\rm n}}^m\bigr)\geq 0.
\]
Using \(t = t_{\rm n}(s)\), we obtain:
\begin{equation}
	\label{eq:pseudoMonotonicity}
		\limsup_{m\to\infty}\int_{\Gamma_C}u_{+,\rm n}^m(s)\,\bigl(t_{+,{\rm n}}^m(s)-t_{\rm n}(s)\bigr) \leq 0,
\end{equation}
where we have used Proposition~\ref{prop:weakConvegences}.  In addition, for all \(s\in [0,S]\), \(\mathbf{u}_+^m(s)\) satisfies:
\[
\left\{\quad
\begin{array}{ll}
	\text{div}\,\boldsymbol{\sigma}(\mathbf{u}_+^m(s)) + \mathbf{F}_+^m(s) = \mathbf{0}, & \text{ in }\Omega,\\[0.4ex]
		\mathbf{u}_+^m(s) = \mathbf{0}, & \text{ on }\Gamma_U,\\[0.4ex]
		\boldsymbol\sigma(\mathbf{u}_+^m(s))\,\mathbf{n} = \mathbf{T}_+^m, & \text{ on }\Gamma_T,\\[0.4ex]
		\boldsymbol\sigma(\mathbf{u}_+^m(s))\,\mathbf{n}=t_{+,{\rm n}}^m(s)\,\mathbf{n}+t_{+,{\rm t}}^m(s)\,\boldsymbol{\tau}, & \text{ on }\Gamma_C,
\end{array}
\right.
\]
for some measure \(t_{+,{\rm t}}^m(s)\) in \(H^{-1/2}(\Gamma_C))\) such that \(|t_{+,{\rm t}}^m(s)|\leq -ft_{+,{\rm n}}^m(s)\leq -\|f\|_{L^\infty}t_{+,{\rm n}}^m(s)\).  The weak convergence of \(\mathbf{u}_+^m(s)\) towards \(\mathbf{u}(s)\) in \(H^1(\Omega;\mathbb{R}^2)\) (by Proposition~\ref{prop:weakConvegences}), together with the strong convergence in \(L^2\) of \(\mathbf{F}_+^m(s)\) and \(\mathbf{T}_+^m(s)\) (by Proposition~\ref{prop:convData}) yield:
\[
t_{+,{\rm t}}^m(s)\rightharpoonup t_{\rm t}(s) \quad\text{ weakly in } H^{-1/2}(\Gamma_C),
\]
where \(t_{\rm t}(s)\) stands for the tangential part of \(\mathbf{t}(s)=\boldsymbol\sigma(\mathbf{u}(s))\,\mathbf{n}\) on \(\Gamma_C\) (using Formula~\eqref{eq:Stokes}).   We introduce the displacement \(\mathbf{U}^m(s)\) defined by:
\[
\left\{\quad
\begin{array}{ll}
	\text{div}\,\boldsymbol{\sigma}(\mathbf{U}^m(s)) + \mathbf{F}_+^m(s) = \mathbf{0}, & \text{ in }\Omega,\\[0.4ex]
		\mathbf{U}^m(s) = \mathbf{0}, & \text{ on }\Gamma_U,\\[0.4ex]
		\boldsymbol\sigma(\mathbf{U}^m(s)\,\mathbf{n} = \mathbf{T}_+^m(s), & \text{ on }\Gamma_T,\\[0.4ex]
		\boldsymbol\sigma(\mathbf{U}^m(s))\,\mathbf{n}=\mathbf{0}, & \text{ on }\Gamma_C,
\end{array}
\right.
\]
so that the sequence \(\mathbf{U}^m(s)\) converges strongly in \(H^1(\Omega;\mathbb{R}^2)\) towards the solution \(\mathbf{U}(s)\) of the corresponding elastic problem with \((\mathbf{F}_+^m(s),\mathbf{T}_+^m(s))\) replaced with \((\mathbf{F}(s),\mathbf{T}(s))\). Using the notation in Definition~\ref{defi:linearOperators}, we have:
\begin{equation}
	\label{eq:decompositionUn}
		\mathbf{u}_{+,{\rm n}}^m(s) = U_{\rm n}^m(s) + L_{\rm n}\bigl(t_{+,{\rm n}}^m(s)\bigr) + L_{\rm t}\bigl(t_{+,{\rm t}}^m(s)\bigr).
\end{equation}
Our compensated-compactness-like Theorem~\ref{theo:BIcompensatedCompactness} now yields:
\[
	\lim_{m\to\infty}\int_{\Gamma_C}t_{+,{\rm n}}^m(s)\,L_{\rm t}\bigl(t_{+,{\rm t}}^m(s)\bigr) = \int_{\Gamma_C}t_{\rm n}(s)\,L_{\rm t}\bigl(t_{\rm t}(s)\bigr).
\]
Taking this into account together with the weak lower semicontinuity of the norm \(\sqrt{\langle t,L_{\rm n}(t)\rangle}\), the decomposition~\eqref{eq:decompositionUn} of \(\mathbf{u}_{+,{\rm n}}^m(s)\) and inequality~\eqref{eq:pseudoMonotonicity}, we obtain:
\begin{multline*}
	\int_{\Gamma_C} t_{\rm n}(s)\,L_{\rm n}\bigl(t_{\rm n}(s)\bigr) \leq \liminf_{m\to\infty}\int_{\Gamma_C} t_{+,{\rm n}}^m(s)\,L_{\rm n}\bigl(t_{+,{\rm n}}^m(s)\bigr) \leq \\ \limsup_{m\to\infty}\int_{\Gamma_C} t_{+,{\rm n}}^m(s)\,L_{\rm n}\bigl(t_{+,{\rm n}}^m(s)\bigr) \leq \int_{\Gamma_C} t_{\rm n}(s)\,L_{\rm n}\bigl(t_{\rm n}(s)\bigr),
\end{multline*}
which yields:
\[
\lim_{m\to\infty}\int_{\Gamma_C} t_{+,{\rm n}}^m(s)\,L_{\rm n}\bigl(t_{+,{\rm n}}^m(s)\bigr) = \int_{\Gamma_C} t_{\rm n}(s)\,L_{\rm n}\bigl(t_{\rm n}(s)\bigr).
\]
As \(t_{+,{\rm n}}^m(s)\) converges weakly in \(H^{-1/2}(\Gamma_C)\) towards \(t_{\rm n}(s)\) for all \(s\in [0,S]\), this identity yields the  convergence of the norms (by Proposition~\ref{prop:bilinearMapping}), and therefore the strong convergence of \(t_{+,{\rm n}}^m(s)\) in \(H^{-1/2}(\Gamma_C)\). The proof for \(t_{-,{\rm n}}^m(s)\) runs exactly along the same lines. 

As the sequence \((t_{+,{\rm n}}^m(s),t_{-,{\rm n}}^m(s))\) is uniformly (with respect to both \(m\) and \(s\in [0,S]\)) bounded in \(H^{-1/2}(\Gamma_C)\times H^{-1/2}(\Gamma_C)\) by the second estimate in Proposition~\ref{prop:uniformEstimatesSweepingProcess}, pointwise strong convergence combined with the dominated convergence theorem yields the strong convergence of \((t_{+,{\rm n}}^m,t_{-,{\rm n}}^m)\) in \(L^2([0,S],\mu;H^{-1/2}(\Gamma_C)\times H^{-1/2}(\Gamma_C))\), for any finite positive measure \(\mu\) on \([0,S]\). \qed

\bigskip
The convergence properties in Theorem~\ref{prop:weakConvegences} and Theorem~\ref{theo:pseudomonotonicity} obviously entail that \((u_{\rm n},t_{\rm n})\) satisfies the contact conditions in terms of the variational inequality in Definition~\ref{def:quasi-static} (second condition in Definition~\ref{def:quasi-static}). There remains to show that the limit \(\mathbf{u}\) satisfies the evolution inequality in Definition~\ref{def:quasi-static} (first condition in Definition~\ref{def:quasi-static}). Set:
\begin{equation}
	\label{eq:defvm}
		\mathbf{v}^m(s) 	:= \left|\begin{array}{ll}
		\mathbf{0},\quad & \text{if }s=0,\\[1.0ex]
		\mathbf{0},\quad & \text{if }s\in \left]s_{i-1}^m,s_i^m\right]\text{ and  }v(s_i^m)= v(s_{i-1}^m),\\[1.0ex]
		\displaystyle
		\frac{\mathbf{u}_{i}^m-\mathbf{u}_{i-1}^m}{v(s_i^m)- v(s_{i-1}^m)},\quad & \text{if }s\in \left]s_{i-1}^m,s_i^m\right]\text{ and  }v(s_i^m)> v(s_{i-1}^m).
	\end{array}\right.	
\end{equation}

\begin{prop}
	\label{prop:convergenceVelocity}
	We denote by \(\dot{v}\) the distributional derivative of the non-decreasing function \(v\) on \([0,S]\). It is a finite positive measure on \([0,S]\) and the measure \(\dot{\mathbf{u}}\in\mathscr{M}([0,S];H^1)\) is absolutely continuous with respect to \(\dot{v}\).  The sequence \(\mathbf{v}^m\) converges weakly in \(L^2([0,S],\dot{v};H^1)\) towards the Radon-Nikodym derivative \(\tilde{\mathbf{v}}={\rm d}\dot{\mathbf{u}}/{\rm d}\dot{v}\in L^2([0,S],\dot{v};H^1)\).
\end{prop}

\noindent\textbf{Proof.} The first estimate of Proposition~\ref{prop:uniformEstimatesSweepingProcess} yields:
\[
	\forall s_1,s_2\in [0,S],\qquad \int_{\left]s_1,s_2\right[}\|\dot{\mathbf{u}}_+^m\|\leq C\Bigl\{\text{var}\bigl(\mathbf{F};\left]s_1,s_2\right[\bigr)+\text{var}\bigl(\mathbf{T};\left]s_1,s_2\right[\bigr)\Bigr\},
\]
which, by Helly's theorem \cite[Corollary 1.128]{BarBuPrecupanu}, entails that \(\mathbf{u}_+^m\) converges weakly-* in \(\textsl{BV}^+([0,S];V)\) and that:
\[
	\forall s_1,s_2\in [0,S],\qquad \int_{\left]s_1,s_2\right[}\|\dot{\mathbf{u}}\| \leq \liminf_{m\to+\infty} \int_{\left]s_1,s_2\right[}\bigl\|\dot{\mathbf{u}}_+^m\bigr\|\leq C\int_{\left]s_1,s_2\right[}\dot{v},
\]
where \(\|\dot{\mathbf{u}}\|\) stands for the total variation measure of \(\dot{\mathbf{u}}\in\mathscr{M}([0,S];H^1)\) and the function \(v\) was defined in formula~\eqref{eq:defVarLoad}. It is sufficient to conclude that \(\dot{\mathbf{u}}\) is absolutely continuous with respect to \(\dot{v}\), (the inequality applies to open intervals, therefore to any open set and finally to any Borel set by exterior regularity of Borel measures).  As the Radon-Nikodym theorem extends to measures that take values in a Hilbert space (see, e.g., \cite[Section VII.7, p.218]{diestelUhl}), there exists a unique \(\tilde{\mathbf{v}}\in L^1([0,S],\dot{v};H^1)\) such that:
\[
	\dot{\mathbf{u}} = \tilde{\mathbf{v}}\,\dot{v}.
\]
By the first estimate of Proposition~\ref{prop:uniformEstimatesSweepingProcess}, we have the estimate:
\[
	\int_{[0,S]}\bigl\|\mathbf{v}^m(s)\bigr\|_{H^1}^2\,\dot{v} \leq C^2v(S),
\]
where \(\mathbf{v}^m\) was defined in~\eqref{eq:defvm}.  We can therefore extract a subsequence, still denoted by \(\mathbf{v}^m\), which converges weakly in \(L^2([0,S],\dot{v};H^1)\) towards some limit \(\mathbf{l}\in L^2([0,S],\dot{v};H^1)\).  Picking arbitrary \(s\in [0,S]\) and \(\mathbf{w}\in H^1(\Omega;\mathbb{R}^N)\), we have:
\[
	\int_{[0,\sigma^m(s)]}\mathbf{w}\cdot\mathbf{v}^m\,\dot{v} = \mathbf{w}\cdot\bigl(\mathbf{u}_-^m(s)-\mathbf{u}_0\bigr),
\]
where the `\(\cdot\)' stands for the inner product in \(H^1(\Omega;\mathbb{R}^N)\), \(\sigma^m(s):=s_i^m\), for \(s\in \left]s_{i-1}^m,s_i^m\right]\), and \(\sigma^m(0):=0\). To get the above identity, we have used that \(v(s_i^m)=v(s_{i-1}^m)\) implies \(\mathbf{u}_i^m=\mathbf{u}_{i-1}^m\) by Proposition~\ref{prop:uniformEstimatesSweepingProcess}. Taking the limit \(m\to\infty\) in this equality, using Proposition~\ref{theo:pseudomonotonicity}, we obtain:
\[
\forall s\in [0,S],\quad \forall\mathbf{w}\in H^1(\Omega;\mathbb{R}^N),\qquad\int_{[0,s]}\mathbf{w}\cdot\mathbf{l}\,\dot{v} = \mathbf{w}\cdot\bigl(\mathbf{u}(s)-\mathbf{u}_0\bigr),
\]
which shows that \(\mathbf{l}\,\dot{v}=\dot{\mathbf{u}}\), and therefore that \(\mathbf{l}=\tilde{\mathbf{v}}\), \(\dot{v}-a.e\).  As each converging subsequence of \(\mathbf{v}^m\) converges to the same limit \(\tilde{\mathbf{v}}\), we conclude that the whole sequence converges weakly in \(L^2([0,S],\dot{v};H^1)\) towards \(\tilde{\mathbf{v}}\).  \qed

\begin{theo}
	\label{prop:finalConvergenceSweepingProcess}
	We suppose that \(N=2\) and \(\boldsymbol{\Lambda}\in W^{1,\infty}(\Omega)\) is \emph{isotropic} at each \(x\in \Omega\). The function \((\mathbf{u},t_{\rm n})\in \textsl{BV}^+([0,S];V\times H^{-1/2}(\Gamma_C))\) is a solution of the quasi-static problem with data \(\mathbf{F}\), \(\mathbf{T}\) and initial condition \(\mathbf{u}_0\), in the sense of Definition~\ref{def:quasi-static}.
\end{theo}

\noindent\textbf{Proof.} By Proposition~\ref{lem:approxSweeping}, taking \(\mathbf{v}+\mathbf{\dot{u}}_+^m\) as test function, we get:
\[
	\int_{[0,S]}\int_{\Omega} \boldsymbol{\varepsilon}(\mathbf{u}_+^m):\boldsymbol{\Lambda}\boldsymbol{\varepsilon}(\mathbf{v}) \geq \int_{[0,S]}\int_{\Omega} \mathbf{F}_+^m\cdot\mathbf{v} + \int_{[0,S]}\int_{\Gamma_T}\mathbf{T}_+^m\cdot\mathbf{v}  + \int_{[0,S]}\int_{\Gamma_C}t_{+,{\rm n}}^m\,v_{\rm n}+\int_{[0,S]}\int_{\Gamma_C}ft_{+,{\rm n}}^m\bigl|v_{\rm t}\bigr|,
\]
for all \(\mathbf{v}\in\mathscr{M}([0,S];V)\) and all \(m\in\mathbb{N}\).   Considering the Radon-Nikodym derivative \(\bar{\mathbf{v}}\) of \(\mathbf{v}\) with respect to its total variation measure as in Definition~\ref{defi:intMeasure}, the convergence properties of the sequence \(\mathbf{u}_+^m\) in Proposition~\ref{prop:weakConvegences} combined with the uniform (with respect to both \(m\) and \(s\in [0,S]\)) estimate for \(\mathbf{u}_+^m\) (second estimate in Proposition~\ref{prop:uniformEstimatesSweepingProcess}) offers the possibility to pass to the limit \(m\to\infty\) in the left-hand side of the above inequality, by dominated convergence. The same argument applies to all the integrals in the right-hand side, thanks to Propositions~\ref{prop:convData},~\ref{prop:uniformEstimatesSweepingProcess} and \ref{prop:weakConvegences}, yielding:
\begin{equation}
	\label{eq:uisintheConvex}
	\int_{[0,S]}\int_{\Omega} \boldsymbol{\varepsilon}(\mathbf{u}):\boldsymbol{\Lambda}\boldsymbol{\varepsilon}(\mathbf{v}) \geq \int_{[0,S]}\int_{\Omega} \mathbf{F}\cdot\mathbf{v} + \int_{[0,S]}\int_{\Gamma_T}\mathbf{T}\cdot\mathbf{v}  + \int_{[0,S]}\int_{\Gamma_C}t_{{\rm n}}\,v_{\rm n}+\int_{[0,S]}\int_{\Gamma_C}ft_{{\rm n}}\bigl|v_{\rm t}\bigr|.
\end{equation}

Taking successively \(\mathbf{v}=\mathbf{u}_{i-1}^m\) and \(\mathbf{v}=2\mathbf{u}_{i}^m-\mathbf{u}_{i-1}^m\) in the inequality~\eqref{eq:optimalityInequality} for \(\mathbf{u}_i^m\), we obtain:
\begin{multline*}
	\int_{\Omega} \boldsymbol{\varepsilon}(\mathbf{u}_i^m):\boldsymbol{\Lambda}\boldsymbol{\varepsilon}(\mathbf{u}_i^m-\mathbf{u}_{i-1}^m) = \int_{\Omega} \mathbf{F}_i^m\cdot(\mathbf{u}_i^m-\mathbf{u}_{i-1}^m) + \int_{\Gamma_T}\mathbf{T}_i^m\cdot(\mathbf{u}_i^m-\mathbf{u}_{i-1}^m) \\\mbox{} + \int_{\Gamma_C}t_{i,{\rm n}}^m\,(v_{i,{\rm n}}^m-v_{i-1,{\rm n}}^m)+\int_{\Gamma_C}ft_{i,{\rm n}}^m\bigl|u_{i,{\rm t}}^m-u_{i-1,{\rm t}}^m\bigr|,
\end{multline*}
which yields:
\begin{multline}
	\sum_{i=1}^{n_m}\int_{\Omega} \boldsymbol{\varepsilon}(\mathbf{u}_i^m):\boldsymbol{\Lambda}\boldsymbol{\varepsilon}(\mathbf{u}_i^m-\mathbf{u}_{i-1}^m) = \int_{[0,S]}\int_{\Omega} \mathbf{F}_-^m\cdot\mathbf{v}^m\,\dot{v} + \int_{[0,S]}\int_{\Gamma_T}\mathbf{T}_-^m\cdot\mathbf{v}^m\,\dot{v} \\ \mbox{} + \int_{[0,S]}\int_{\Gamma_C}\Bigl(t_{-,{\rm n}}^m\,v_{\rm n}^m+ft_{-,{\rm n}}^m\bigl|v_{\rm t}^m\bigr|\Bigr)\,\dot{v}. \label{eq:termDotU}
\end{multline}
We are now going to analyze the limit \(m\to\infty\) of each side of~\eqref{eq:termDotU}. Let us first focus on the right-hand side. As \(t_{-,{\rm n}}^m\leq 0\), given an arbitrary function \(\varphi\in C^1(\Gamma_C;[-1,1])\), the right-hand side of~\eqref{eq:termDotU} is smaller than:
\[
\int_{[0,S]}\int_{\Omega} \mathbf{F}_-^m\cdot\mathbf{v}^m\,\dot{v} + \int_{[0,S]}\int_{\Gamma_T}\mathbf{T}_-^m\cdot\mathbf{v}^m\,\dot{v} \\ \mbox{} + \int_{[0,S]}\int_{\Gamma_C}\Bigl(t_{-,{\rm n}}^m\,v_{\rm n}^m+ft_{-,{\rm n}}^m\varphi v_{\rm t}^m\Bigr)\,\dot{v},
\]
which can be rewritten as a term of the form:
\[
	\int_{[0,S]}\int_{\Omega}\boldsymbol{\varepsilon}(\mathbf{w}_\varphi^m):\boldsymbol{\Lambda}\boldsymbol{\varepsilon}(\mathbf{v}^m)\,\dot{v},
\]
for some uniquely defined piecewise constant in time sequence \(\mathbf{w}_\varphi^m\) in \(L^2([0,S],\dot{v};V)\), which is a Cauchy sequence in \(L^2([0,S],\dot{v};V)\), thanks to the convergence properties of the sequences \(\mathbf{F}_-^m\), \(\mathbf{T}_-^m\) and \(t_{-,{\rm n}}^m\) provided by Proposition~\ref{prop:convData} and Theorem~\ref{theo:pseudomonotonicity}. Its strong limit \(\mathbf{w}_\varphi\) solves the elastic problem with data \((\mathbf{F},\mathbf{T}, t_{\rm n}, \varphi t_{\rm n})\). Hence, by the weak convergence of \(\mathbf{v}^m\) in \(L^2([0,S],\dot{v};V)\) obtained in Proposition~\ref{prop:convergenceVelocity}, we have, for all \(\varphi\in C^1(\Gamma_C;[-1,1])\):
\begin{multline*}
	\limsup_{m\to+\infty}\Biggl\{\int_{[0,S]}\int_{\Omega} \mathbf{F}_-^m\cdot\mathbf{v}^m\,\dot{v} + \int_{[0,S]}\int_{\Gamma_T}\mathbf{T}_-^m\cdot\mathbf{v}^m\,\dot{v} + \int_{[0,S]}\int_{\Gamma_C}\Bigl(t_{-,{\rm n}}^m\,v_{\rm n}^m+ft_{-,{\rm n}}^m\bigl| v_{\rm t}^m\bigr|\Bigr)\,\dot{v}\Biggr\} \leq  \mbox{} \\ \int_{[0,S]}\int_{\Omega} \mathbf{F}\cdot\dot{\mathbf{u}} + \int_{[0,S]}\int_{\Gamma_T}\mathbf{T}\cdot\dot{\mathbf{u}} + \int_{[0,S]}\int_{\Gamma_C}t_{\rm n} \dot{u}_{\rm n}+\int_{[0,S]}\int_{\Gamma_C}ft_{\rm n}\varphi\dot{u}_{\rm t}.
\end{multline*}
As \(\varphi\in C^1(\Gamma_C;[-1,1])\) is arbitrary, we have proved that:
\begin{multline}
	\limsup_{m\to+\infty}\Biggl\{\int_{[0,S]}\int_{\Omega} \mathbf{F}_-^m\cdot\mathbf{v}^m\,\dot{v} + \int_{[0,S]}\int_{\Gamma_T}\mathbf{T}_-^m\cdot\mathbf{v}^m\,\dot{v} + \int_{[0,S]}\int_{\Gamma_C}\Bigl(t_{-,{\rm n}}^m\,v_{\rm n}^m+ft_{-,{\rm n}}^m\bigl| v_{\rm t}^m\bigr|\Bigr)\,\dot{v}\Biggr\} \leq  \mbox{} \\ \int_{[0,S]}\int_{\Omega} \mathbf{F}\cdot\dot{\mathbf{u}} + \int_{[0,S]}\int_{\Gamma_T}\mathbf{T}\cdot\dot{\mathbf{u}} + \int_{[0,S]}\int_{\Gamma_C}t_{\rm n} \dot{u}_{\rm n}+\int_{[0,S]}\int_{\Gamma_C}ft_{\rm n}\bigl|\dot{u}_{\rm t}\bigr|.\label{eq:rightHandSide}
\end{multline}
Tackling the limit \(m\to\infty\) in the left-hand side of~\eqref{eq:termDotU} is going to be more involved. First note that:
\begin{multline*}
	\int_{\Omega} \boldsymbol{\varepsilon}(\mathbf{u}_i^m):\boldsymbol{\Lambda}\boldsymbol{\varepsilon}(\mathbf{u}_i^m-\mathbf{u}_{i-1}^m) = \frac{1}{2}\int_{\Omega} \boldsymbol{\varepsilon}(\mathbf{u}_i^m):\boldsymbol{\Lambda}\boldsymbol{\varepsilon}(\mathbf{u}_i^m) - \frac{1}{2}\int_{\Omega} \boldsymbol{\varepsilon}(\mathbf{u}_{i-1}^m):\boldsymbol{\Lambda}\boldsymbol{\varepsilon}(\mathbf{u}_{i-1}^m) \\ \mbox{}+ \frac{1}{2}\int_{\Omega} \boldsymbol{\varepsilon}(\mathbf{u}_i^m-\mathbf{u}_{i-1}^m):\boldsymbol{\Lambda}\boldsymbol{\varepsilon}(\mathbf{u}_i^m-\mathbf{u}_{i-1}^m),
\end{multline*}
so that:
\begin{multline*}
	\sum_{i=1}^{n_m}\int_{\Omega} \boldsymbol{\varepsilon}(\mathbf{u}_i^m):\boldsymbol{\Lambda}\boldsymbol{\varepsilon}(\mathbf{u}_i^m-\mathbf{u}_{i-1}^m) = \frac{1}{2}\int_{\Omega} \boldsymbol{\varepsilon}\bigl(\mathbf{u}_+^m(S)\bigr):\boldsymbol{\Lambda}\boldsymbol{\varepsilon}\bigl(\mathbf{u}_+^m(S)\bigr)  - \frac{1}{2}\int_{\Omega} \boldsymbol{\varepsilon}(\mathbf{u}_0):\boldsymbol{\Lambda}\boldsymbol{\varepsilon}(\mathbf{u}_0)\\ \mbox{} + \frac{1}{2}\sum_{i=1}^{n_m}\bigl(v(s_i^m)-v(s_{i-1}^m)\bigr)\int_{\left]s_{i-1}^m,s_i^m\right]}\int_{\Omega} \boldsymbol{\varepsilon}(\mathbf{v}^m):\boldsymbol{\Lambda}\boldsymbol{\varepsilon}(\mathbf{v}^m)\,\dot{v},
\end{multline*}
which can be rewritten as:
\begin{multline*}
	\sum_{i=1}^{n_m}\int_{\Omega} \boldsymbol{\varepsilon}(\mathbf{u}_i^m):\boldsymbol{\Lambda}\boldsymbol{\varepsilon}(\mathbf{u}_i^m-\mathbf{u}_{i-1}^m) = \frac{1}{2}\int_{\Omega} \boldsymbol{\varepsilon}\bigl(\mathbf{u}_+^m(S)\bigr):\boldsymbol{\Lambda}\boldsymbol{\varepsilon}\bigl(\mathbf{u}_+^m(S)\bigr)  - \frac{1}{2}\int_{\Omega} \boldsymbol{\varepsilon}(\mathbf{u}_0):\boldsymbol{\Lambda}\boldsymbol{\varepsilon}(\mathbf{u}_0)\\ \mbox{} + \frac{1}{2}\int_{\left[0,S\right]}\int_{\Omega} \boldsymbol{\varepsilon}(\mathbf{v}^m):\boldsymbol{\Lambda}\boldsymbol{\varepsilon}(\mathbf{v}^m)\;\biggl\{\sum_{i=1}^{n_m}\bigl(v(s_i^m)-v(s_{i-1}^m)\bigr)\chi_{\left]s_{i-1}^m,s_i^m\right]}\biggr\}\,\dot{v},
\end{multline*}
where \(\chi_{\left]s_{i-1}^m,s_i^m\right]}\) stands for the characteristic function of the interval \(\left]s_{i-1}^m,s_i^m\right]\). We denote by \(A\) the (at most countable) set of discontinuity points of the function \(v\) on \([0,S]\), which is also the set of atoms of the measure \(\dot{v}\). For all \(s\in\left]s_{i-1}^m,s_i^m\right]\), we have:
\[
	\dot{v}\bigl(\{s\}\bigr) = v(s+)-v(s-) \leq v(s_i^m)-v(s_{i-1}^m),
\]
so that:
\[
	\forall s\in [0,S],\qquad \sum_{i=1}^{n_m}\bigl(v(s_i^m)-v(s_{i-1}^m)\bigr)\chi_{\left]s_{i-1}^m,s_i^m\right]}(s) \geq \dot{v}\bigl(\{s\}\bigr)\chi_A(s).
\]
Finally, we have proved that:
\begin{multline}
	\sum_{i=1}^{n_m}\int_{\Omega} \boldsymbol{\varepsilon}(\mathbf{u}_i^m):\boldsymbol{\Lambda}\boldsymbol{\varepsilon}(\mathbf{u}_i^m-\mathbf{u}_{i-1}^m) \geq \frac{1}{2}\int_{\Omega} \boldsymbol{\varepsilon}\bigl(\mathbf{u}_+^m(S)\bigr):\boldsymbol{\Lambda}\boldsymbol{\varepsilon}\bigl(\mathbf{u}_+^m(S)\bigr)  - \frac{1}{2}\int_{\Omega} \boldsymbol{\varepsilon}(\mathbf{u}_0):\boldsymbol{\Lambda}\boldsymbol{\varepsilon}(\mathbf{u}_0)\\ \mbox{} + \frac{1}{2}\int_{A}\int_{\Omega} \boldsymbol{\varepsilon}\bigl(\mathbf{v}^m(s)\bigr):\boldsymbol{\Lambda}\boldsymbol{\varepsilon}\bigl(\mathbf{v}^m(s)\bigr)\;\dot{v}\bigl(\{s\}\bigr)\,\dot{v}.
	\label{eq:ineqDeLaMort}
\end{multline}
As \(\mathbf{v}^m\) converges weakly in \(L^2([0,S],\dot{v};H^1)\) towards the Radon-Nikodym derivative \(\tilde{\mathbf{v}}={\rm d}\dot{\mathbf{u}}/{\rm d}\dot{v}\in L^2([0,S],\dot{v};H^1)\), the sequence \(\mathbf{v}^m(s)\,\sqrt{\dot{v}\bigl(\{s\}\bigr)}\,\chi_A(s)\) also converges weakly in \(L^2([0,S],\dot{v};H^1)\) towards \(\tilde{\mathbf{v}}(s)\,\sqrt{\dot{v}\bigl(\{s\}\bigr)}\,\chi_A(s)\).  Using, the weak lower semicontinuity of the norm, we have:
\begin{multline*}
\liminf_{m\to +\infty}\int_{A}\int_{\Omega} \boldsymbol{\varepsilon}\bigl(\mathbf{v}^m(s)\bigr):\boldsymbol{\Lambda}\boldsymbol{\varepsilon}\bigl(\mathbf{v}^m(s)\bigr)\;\dot{v}\bigl(\{s\}\bigr)\,\dot{v}\geq \int_{A}\int_{\Omega} \boldsymbol{\varepsilon}\bigl(\tilde{\mathbf{v}}(s)\bigr):\boldsymbol{\Lambda}\boldsymbol{\varepsilon}\bigl(\tilde{\mathbf{v}}(s)\bigr)\;\dot{v}\bigl(\{s\}\bigr)\,\dot{v} \\ \mbox{} = \sum_{s\in A}\int_{\Omega} \boldsymbol{\varepsilon}\bigl(\mathbf{u}(s+)-\mathbf{u}(s-)\bigr):\boldsymbol{\Lambda}\boldsymbol{\varepsilon}\bigl(\mathbf{u}(s+)-\mathbf{u}(s-)\bigr).
\end{multline*}
Recalling that \(\mathbf{u}_+^m(S)\) converges weakly in \(H^1(\Omega;\mathbb{R}^2)\) towards \(\mathbf{u}(S)\) by Proposition~\ref{prop:weakConvegences}, we can pass to the limit \(m\to+\infty\) in inequality~\eqref{eq:ineqDeLaMort} and obtain:
\begin{multline*}
\liminf_{m\to +\infty}\;\sum_{i=1}^{n_m}\int_{\Omega} \boldsymbol{\varepsilon}(\mathbf{u}_i^m):\boldsymbol{\Lambda}\boldsymbol{\varepsilon}(\mathbf{u}_i^m-\mathbf{u}_{i-1}^m) \geq \frac{1}{2}\int_{\Omega} \boldsymbol{\varepsilon}\bigl(\mathbf{u}(S)\bigr):\boldsymbol{\Lambda}\boldsymbol{\varepsilon}\bigl(\mathbf{u}(S)\bigr)  - \frac{1}{2}\int_{\Omega} \boldsymbol{\varepsilon}(\mathbf{u}_0):\boldsymbol{\Lambda}\boldsymbol{\varepsilon}(\mathbf{u}_0)\\ \mbox{} + \frac{1}{2}\sum_{s\in A}\int_{\Omega} \boldsymbol{\varepsilon}\bigl(\mathbf{u}(s+)-\mathbf{u}(s-)\bigr):\boldsymbol{\Lambda}\boldsymbol{\varepsilon}\bigl(\mathbf{u}(s+)-\mathbf{u}(s-)\bigr).
\end{multline*}
Recalling that \(\mathbf{u}\) is right continuous, and invoking the identity in formula~\eqref{eq:bvDissipation}, the right-hand side of the above inequality can be rewritten as:
\begin{align*}
	\liminf_{m\to +\infty}\;\sum_{i=1}^{n_m}\int_{\Omega} \boldsymbol{\varepsilon}(\mathbf{u}_i^m):\boldsymbol{\Lambda}\boldsymbol{\varepsilon}(\mathbf{u}_i^m-\mathbf{u}_{i-1}^m) & \geq \frac{1}{2}\int_{\left[0,S\right]}\int_{\Omega} \boldsymbol{\varepsilon}\bigl(\mathbf{u}(s+)+\mathbf{u}(s-)\bigr):\boldsymbol{\Lambda}\boldsymbol{\varepsilon}(\dot{\mathbf{u}})\\
	& \qquad \mbox{}+\frac{1}{2}\int_{\left[0,S\right]}\int_{\Omega} \boldsymbol{\varepsilon}\bigl(\mathbf{u}(s+)-\mathbf{u}(s-)\bigr):\boldsymbol{\Lambda}\boldsymbol{\varepsilon}(\dot{\mathbf{u}}),\\
	& = \int_{\left[0,S\right]}\int_{\Omega} \boldsymbol{\varepsilon}\bigl(\mathbf{u}\bigr):\boldsymbol{\Lambda}\boldsymbol{\varepsilon}\bigl(\dot{\mathbf{u}}\bigr) .
\end{align*}
Gathering this last inequality with inequality~\eqref{eq:rightHandSide}, we can take the inferior limit \(m\to+\infty\) in identity~\eqref{eq:termDotU} and get:
\begin{equation}
	\label{eq:finalInequality}
	\int_{\left[0,S\right]}\int_{\Omega} \boldsymbol{\varepsilon}\bigl(\mathbf{u}\bigr):\boldsymbol{\Lambda}\boldsymbol{\varepsilon}\bigl(\dot{\mathbf{u}}\bigr) \leq \int_{[0,S]}\int_{\Omega} \mathbf{F}\cdot\dot{\mathbf{u}} + \int_{[0,S]}\int_{\Gamma_T}\mathbf{T}\cdot\dot{\mathbf{u}} + \int_{[0,S]}\int_{\Gamma_C}t_{\rm n} \dot{u}_{\rm n}+\int_{[0,S]}\int_{\Gamma_C}ft_{\rm n}\bigl|\dot{u}_{\rm t}\bigr|.
\end{equation}
Combining this inequality with~\eqref{eq:uisintheConvex}, we have proved that \(\mathbf{u}\) fulfills the evolution inequality in Definition~\ref{def:quasi-static} (first condition in Definition~\ref{def:quasi-static}). 

The convergence properties in Theorem~\ref{prop:weakConvegences} and Theorem~\ref{theo:pseudomonotonicity} obviously 
enable to pass to the limit \(m\to+\infty\) in the contact conditions (second condition in Proposition~\ref{lem:approxSweeping}), showing that \((u_{\rm n},t_{\rm n})\) satisfies the second condition in Definition~\ref{def:quasi-static}.\qed

\bigskip
\noindent{\textbf{Remark.}} Taking successively \(\mathbf{v}=\mathbf{0}\) and \(\mathbf{v}=2\dot{\mathbf{u}}\) in the first condition in Definition~\ref{def:quasi-static}, we see that inequality~\eqref{eq:finalInequality} is actually an equality. Rewinding the proof of that inequality, this entails that we actually have:
\[
\lim_{m\to +\infty}\int_{\Omega} \boldsymbol{\varepsilon}\bigl(\mathbf{u}_+^m(S)\bigr):\boldsymbol{\Lambda}\boldsymbol{\varepsilon}\bigl(\mathbf{u}_+^m(S)\bigr) = \int_{\Omega} \boldsymbol{\varepsilon}\bigl(\mathbf{u}(S)\bigr):\boldsymbol{\Lambda}\boldsymbol{\varepsilon}\bigl(\mathbf{u}(S)\bigr),
\]
which shows that the weak convergence in \(H^1(\Omega;\mathbb{R}^2)\) of \(\mathbf{u}_+^m(S)\) towards \(\mathbf{u}(S)\) proved in Proposition~\ref{prop:weakConvegences}, is actually strong, and the same is of course true for all \(\mathbf{u}_+^m(s)\) with \(s\in [0,S]\).

\begin{prop}
	\label{prop:absolutelyContinuousData}
	In the case where the data \(\mathbf{F}\) and \(\mathbf{T}\) are absolutely continuous, then so is the above constructed solution \(\mathbf{u}\) of the evolution inequality in Definition~\ref{def:quasi-static} (first condition in Definition~\ref{def:quasi-static}).  In this case, \(\mathbf{u}\) fulfills the equivalent evolution inequality, for all \(\mathbf{v}\in V\),  and almost all \(s\in [0,S]\):
	\begin{multline*}
		\int_{\Omega} \boldsymbol{\varepsilon}(\mathbf{u}(s)):\boldsymbol{\Lambda}\boldsymbol{\varepsilon}\bigl(\mathbf{v}-\mathbf{\dot{u}}(s)\bigr) \geq \int_{\Omega} \mathbf{F}(s)\cdot\bigl(\mathbf{v}-\mathbf{\dot{u}}(s)\bigr) + \int_{\Gamma_T}\mathbf{T}(s)\cdot\bigl(\mathbf{v}-\mathbf{\dot{u}}(s)\bigr) +\mbox{}\\\mbox{}+ \int_{\Gamma_C}t_{\rm n}(s)\,\Bigl[v_{\rm n}-\dot{u}_{\rm n}(s)\Bigr] +\int_{\Gamma_C} ft_{\rm n}(s)\,\Bigl[\bigl|\mathbf{v}_{\rm t}\bigr|-|\mathbf{\dot{u}}_{\rm t}(s)\bigr|\Bigr].
	\end{multline*}
\end{prop}

\noindent\textbf{Proof.} The first estimate of Proposition~\ref{prop:uniformEstimatesSweepingProcess} yields:
\[
	\forall s_1,s_2\in [0,S],\qquad \int_{[s_1,s_2]}\|\dot{\mathbf{u}}_+^m\| \leq C\Bigl\{\text{var}\bigl(\mathbf{F};[s_1,s_2]\bigr)+\text{var}\bigl(\mathbf{T};[s_1,s_2]\bigr)\Bigr\},
\]
so that the same inequality applies to the total variation measure \(\|\dot{\mathbf{u}}\|\) of \(\dot{\mathbf{u}}\). Therefore, absolute continuity of the data entails that of the solution \(\mathbf{u}\). Picking arbitrary \(s_0\in\left]0,S\right[\), \(v\in V\), taking a test function of the form \(\mathbf{\dot{u}}(s) + (\mathbf{v}-\mathbf{\dot{u}}(s))m\chi_{[s_0-1/m,s_0+1/m]}(s)\) (where \(\chi_{[a,b]}\) stands for the characteristic function of the interval \([a,b]\)) in the evolution inequality in Definition~\ref{def:quasi-static}, and taking the limit \(m\to+\infty\), we obtain the claimed inequality for every Lebesgue point \(s_0\in [0,S]\) of \(\dot{\mathbf{u}}\in L^1(0,S;H^1)\). \qed

\bigskip
Theorem~\ref{prop:finalConvergenceSweepingProcess} and Proposition~\ref{prop:absolutelyContinuousData} yield a proof of Theorem~\ref{theo:mainExistenceTheorem}, \emph{in the isotropic case}.  The extension to the anisotropic case under condition~\eqref{eq:ConditionFrictionAnisotropic} can be done exactly along the same lines, using the same trick as in Section~3.2 of \cite{BI-ma3as} and is left to the reader.

\end{document}